\documentclass[letterpaper,12pt,reqno]{amsart}
\RequirePackage[utf8]{inputenc}
\usepackage[portrait,margin=2.93cm]{geometry}
\usepackage{mathrsfs,soul}
\usepackage{hyperref}
\usepackage[foot]{amsaddr}
\usepackage{amsmath,amssymb,amsthm,amsfonts,amsbsy,latexsym,dsfont,mathtools,nccmath,scalerel}

\usepackage{relsize}

\makeatletter 
\@mparswitchfalse%
\makeatother
\normalmarginpar 
\usepackage{graphicx}
\usepackage[numeric,initials,nobysame]{amsrefs}
\usepackage{upref,setspace}
\usepackage{xcolor,colortbl}

\usepackage{fourier}

\usepackage{xcolor}

\usepackage{mathrsfs,soul}

\usepackage{enumerate}
\newenvironment{enumeratei}{\begin{enumerate}[\upshape (i)]}{\end{enumerate}}
\newenvironment{enumeratea}{\begin{enumerate}[\upshape (a)]}{\end{enumerate}}

\usepackage{paralist}

\newenvironment{inparaenuma}{\begin{inparaenum}[\upshape \bfseries (a) ]}{\end{inparaenum}}

\newenvironment{inparaenumi}{\begin{inparaenum}[\upshape  (i) ]}{\end{inparaenum}}

\definecolor{refkey}{gray}{.75}
\definecolor{labelkey}{gray}{.75}

\newtheorem{thm}{Theorem}[section]
\newtheorem{lem}[thm]{Lemma}

\newtheorem{prop}[thm]{Proposition}

\newtheorem{ass}[thm]{Condition}

\newtheorem{conj}[thm]{Open Problem}

\theoremstyle{remark}

\theoremstyle{definition}
\newtheorem{rem}{Remark}[section]

\newcommand{\ind}{\mathbf{1}}
\newcommand{\eps}{\varepsilon}

\newcommand{\set}[1]{\left\{#1\right\}}

\newcommand{\equald}{\stackrel{\mathrm{d}}{=}}
\newcommand{\probc}{\stackrel{\mathrm{P}}{\longrightarrow}}

\newcommand{\weakc}{\stackrel{\mathrm{d}}{\longrightarrow}}

\def\qed{ \hfill $\blacksquare$}

\newcommand{\cA}{\mathcal{A}}\newcommand{\cC}{\mathcal{C}}
\newcommand{\cD}{\mathcal{D}}\newcommand{\cE}{\mathcal{E}}\newcommand{\cF}{\mathcal{F}}
\newcommand{\cG}{\mathcal{G}}\newcommand{\cH}{\mathcal{H}}\newcommand{\cI}{\mathcal{I}}
\newcommand{\cL}{\mathcal{L}}

\newcommand{\cP}{\mathcal{P}}
\newcommand{\cS}{\mathcal{S}}\newcommand{\cT}{\mathcal{T}}
\newcommand{\cW}{\mathcal{W}}

\newcommand{\vC}{\mathbf{C}}
\newcommand{\vF}{\mathbf{F}}

\newcommand{\vM}{\mathbf{M}}

\newcommand{\vX}{\mathbf{X}}

\newcommand{\ve}{\mathbf{e}}

\newcommand{\vt}{\mathbf{t}}
\newcommand{\vw}{\mathbf{w}}\newcommand{\vx}{\mathbf{x}}

\newcommand{\mva}{\boldsymbol{a}}

\newcommand{\mvy}{\boldsymbol{y}}

\newcommand{\mvxi}{\boldsymbol{\xi}}\newcommand{\mvXi}{\boldsymbol{\Xi}}

\newcommand{\fT}{\mathfrak{T}}
\newcommand{\fX}{\mathfrak{X}}

\newcommand{\fp}{\mathfrak{p}}

\newcommand{\ft}{\mathfrak{t}}

\newcommand{\bE}{\mathbb{E}}
\newcommand{\bH}{\mathbb{H}}

\newcommand{\bN}{\mathbb{N}}
\newcommand{\bR}{\mathbb{R}}
\newcommand{\bT}{\mathbb{T}}

\newcommand{\bZ}{\mathbb{Z}}

\newcommand{\sS}{\mathscr{S}}

 \DeclareMathOperator{\height}{ht}

\DeclareMathOperator{\E}{\mathbb{E}}
\DeclareMathOperator{\pr}{\mathbb{P}}

\DeclareMathOperator{\var}{Var}
\DeclareMathOperator{\cov}{Cov}

\DeclareMathOperator{\GHP}{GHP}

\DeclareMathOperator{\ERRG}{ERRG}

\DeclareMathOperator{\scl}{scl}
\DeclareMathOperator{\er}{er}
\DeclareMathOperator{\mass}{mass}

\DeclareMathOperator{\spls}{spls}

\DeclareMathOperator{\pa}{pa}

\DeclareMathOperator{\crit}{Crit}

\DeclareMathOperator{\pois}{Pois}

\newcommand{\sss}{\scriptscriptstyle}
\newcommand{\erdos}{Erd\H{o}s-R\'enyi }
\newcommand{\ldown}{l^2_{\downarrow}}

\newcommand{\diam}{\mathrm{diam}}

\newcommand{\vCrit}{\mathbf{Crit}}

\newcommand{\bars}{\bar{s}}

\newcommand{\diamax}{{\diam_{\max}}}

\definecolor{jrnl}{rgb}{0.0, 0.5, 0.0}
\definecolor{jrnl1}{rgb}{0.0, 0.72, 0.6}
\definecolor{aqua}{rgb}{0.0, 1.0, 1.0}
\definecolor{webbrown}{rgb}{.6,0,0}
\definecolor{pinegreen}{rgb}{0.0, 0.47, 0.44}
\definecolor{ultramarineblue}{rgb}{0.25, 0.4, 0.96}
\definecolor{lincolngreen}{rgb}{0.11, 0.35, 0.02}
\definecolor{green(html/cssgreen)}{rgb}{0.0, 0.5, 0.0}
\definecolor{airforceblue}{rgb}{0.36, 0.54, 0.66}
\definecolor{azure}{rgb}{0.0, 0.5, 1.0}
\definecolor{bleudefrance}{rgb}{0.19, 0.55, 0.91}
\definecolor{cobalt}{rgb}{0.0, 0.28, 0.67}

\hypersetup{colorlinks=true, linktocpage=true, pdfstartpage=1, pdfstartview=FitV,
	breaklinks=true, pdfpagemode=UseNone, pageanchor=true, pdfpagemode=UseOutlines,
	plainpages=false, bookmarksnumbered, bookmarksopen=true, bookmarksopenlevel=1,
	hypertexnames=true, pdfhighlight=/O,
	urlcolor=webbrown, linkcolor=cobalt, citecolor=jrnl
}

\usepackage[textsize=small]{todonotes}       

\newcommand{\chh}[1]{{#1}}

\def\XXint#1#2#3{{\setbox0=\hbox{$#1{#2#3}{\int}$ }
		\vcenter{\hbox{$#2#3$ }}\kern-.6\wd0}}

\makeatletter
\providecommand{\leftsquigarrow}{%
	\mathrel{\mathpalette\reflect@squig\relax}%
}
\newcommand{\reflect@squig}[2]{%
	\reflectbox{$\m@th#1\rightsquigarrow$}%
}
\makeatother

\DeclareFontFamily{OT1}{pzc}{}
\DeclareFontShape{OT1}{pzc}{m}{it}{<-> s * [1.10] pzcmi7t}{}
\DeclareMathAlphabet{\mathpzc}{OT1}{pzc}{m}{it}

\newcommand{\pmtr}{{\mathpzc{Pmtr}}}

\usetikzlibrary{arrows.meta}

\colorlet{permBlue}{blue!65!black}

\tikzset{
 vtx/.style      = {circle,shade,ball color=gray!35!white,draw=black,
                     very thick,minimum size=7pt,inner sep=0pt},
  permedge/.style = {ultra thick,line cap=round,line join=round,color=permBlue},
  cand/.style     = {dashed,ultra thick,line cap=round,
                     dash pattern=on 6pt off 3pt},
  accept/.style   = {ultra thick,line cap=round,line join=round},
}

\newcommand{\HL}{G^{\sss \text{HL}}}
\newcommand{\torus}{G^{\sss \mathbb{T}}}

\usepackage{tikz}
\usetikzlibrary{positioning, trees}

\usepackage{filemod}

\newcommand{\disjt}{
	\raisebox{-0.9pt}{\scalebox{1.3}{\hskip0.6pt$\circ$\hskip0.5pt}}
}

\newcommand{\annlarge}{
		\raisebox{0.8pt}{\scalebox{.7}{\boxed{\square}}}
}

\newcommand{\annsmall}{
	\raisebox{0.8pt}{\scalebox{.5}{\boxed{\square}}}
}

\newcommand{\edge}{\circ\hskip-3.5pt -\hskip-3.2pt\circ}

\newcommand{\connects}{\leftrightarrow}

\begin{document}

\title[Percolation on the hierarchical lattice and the discrete torus]{Geometry of critical discrete structures: long-range percolation on the hierarchical lattice and the discrete torus}
	
\date{}

\author[Chatterjee]{Arghyadeep Chatterjee}
\author[Maniyar]{Sourish Maniyar}
\author[Sen]{Sanchayan Sen}
\address{Department of Mathematics, Indian Institute of Science}
\email{sanchayan.sen1@gmail.com}
\begin{abstract}
Consider 
(a) balls $\Lambda_n$ of growing volumes in the $d$-dimensional hierarchical lattice, and
(b) the $d$-dimensional discrete torus $\bT_n^d$ on $n^d$ vertices.
Place edges independently between each pair of vertices 
$x\neq y\in\Lambda_n$ or $\bT_n^d$ with probability
$1-\exp(-\beta J(x, y) )$ where 
$J(x, y) \asymp  \| x-y \|^{-\alpha}$ for some $0<\alpha<d$.
For both of these models, we prove the following:
(i)
We obtain tight bounds, up to constants, on the two-point function in the barely subcritical regime.
We show that in part of the barely subcritical regime, the two-point function has a plateau \cite{hutchcroft-michta-slade-torus-plateau, slade-plateau}.
(ii)
We identify the critical window when $0<\alpha<5d/6$. 
Further, using the bound on the two-point function mentioned in (i) together with a universality principle proven in \cites{sbssxw-graphon, sbssxw-universal}, we establish the scaling limit of the maximal components, viewed as metric measure spaces, within the critical window. 
More precisely, we show that the metric scaling limit of the maximal components is Brownian, and that these models belong to the \erdos universality class when $0<\alpha<5d/6$.
It was recently conjectured by Hutchcroft \cite{hutchcroft-critical-cluster-volume}*{Section~7.1} that the model of critical hierarchical percolation with $\alpha\in(d, 4d/3]$ is a member of the \erdos universality class, and we believe that this is also true for all $\alpha\in (0, d]$.
Similarly, critical long-range percolation on the discrete torus is expected to be in this universality class when the effective dimension is high enough. 
These results take a first step in that direction.
(iii)
We show that when $0<\alpha<2d/3$, the girth of each maximal component in the critical window is $\Omega_P(|\Lambda_n|^{1/3})$ and $\Omega_P(n^{d/3})$ respectively for these two models, contrary to the situation when $d<\alpha$ where the girth would equal $3$ .
\end{abstract}
\maketitle


\section{Introduction}\label{sec:intro}

We study the nature of the phase transition in long-range percolation, driven by a large family of kernels, on the hierarchical lattice and the discrete torus.
Our focus will be on understanding the intrinsic geometry inside the critical window, and the closely related problem of getting tight bounds on the two-point function in the barely subcritical regime.
In particular, for a certain class of underlying kernels, we verify the membership of these models in the \erdos universality class.
Let us first describe this universality class.

In an influential paper \cite{aldous1997brownian}, Aldous established the scaling limit of the sequence of rescaled component sizes of the \erdos random graph inside the critical window, and demonstrated the connections between random graph dynamics near criticality and a stochastic process called the multiplicative coalescent.
Then, in \cite{addario2012continuum}, the scaling limit of the sequence of components of the critical \erdos random graph was obtained after incorporating the intrinsic geometry of the components by viewing them as random metric spaces.
This was subsequently strengthened to weak convergence of the components viewed as random metric measure spaces in \cite{addario2013scaling}.
More precisely, this result from \cites{addario2012continuum, addario2013scaling} states the following:
Fix $\lambda\in\bR$, and consider the \erdos random graph $\ERRG(n, n^{-1} + \lambda n^{-4/3})$ on $n$ vertices where edges are placed independently between each pair of vertices with probability $n^{-1} + \lambda n^{-4/3}$.
For each $i\geq 1$, consider the $i$-th maximal component of this graph, and view it as a random metric measure space by endowing it with $n^{-1/3}$ times the graph distance, and $n^{-2/3}$ times the counting measure; denote the resulting metric measure space by $\bar\cC_i^{\sss\er}(\lambda)$.
Then 
\begin{equation}\label{eqn:lim-add-br-go}
	\big(
	\bar\cC_i^{\sss\er}(\lambda) ;\, i\geq 1
	\big)
	\weakc 
	\vCrit(\lambda) := \big( \crit_i(\lambda);\, i\geq 1 \big) 
	\ \text{ as }\ n\to\infty\, ,
\end{equation}
for a sequence of limiting random fractals $\crit_i(\lambda)$ that are described in more detail in Section~\ref{sec:1}; 
the topology of the convergence in \eqref{eqn:lim-add-br-go} will be specified in Section~\ref{sec:topologies}.
The limiting sequence of spaces $\vCrit(\cdot)$ is a fundamental object that is expected to be the universal scaling limit for a host of models of random discrete structures at criticality.
Further, establishing a convergence result as in \eqref{eqn:lim-add-br-go} is intertwined with the analysis of several other probabilistic objects. 
An important example of such an object is the minimal spanning tree (MST): 
As shown in \cite{addario2013scaling}, establishing the scaling limit inside the critical window is a key step in analyzing the MST on the giant component in the supercritical regime;
see \cites{BraBulCohHavSta03,wu2006transport,braunstein2007optimal,chen2006universal} for predictions based on empirical evidence and  \cites{addario2013scaling, addarioberry-sen, bhamidi-sen-heavy-tailed-mst, broutin-marckert-convex-minorant-mst} for some rigorous results. 
We will call the models of random discrete structures that have $\vCrit(\cdot)$ as the metric scaling limit inside the critical window members of the \erdos universality class.

Going back to the models of interest in this paper, hierarchical lattice versions of fundamental models in statistical physics were originally introduced in the context of spin systems in \cite{dyson1969existence,baker1972ising} and as a test bed for renormalization group techniques to study  corresponding questions on the classical lattices \cite{wilson1971renormalization}; 
see \cites{BRBDSG-2019-renormal} for a recent comprehensive survey describing the important role played by hierarchical models in understanding weakly self-avoiding walks and the $\varphi^4$ model.
Of particular relevance to our work is the phase transition in long-range percolation on the hierarchical lattice.
This has been studied \cites{koval-meester-trapman, dawson-gorostiza-percolation-in-ultrametric-space} in the literature when the corresponding kernel is translation-invariant, symmetric, and summable having polynomial decay and, in particular, continuity of the phase transition and uniqueness of the infinite cluster has been established.
The critical behavior in this setting has been more closely investigated in a series of recent papers by Hutchcroft \cites{hutchcroft2022sharp, hutchcroft-two-point, hutchcroft-critical-cluster-volume}, where the two-point function and the cluster sizes at criticality have been studied.
The work \cite{hutchcroft-critical-cluster-volume} considers long-range percolation driven by the kernel $J(x, y) = \| x - y \|^{-\alpha}$, where $d<\alpha<2d$ (see Section~\ref{sec:results-HL} for the relevant notation).
In this regime, establishing the metric scaling limit of the critical clusters inside large balls remains an open problem, although \cite{hutchcroft-critical-cluster-volume}*{Corollary~1.15} proves a tightness result for critical cluster sizes when $4d/3<\alpha<2d$ for a specially defined percolation configuration (see the definition of ``$\eta_{\Lambda}$'' in \cite{hutchcroft-critical-cluster-volume}*{Section~1.1}).
It is conjectured \cite{hutchcroft-critical-cluster-volume}*{Section~7.1} that this model belongs to the \erdos universality class when $d<\alpha\leq 4d/3$.

Percolation on the high-dimensional discrete torus is an important model with nontrivial internal geometry where one expects to see \erdos asymptotics at criticality; see \cite{heydenreich-hofstad-book-high-d-percolation} for a survey of the relevant literature.
It was conjectured by Aizenman \cite{aizenman-nuclear-phys-1997} that when the dimension $d\geq 7$, the size of the critical maximal cluster in nearest neighbor percolation on the discrete torus of volume $V$ scales like $V^{2/3}$, similar to what one observes for the critical \erdos random graph.
This conjecture was settled in \cite{heydenreich-hofstad-high-d-torus-II} for nearest neighbor percolation and for the spread-out model under the following assumption: $d$ is sufficiently large for nearest neighbor percolation, and $d\geq 7$ and the spread is sufficiently large for the spread-out model.
This paper improved upon the previous work \cite{heydenreich-hofstad-high-d-torus}, which contained the correct upper bound, but a suboptimal lower bound on the size of the critical maximal cluster.
Further, under the same assumption, \cite{heydenreich-hofstad-high-d-torus-II} showed that the diameter of the critical maximal clusters scale like $V^{1/3}$, again similar to the critical \erdos random graph.
More recently, it was shown in \cite{blanc-renaudie-nachmias-torus-component-size} that under the same assumption, the scaling limit of the ordered cluster sizes, normalized by $V^{2/3}$, is the same as that of the critical \erdos component sizes identified by Aldous \cite{aldous1997brownian}.
The cycle structure for critical high-dimensional percolation on the discrete torus was studied in \cite{hofstad-sapozhnikov}, where it was shown that under the same assumption as above, the number of cycles with length of the order of $V^{1/3}$ in each maximal cluster forms a tight sequence, provided the ``number'' of such cycles is interpreted in an appropriate way.
Establishing the scaling limit at the level of metric measure spaces for high-dimensional nearest neighbor/spread-out percolation, as well as for long-range percolation on the discrete torus when the effective dimension is high remains an important open problem.

In this paper, we consider long-range percolation on the hierarchical lattice and the discrete torus driven by an underlying kernel $J$ that satisfies $J(x, y)\asymp \| x - y \|^{-\alpha}$, where $0<\alpha<d$ (the precise definitions of the models and the notation used here will be given in Section~\ref{sec:def-and-res}), and identify the point of phase transition for both these models.
When $\alpha\in(0, 5d/6)$, we establish the existence of a critical window, and prove membership of these models in the \erdos universality class by establishing the scaling limit of the critical clusters viewed as metric measure spaces.
We further study the cycle structure in the maximal clusters and show that when $\alpha\in(0, 2d/3)$, any maximal critical cluster inside a large ball is either a tree or has girth of the order of the diameter of the cluster.
This behavior is different from what happens when $\alpha\in (d, 2d)$ where each maximal cluster would have many triangles.
A key component of our proof is obtaining bounds on the two-point function in the barely subcritical regime.
In particular, we show that the two-point function has a plateau in part of the barely subcritical regime; we revisit this point in Remark~\ref{rem:plateau} below.

There are two standard techniques for establishing convergences as in \eqref{eqn:lim-add-br-go}.
One is to explore the random graph in a suitable fashion, and then prove convergence of the associated exploration process and the height process, and deduce the metric scaling limit from there.
The other one is to identify a suitable spanning tree inside each maximal component, establish its scaling limit, and use this to deduce convergence of the components.
These two methods have been successfully applied to the following models of random graphs:
the inhomogeneous random graph with a rank-one kernel (including the \erdos random graph), the configuration model and (simple) random graphs with a given degree sequence,
and random intersection graphs.
See 
\cites{goldschmidt-conchon-stable-graph,addario2012continuum,goldschmidt-haas-senizergues,bhamidi-sen-vacant-set,dhara-hofstad-leeuwaarden-sen-cm-finite-third-moment,riordan2012phase,joseph2010component,sd-rvdh-jvl-ss-cm-component-scaling-infinite-third-moment,SB-SD-vdH-SS-cm-metric-infinite-third-moment,SB-SD-vdH-SS-global-lower-mass,bhamidi-hofstad-sen-multiplicative,broutin-duquesne-wang,miermont-sen,aldous1997brownian,addario2010critical, sbssxw-rank-1-scaling-limit, wang-minmin-random-intersection-graphs}
for various results on the critical behavior in these models proved under different sets of assumptions.
The two techniques above, albeit powerful, have limited applicability in the sense that the model being studied needs to have certain nice features in order to be amenable to analysis via these techniques.
In \cites{sbssxw-graphon, sbssxw-universal}, a more robust technique for establishing such metric scaling limits was developed.
This general technique has found applications in the following problems:
in the study of noise sensitivity in critical random graphs \cite{lubetzky-peled-noise},
in establishing the critical scaling limit of edge-weighted graphs converging to an $L^3$ graphon \cite{sbssxw-graphon}, and in obtaining the critical percolation scaling limit of the hypercube \cite{blanc2024scaling}.
This universality result will be our key tool.

So, the aim of this paper is twofold:
(i) For long-range percolation driven by a class of non-summable kernels on the hierarchical lattice and the discrete torus, we establish the critical percolation (metric) scaling limit of the maximal clusters.
This proves the first such scaling limit results, at the level of metric measure spaces, for these models.
(ii) A second goal is to demonstrate the applicability of the universality result proven in \cites{sbssxw-graphon, sbssxw-universal} to such models.

\subsection{Organization of the paper}
In Section~\ref{sec:topologies}, we briefly describe the notions of weak convergence on the space of compact metric measure spaces needed to state the main results.
Section~\ref{sec:def-and-res} contains the definitions related to the two models we are interested in and the statements of our main results.
We also list a few open problems related to this work in Section~\ref{sec:def-and-res}.
Section~\ref{sec:prelim} contains the notation used in the proofs as well as the construction of the scaling limit $\vCrit(\cdot)$ appearing in \eqref{eqn:lim-add-br-go}.
Finally, Sections~\ref{sec:HL-proofs} and \ref{sec:torus-proofs} contain the proofs of our results on the hierarchical lattice and the discrete torus respectively.

\subsection{Relevant topologies}\label{sec:topologies}
We will work with the Gromov-Hausdorff-Prokhorov (GHP) topology induced by the GHP metric (denoted by $d_{\GHP}$) on the space (say $\cS$) of isometry equivalence classes of compact metric measure spaces.
Note that $(\sS, d_{\GHP})$ is a complete separable metric space  \cite{EJP2116}. 
We refer the reader to \cites{EJP2116,addario2013scaling,burago2001course} for the relevant definitions and background on this topology.

We will write $X = (X, d, \mu)$ to denote both the metric measure space and the corresponding equivalence class. 
In this paper, typically for each $n\geq 1$ we will have a sequence $\vC^{\sss(n)} = (\cC_1^{\sss(n)}, \cC_2^{\sss(n)}, \ldots)$, where $\cC_i^{\sss(n)}$ is a (random) compact metric measure space.  We will view such objects as elements in $\sS^{\bN}$.
 Distributional convergence for a sequence $( \vC^{\sss(n)}\, ;\,  n\geq 1 )$ of $\sS^{\bN}$-valued random objects will be understood with respect to the the product topology on $\sS^{\bN}$ inherited from $d_{\GHP}$.

\vskip4pt

\noindent\textbf{The scaling operator:} 
For $a, b> 0$, let $\scl(a, b):  \sS \to \sS$ be the scaling operator
\begin{align*}
\scl(a, b)\big( (X , d , \mu) \big):= (X, d',\mu'),
\end{align*}
where $d'(x,y) := a d(x,y)$ for all $x,y \in X$, and $\mu'(A) := b \mu(A)$ for all Borel subsets $A \subseteq X$. 
For simplicity, we write $\scl(a, b) X := \scl(a, b)\big( (X , d , \mu) \big)$ and $a X := \scl(a, 1) X$.

\section{Model description, main results, and related open problems}\label{sec:def-and-res}
Sections~\ref{sec:results-HL} and \ref{sec:results-torus} below respectively present the results on the hierarchical lattice and the discrete torus.

\subsection{Main results: hierarchical lattice}\label{sec:results-HL}
We use the formulation in \cites{hutchcroft-critical-cluster-volume, hutchcroft-two-point}.
Fix $d\geq 1, L\geq 2$, and let $\Omega_L^d := (\bZ/ L\bZ)^d$.
Note that $\Omega_L^d$ can be naturally viewed as an abelian group.
Let $\bH_L^d$ be the direct sum of infinitely many copies of $\Omega_L^d$:
\[
\bH_L^d 
:= 
\bigoplus_{i=1}^{\infty} \Omega_L^d
=
\big\{
(x_i\, ;\, i\geq 1) \in \big( \Omega_L^d \big)^{\bN} 
\, :\, 
x_i=0 \text{ for all but finitely many } i
\big\} .
\]
The element $(0, 0, \ldots)\in\bH_L^d$ will simply be denoted by $0$.
For $x=(x_i\, ;\, i\geq 1) \in\bH_L^d$, we let 
\[
\|x\| :=
\begin{cases}
0\, ,\ \text{ if }x_j=0\text{ for all } j\, , \\[3pt]
L^i\, , \ \text{ if }x_i\neq 0=x_{i+1}=x_{i+2}=\ldots \, ,
\end{cases}
\]
and let $\|x-y\|$ be the distance between $x, y\in\bH_L^d$.

\begin{rem}
The symbol ``$\| \cdot \|$'' does not denote a norm here, but it is used so that it is easy to draw parallels with long-range percolation on $\bZ^d$.
Note also that $\Omega_L^d$ and the discrete torus on $L^d$ vertices share the same underlying set.
However, the internal structure of the discrete torus does not play a role in the definition of the distance $\| x - y\|$.
Section~\ref{sec:results-torus} below concerns the results on critical long-range percolation on the discrete torus, and the intrinsic geometry of the torus is taken into account there.
\qed
\end{rem}

For $x\in\bH_L^d$ and $r>0$, let 
$B(x,r) := \{y\in\bH_L^d\, :\, \|y-x\| \leq r \}$.
In this formulation, the space $\bH_L^d$ is ``$d$-dimensional'' in the sense that for any $x\in\bH_L^d$ and any $r\geq 1$, the number of points in $B(x, r)$ satisfies
$r^d L^{-d} \leq |B(x, r)| \leq r^d$. 
Define $\Lambda_0 := \{0\}$ and
\begin{align}\label{eqn:HL-def-Lambda-n}
\Lambda_j := B(0, L^j)\, , \ \  j\geq 1. 
\end{align}
We note two properties of the space $\bH_L^d$ here that will be useful in our proof:
(i) The space $\bH_L^d$ is ultrametric. 
In particular, for any $x, y\in\bH_L^d$ and $r>0$, if $y\in B(x, r)$, then $B(y, r) = B(x, r)$.
(ii) If $x, y\in\Lambda_i$ for some $i\geq 1$, then there exists an isometry of $\Lambda_i$ that maps $x$ to $y$.
In other words, $\Lambda_i$ looks the same viewed from the perspective of $x$ and from the perspective of $y$.

\begin{figure}[t]
\centering
\includegraphics[trim=1.8cm 9.5cm 4.4cm 3.4cm, clip=true, angle=0, scale=.65]{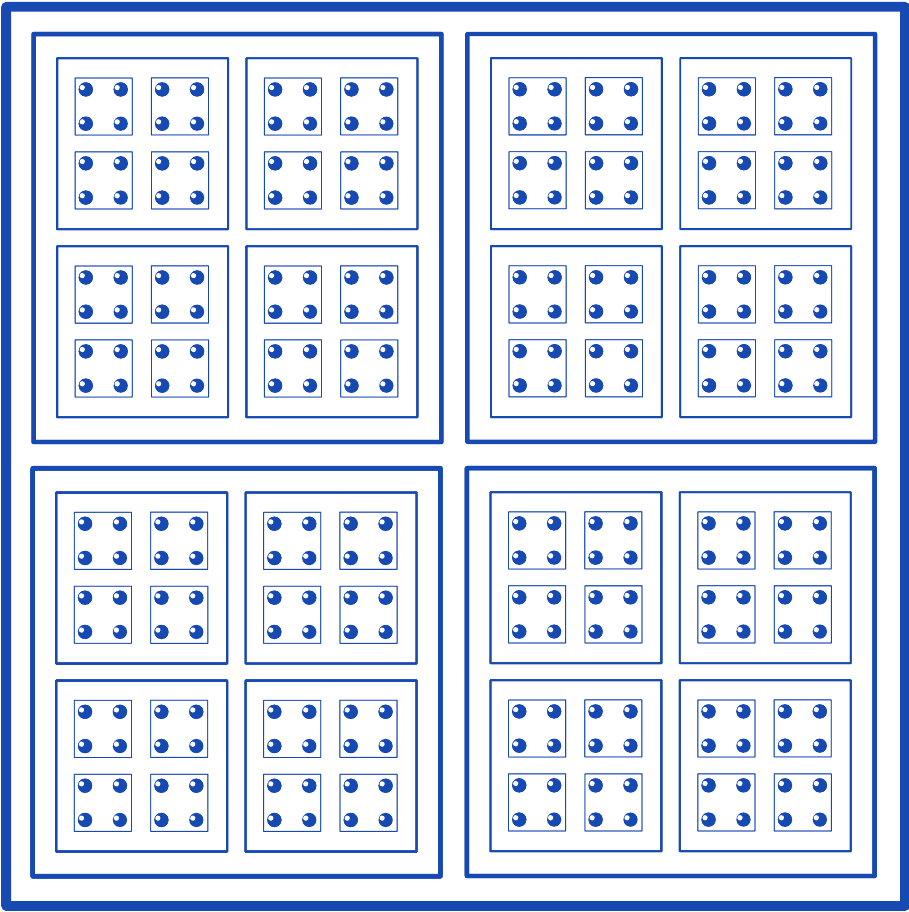}
\caption{
A visualization of the ball $\Lambda_4$ in $\bH_2^2$, where the distance between any two points is the side length ($2, 4, 8$, or $16$) of the smallest outlined square that contains both points.}
\end{figure}

For any symmetric function $\kappa:\bH_L^d\times\bH_L^d \to [0, \infty)$, we write $\HL(\kappa)$ for the random graph with vertex set $\bH_L^d$ obtained by placing edges independently between each pair $x\neq y\in\bH_L^d$ with probability $1-\exp(-\kappa(x, y))$.
For $A\subseteq\bH_L^d$, we will write $\HL_{\sss A}(\kappa)$ for the induced subgraph of $\HL(\kappa)$ on the vertex set $A$.
We denote the $i$-th largest component of $\HL_{\sss\Lambda_n}(\kappa)$ by $\cC_i^{\sss(n)}(\kappa)$, $i\geq 1$.

Now, fix $0<\alpha<d$.
Let $\rho: [0,\infty)\to [0, \infty)$ be such that 
there exist $0<A_1\leq A_2$ for which
\begin{align}\label{eqn:1a}
A_1 r^{-\alpha}\leq \rho(r)\leq A_2  r^{-\alpha} \ \ \text{ for all }\ \ r\geq 1,
\end{align}
and $\rho(0) = 0$.
Let $J:\bH_L^d\times\bH_L^d\to [0, \infty)$ be given by 
$J(x, y) : = \rho( \|x-y\| )$.
Since $\alpha<d$, the kernel $J$ is non-summable, i.e.,
\begin{align}\label{eqn:1}
\sum_{x\in\bH_L^d} J(0, x) = \infty.
\end{align}
Note that \eqref{eqn:1} implies that for any fixed $\beta>0$, each vertex has infinitely many neighbors in $\HL(\beta J)$ almost surely.
Further, for any fixed $\beta>0$, the graphs $\HL_{\sss\Lambda_n}(\beta J)$, $n\geq 1$, are sufficiently dense, so that one has to consider a sequence $\beta_n\to 0$ in order to observe critical phenomena.

Let $\zeta_n$ be the unique positive number that satisfies
\begin{align}\label{eqn:2}
(L^d - 1)\cdot\sum_{i=1}^{n} L^{(i-1)d} \exp\big( - \zeta_n^{-1}\rho(L^i) \big)
=
L^{nd} - 2
\, .
\end{align}
Then there exist (see Lemma~\ref{lem:HL-zeta-n-m-j-n-asymptotics}) $C, C'>0$ depending only on $A_1, A_2, L, d$, and $\alpha$ such that
\begin{align}\label{eqn:3}
C L^{n(d-\alpha)}\leq \zeta_n \leq C' L^{n(d-\alpha)} 
\ \ \text{ for all }\ \ n\, .
\end{align}
The next result shows that $\HL_{\sss\Lambda_n}(\beta J)$ undergoes a phase transition in the largest component around $\beta=\zeta_n^{-1}$.
Recall that $\HL_{\sss\Lambda_n}(\beta J)$ is a random graph on $|\Lambda_n| = L^{nd}$ vertices.

\begin{prop}\label{prop:HL-phase-transition}
Fix $0< \alpha <d$. 
For $t>-1$, let
$
r_{t}^{\sss(n)}(x, y)
:=
(1+t) \cdot \zeta_n^{-1} \cdot J(x, y)
$, 
$x, y\in\bH_L^d$.
\begin{enumeratea}
\item 
Purely subcritical regime:
There exists $C_{\ref{eqn:4}}>0$ depending only on $d$ and $L$ such that for any $\eps\in(0, 1)$,
\begin{align}\label{eqn:4}
\pr\big(\hskip0.2pt
\big| \cC_1^{\sss(n)}\big( r^{\sss(n)}_{-\eps} \big)  \big|
\leq
C_{\ref{eqn:4}} n\eps^{-2}
\big)
\to 1 \ \ \text{ as }\ \ n\to\infty.
\end{align}

\item
Purely supercritical regime:
Fix $\eps>0$. 
Then there exists $C>0$ depending only on $\eps,$ $\alpha,$ $A_1, A_2, d$, and $L$ such that
\begin{align}\label{eqn:4a}
\pr\big(\hskip0.2pt
\big| \cC_1^{\sss(n)}\big( r^{\sss(n)}_{\eps} \big)  \big|
\geq
C L^{nd}
\big)
\to 1 \ \ \text{ as }\ \ n\to\infty.
\end{align}

\end{enumeratea}
\end{prop}

The proof of Proposition~\ref{prop:HL-phase-transition} will be given in Section~\ref{sec:HL-proof-phase-transition}. 
Our aim is to obtain more precise information about the nature of the phase transition on a finer scale.
To this end, we will consider the following kernel to capture the critical window:
For $\lambda\in\bR$ define 
\begin{align}\label{eqn:HL-def-kernel}
\kappa_{\lambda}^{\sss(n)}(x, y)
:= 
\max\bigg\{
\frac{J(x, y)}{\zeta_n} +\frac{ \lambda }{ L^{4nd/3} }
\ , \
0 \bigg\} 
\, , \  \  x, y\in\bH_L^d\, .
\end{align}
Note that \eqref{eqn:1a} and \eqref{eqn:3} imply that even for $\lambda<0$,
$
\min\big\{ \kappa_{\lambda}^{\sss(n)}(x, y)\, :\, x, y\in\bH_L^d,\ y\in x+\Lambda_n  \big\} > 0
$
for all large $n$
(thus, we can ignore the maximum with $0$ appearing in \eqref{eqn:HL-def-kernel} for our purposes). 
That \eqref{eqn:HL-def-kernel} indeed captures the critical window will become clear from Theorem~\ref{thm:HL-scaling-limit} stated below.
An important step in the proof of Theorem~\ref{thm:HL-scaling-limit} will be the study of the two-point function in the barely subcritical regime; 
we will define the kernel corresponding to the barely subcritical regime next.
Consider 
\begin{align}\label{eqn:HL-theta-full-range}
\theta\in (\alpha, \alpha + d/3)\, .
\end{align}
Let
$\rho^-=\rho^{\sss (n), -}:[0, \infty)\to[0, \infty)$ be given by
$
\rho^-(r) := \max\{\, \rho(r) - L^{-n\theta} \, ,\, 0 \,\}
$.
Define
\begin{align}\label{eqn:HL-def-kappa-minus}
\kappa^{\sss (n), -}(x, y) 
:= \zeta_n^{-1}\cdot\rho^-( \|x-y\| ) \, ,
\ \ x, y\in\bH_L^d\, .
\end{align}
For $x, y\in\Lambda_n$, we let $\{x\connects y\}$ denote the event that there is a path between $x$ and $y$ in $\HL_{\sss\Lambda_n}(\kappa^{\sss (n), -})$; 
in particular, we always have $x\connects x$.

\begin{rem}
It can be shown that if $0<\theta\leq\alpha$, then $\HL_{\sss\Lambda_n}(\kappa^{\sss (n), -})$ is purely subcritical satisfying 
$
| \cC_1^{\sss(n)}( \kappa^{\sss (n), -} )  |
=
O_P(n)
$
as $n\to\infty$. 
On the other hand, for $\HL_{\sss\Lambda_n}(\kappa^{\sss (n), -})$ to be outside the critical window, we must have 
$\zeta_n L^{n\theta} \ll L^{4nd/3}$, which in view of \eqref{eqn:3} is equivalent to $\theta < \alpha + d/3$.
Thus, $\HL_{\sss\Lambda_n}(\kappa^{\sss (n), -})$ is barely subcritical for precisely the values of $\theta$ given in \eqref{eqn:HL-theta-full-range}.
\qed
\end{rem}

We will write 
\begin{align}\label{eqn:HL-def-parameters}
\pmtr = (A_1, A_2, L, d, \alpha, \theta)
\ \text{ and }\
\pmtr^{\ast} = (A_1, A_2, L, d, \alpha)	
\end{align}
for the vectors of parameters involved.

\vskip5pt

\noindent{\bf Convention about constants:}
We describe here a convention that we will follow in the rest of Section~\ref{sec:results-HL}, and in Section~\ref{sec:HL-proofs} while writing the proofs of our results on hierarchical percolation.
We will write $C, C'$ etc. to denote positive constants whose values may change from line to line; the values of these constants will depend only on $\pmtr^{\ast}$ as given in \eqref{eqn:HL-def-parameters}.
Some special constants will be indexed by the display numbers where they are first introduced, and if these constants depend on a parameter other than those in $\pmtr^{\ast}$, then that will be clearly mentioned; e.g., the constant $C_{\ref{eqn:767}}(r)$ introduced right before \eqref{eqn:767} depends on $\pmtr^{\ast}$ and $r$, and the constant $C_{\ref{eqn:898}}$ appearing in \eqref{eqn:898} depends only on $\pmtr^{\ast}$.
For a sequence $(\cI_n\, ;\, n\geq 1)$ of index sets and two sets of real numbers
$\{ a_n^{\sss (i)}  : i\in \cI_n, n\geq 1 \}$ and $\{ b_n^{\sss (i)}  : i\in \cI_n, n\geq 1 \}$,
we will write 
``$a_n^{\sss (i)} \lesssim b_n^{\sss (i)}$ for $i\in \cI_n$"
to mean that there exist $n_0\geq 1$ depending only on $\pmtr$ and $C>0$ depending only on $\pmtr^{\ast}$ such that for all $n\geq n_0$, 
$a_n^{\sss (i)} \leq C b_n^{\sss (i)}$ for all $i\in \cI_n$.
We will write 
``$a_n^{\sss (i)} \asymp b_n^{\sss (i)}$ for $i\in \cI_n$"
if $a_n^{\sss (i)} \lesssim b_n^{\sss (i)}$ for $i\in \cI_n$ and $b_n^{\sss (i)} \lesssim a_n^{\sss (i)}$ for $i\in \cI_n$.
When $\cI_n$ is a singleton for each $n$ and we just have two sequences 
$(a_n\, ;\, n\geq 1)$ and $(b_n\, ;\, n\geq 1)$, we will simply write ``$a_n\lesssim b_n$" and ``$a_n\asymp b_n$."

\vskip5pt

The next result concerns the two-point function in the barely subcritical regime.
Let $\diam(\Lambda_n) = L^n$ denote the diameter of $\Lambda_n$, and recall that $|\Lambda_n| = L^{nd}$.

\begin{thm}\label{thm:HL-two-point-}[Two-point function in the barely subcritical regime]
\begin{enumeratea}
\item 
Let $0<\alpha< d/3$ and $\alpha<\theta< 2\alpha$.
Then for $x\in\Lambda_n\setminus\{0\}$,
\begin{align}\label{eqn:788-a}
	\pr(0 \connects x)
	\asymp
	\begin{cases}
		\zeta_n^{-1} \| x\|^{-\alpha}\, ,
		\ \ \text{ if }\ \ 0< \|x\| \leq \diam(\Lambda_n)^{a_0} \, ,
		\\[4pt]
		|\Lambda_n|^{-(d+\alpha - \theta)/d}\, , 
		\text{ if }\ \  \diam(\Lambda_n)^{a_0} < \|x\| \leq \diam(\Lambda_n) \, ,
	\end{cases}
\end{align}
where $a_0 := 2 - \theta\alpha^{-1}$.

\vskip2pt

\item
Let $0 < \alpha < d/3$ and $2\alpha\leq \theta < \alpha + d/3$.
Then
\[
\pr(0 \connects x)
\asymp
|\Lambda_n|^{-(d+\alpha - \theta)/d} 
\ \text{ for }\ x\in\Lambda_n\setminus\{ 0 \}\, .
\]

\item 
Let $d/3\leq \alpha\leq 5d/6$ and $\alpha<\theta<\alpha + d/3$.
Then \eqref{eqn:788-a} continues to hold.

\vskip2pt

\item 
Let $5d/6 < \alpha < d$ and $\alpha < \theta < 2d - \alpha$. 
Then \eqref{eqn:788-a} continues to hold.
\end{enumeratea}
\end{thm}

\begin{rem}\label{rem:plateau}
In the setting of Theorem~\ref{thm:HL-two-point-}(b), the two-point function only varies within a constant multiplier throughout $\Lambda_n\setminus\{0\}$.
In other words, the topology of $\Lambda_n$ does not have a strong influence on the behavior of the two-point function, similar to what one observes for the \erdos random graph.

On the other hand, in the regimes considered in Theorem~\ref{thm:HL-two-point-}(a), (c), and (d), the two-point function has a plateau, as is also known \cite{hutchcroft-michta-slade-torus-plateau} to be the case for high-dimensional percolation (nearest neighbor or spread-out) on the discrete torus.
See also the recent work \cite{slade-plateau} where a general approach for identifying the plateau for the torus two-point function is presented.
The long-range models considered in \cite{slade-plateau} correspond to summable kernels, as opposed to the regime $\alpha\in(0, d)$ studied in this paper.
The result concerning the torus two-point function given in \cite{slade-plateau}*{Theorem~1} is proved under two conditions \cite{slade-plateau}*{Hypotheses~1 and 2} that assume a decay rate on the $\bZ^d$ two-point function and a comparison inequality between the torus two-point function and the corresponding ``unwrapped" two-point function, a quantity that arises naturally while coupling percolation on the discrete torus with percolation on the cubic lattice (see, e.g., \cite{benjamini-schramm-ecp-1996}*{Theorem~1} and \cite{heydenreich-hofstad-high-d-torus}*{Proposition~2.1}).
Note that this approach does not apply to the models being considered in this paper, and the proof of Theorem~\ref{thm:HL-two-point-} and of Theorem~\ref{thm:torus-two-point} stated below use a different argument.
\qed
\end{rem}

\begin{rem}
When $5d/6 < \alpha < d$ and $\alpha < \theta < \alpha + d/3$, our proof will show that 
$\pr( 0 \connects x) \lesssim $ the right side of \eqref{eqn:788-a} for $x\in\Lambda_n\setminus\{0\}$.
However, our proof technique provides a matching lower bound only when $\alpha < \theta < 2d - \alpha$, as indicated in Theorem~\ref{thm:HL-two-point-}(d).
Such a lower bound for $5d/6 < \alpha < d$ and $2d - \alpha \leq \theta < \alpha + d/3$, combined with Theorem~\ref{thm:HL-two-point-}(a)--(d), would complete the picture for the two-point function, up to constants, in the entire barely subcritical regime.
\qed
\end{rem}

The next result deals with the scaling limit in the critical window.
Recall the kernel $\kappa^{\sss (n)}_{\lambda}$ from \eqref{eqn:HL-def-kernel} and the scaling limit $\vCrit(\cdot)$ from \eqref{eqn:lim-add-br-go}.

\begin{thm}\label{thm:HL-scaling-limit}[Scaling limit in the critical window]
Fix $0<\alpha<5d/6$ and $\lambda\in\bR$.
For $i\geq 1$, view $\cC_i^{\sss (n)}\big(\kappa^{\sss (n)}_{\lambda} \big)$ as a metric measure space using the graph distance and the counting measure. 
	Then 
	\begin{equation}\label{eqn:65}
		\Big(
		\scl\big( |\Lambda_n|^{-1/3} ,\, |\Lambda_n|^{-2/3}\big) \,
		\cC_i^{\sss (n)}\big(\kappa^{\sss (n)}_{\lambda} \big)\, ;
		\ i \geq 1 \Big) 
		\weakc 
		\vCrit(\lambda)
		\ \ \text{ as }\ \ n\to\infty .
	\end{equation}
\end{thm}

Note that the convergence in \eqref{eqn:65} holds without imposing any additional boundary conditions on $\Lambda_n$.
Theorem~\ref{thm:HL-scaling-limit} establishes the metric scaling limit for $0<\alpha<5d/6$, whereas the kernel $J$ is non-summable (i.e., \eqref{eqn:1} holds) for $\alpha\leq d$.
It would be interesting to see if the gap can be closed.

\begin{conj}\label{conj:HL-5d/6<alpha<d}
Show that this model is in the \erdos universality class when $5d/6\leq\alpha\leq d$ by proving that the metric scaling limit is Brownian.
\end{conj}

When $d<\alpha<2d$, the kernel $J$ is summable, i.e., $\sum_{x\in\bH_L^d} J(0, x)<\infty$.
In this regime, there exists $0<\beta_c<\infty$ such that a phase transition occurs in $\HL(\beta J)$ at $\beta = \beta_c$
\cites{koval-meester-trapman, dawson-gorostiza-percolation-in-ultrametric-space}.
As mentioned in Section~\ref{sec:intro}, the critical behavior in this regime has recently been studied by Hutchcroft 
\cites{hutchcroft-critical-cluster-volume, hutchcroft-two-point}.
It is conjectured in \cite[Section~7.1]{hutchcroft-critical-cluster-volume} that the model is a member of the \erdos universality class when $d<\alpha<4d/3$.
(Our notation differs from the one used in \cite{hutchcroft-critical-cluster-volume} in that the exponent $\alpha$ in our paper is replaced by $d+\alpha$ in \cite{hutchcroft-critical-cluster-volume}.)
We paraphrase this conjecture below:

\begin{conj}[\cite{hutchcroft-critical-cluster-volume}*{Section~7.1}]\label{conj:HL-summable-kernel}
Show that this model is in the \erdos universality class when $d<\alpha<4d/3$ by proving that the metric scaling limit is Brownian.
\end{conj}

It is further predicted in \cite[Section~7.1]{hutchcroft-critical-cluster-volume} that one still gets the \erdos scaling limit when $\alpha=4d/3$.
However, the mass and the distance scaling for the maximal components, instead of being $L^{2nd/3}$ and $L^{nd/3}$ as in \eqref{eqn:65}, will possibly have some correction factors.
We believe that Conjecture~\ref{conj:HL-summable-kernel} can be proved, at least for $\alpha$ in a suitable subset of $(d,\, 4d/3]$, with the help of our main universality theorem stated in Theorem~\ref{thm:gen-2} below.

For $\alpha\in(d,\, 4d/3)$ and $i\geq 1$, the $i$-th maximal component of $\HL_{\sss\Lambda_n}(\beta_c J)$, with the intrinsic distance scaled by $L^{-nd/3}$, is expected to converge to $\crit_i(0)$ which has a finite (random) number of cycles; see Section~\ref{sec:1} for the construction the spaces $\crit_i(\cdot)$.
Loosely speaking, the number of ``long'' cycles (i.e., those with lengths of the order of $L^{nd/3} = |\Lambda_n|^{1/3}$) in each of the maximal components of $\HL_{\sss\Lambda_n}(\beta_c J)$ is expected to form a tight sequence, if the ``number'' of such cycles is suitably defined.
(As mentioned in Section~\ref{sec:intro}, a similar tightness result for the number of long cycles in critical nearest neighbor/spread-out percolation on the high dimensional discrete torus was proved in \cite{hofstad-sapozhnikov}.)
However, there will, in addition, be many ``short'' cycles with length $o_P(L^{nd/3})$ in the maximal components (in fact, there will be many triangles), and the number of such short cycles in each maximal component will tend to infinity in distribution.
The situation is different for $\alpha\in(0, 2d/3)$ in that only long cycles are found in the maximal components; this is the content of the next result.

For any finite, connected graph $G=(V, E)$, we write $\spls(G)$ for the number of surplus edges in $G$, i.e.,
$\spls(G):=|E|-|V|+1$.
In a similar manner, we can define the random variables $\spls(\crit_i(\cdot))$;
the precise definition will be given in \eqref{eqn:spls-crit-defn} below.
Let $\cL_{\circ}(G)$ (resp. $\cL^{\circ}(G)$) denote the length of the shortest (resp. longest) cycle in $G$ where $G$ is either a finite, connected graph or $\crit_i(\lambda)$, $i\geq 1, \lambda\in\bR$, 
with the convention that $\cL_{\circ}(G)=\infty$ (resp. $\cL^{\circ}(G)=\infty$) if there are no cycles in $G$.

\begin{thm}\label{thm:HL-surplus-convergence}[Cycle structure and complexity]
Fix $\alpha\in(0, 2d/3)$ and $\lambda\in\bR$.
Then as $n\to\infty$,
\begin{gather}
\Big( 
|\Lambda_n|^{-1/3} \cL_{\circ}\big( \cC_i^{\sss(n)}(\kappa_{\lambda}^{\sss(n)}) \big)
\, ;\ i\geq 1
\Big)
\weakc
\Big( 
\cL_{\circ}\big( \crit_i(\lambda) \big) \big)
\, ;\ i\geq 1
\Big)\, ,
\nonumber
\\
\Big( 
|\Lambda_n|^{-1/3} \cL^{\circ}\big( \cC_i^{\sss(n)}(\kappa_{\lambda}^{\sss(n)}) \big)
\, ;\ i\geq 1
\Big)
\weakc
\Big( 
\cL^{\circ}\big( \crit_i(\lambda) \big) \big)
\, ;\ i\geq 1
\Big)\, ,\text{ and}
\nonumber
\\
\big(
\spls\big( \cC_i^{\sss(n)}(\kappa_{\lambda}^{\sss(n)}) 
\big)
\, ;\, 
i\geq 1
\big)
\weakc
\big(
\spls\big( \crit_i(\lambda) \big) 
\, ;\, 
i\geq 1
\big)
\label{eqn:66}
\end{gather}
with respect to the product topology, jointly with the convergence in \eqref{eqn:65}.
\end{thm}

The distribution of the sequence on the right side of \eqref{eqn:66} has an explicit description: 
$\big(
\spls\big( \crit_i(\lambda) \big) 
\, ;\, 
i\geq 1
\big)
\equald
\big(
\pois_i(\lambda)\, ;\, i\geq 1
\big)
$,
where $\pois_i(\cdot)$ is as defined right before Theorem~\ref{thm:aldous-review} stated below.
Note that the claim in Theorem~\ref{thm:HL-surplus-convergence} does not follow directly from Theorem~\ref{thm:HL-scaling-limit}.
We will prove Theorem~\ref{thm:HL-surplus-convergence} by separately showing that for each $i\geq 1$, there are no short cycles in 
$\cC_i^{\sss(n)}(\kappa_{\lambda}^{\sss(n)})$ with probability tending to one.

It remains open to determine the threshold $\alpha_{\sss\text{cy}}(d)$ such that with probability tending to one, there are no short cycles in the maximal critical clusters when $0<\alpha<\alpha_{\sss\text{cy}}(d)$, while for $\alpha_{\sss\text{cy}}(d)<\alpha<4d/3$, the maximal critical clusters contain short cycles.
Theorem~\ref{thm:HL-surplus-convergence} shows that $\alpha_{\sss\text{cy}}(d)\geq 2d/3$.
One can in particular ask if $\alpha_{\sss\text{cy}}(d) = d$ or $\alpha_{\sss\text{cy}}(d) < d$; 
i.e., whether the transition from ``only long cycles in the maximal critical clusters'' to ``there exist short cycles in the maximal clusters'' happens at the cutoff where the kernel goes from being non-summable to summable, or if there is an intermediate regime where, inside the critical window, the maximal clusters contain short cycles while the maximum edge probability goes to zero.

\begin{conj}\label{conj:HL-short-cycles}
Determine the value of $\alpha_{\sss\text{cy}}(d)$.
\end{conj}

\begin{rem}
It is possible for critical systems to have metric scaling limits different from $\vCrit(\cdot)$, and there are several results
\cites{bhamidi-hofstad-sen-multiplicative, broutin-duquesne-wang, goldschmidt-conchon-stable-graph, SB-SD-vdH-SS-cm-metric-infinite-third-moment, wang-minmin-random-intersection-graphs}
in the existing literature identifying such classes of limiting random metric measure spaces.
It is natural to ask if the critical scaling limits for $\alpha\in(4d/3,\, 2d)$, if they exist, are rescaled versions of these spaces, or if these lead to completely new universality classes.
\qed
\end{rem}

\subsection{Main results: discrete torus}\label{sec:results-torus}
Fix $d\geq 1$, and for each $n\geq 3$ consider the equivalence relation $\stackrel{\sss n}{\sim}$ on $\bZ^d$ where $x, y\in\bZ^d$ are equivalent (i.e., $x\stackrel{\sss n}{\sim} y$) if and only if $x-y\in n\bZ^d$.
Let $\bT_n^d := \bZ^d/\stackrel{\sss n}{\sim}$ be the $d$-dimensional discrete torus on $n^d$ vertices.
Equivalently, $\bT_n^d$ can be viewed as the cube $\{0, 1, \ldots, n-1 \}^d$ with periodic boundary conditions.
Let $\|x-y\|$ denote the distance between $x, y\in\bT_n^d$ under the quotient metric on $\bT_n^d$ inherited from the $\ell^{\infty}$ distance on $\bZ^d$, and write $\| x \| = \| x - 0 \|$.

For any symmetric function $\kappa:\bT_n^d\times\bT_n^d \to [0, \infty)$, we write $\torus_n(\kappa)$ for the random graph with vertex set $\bT_n^d$ obtained by placing edges independently between each pair $x\neq y\in\bT_n^d$ with probability $1-\exp(-\kappa(x, y))$.
We denote the $i$-th largest component of $\torus_n(\kappa)$ by $\cC_i^{\sss(n)}(\kappa)$, $i\geq 1$.

Now, fix $\alpha\in (0, d)$, and let $\rho$ be as introduced around \eqref{eqn:1a}.
Let $J(x, y) : = \rho( \| x-y\|)$, $x, y\in\bT_n^d$.
Note that the diameter of $\bT_n^d$ is 
\begin{align}\label{eqn:torus-diameter}
d_n := \lfloor n/2 \rfloor\, .
\end{align}
For each $n\geq 3$, let $\zeta_n$ be the unique positive number satisfying
\[
\sum_{i=1}^{d_n - 1}
\big( (2i + 1)^d - (2i - 1)^d \big) \exp\big( - (\zeta_n)^{-1} \rho( i ) \big)
+ 
\big( n^d - (2d_n - 1)^d \big) \exp\big( - (\zeta_n)^{-1} \rho( d_n ) \big)
=
n^d - 2\, .
\]
Then there exist (see Lemma~\ref{lem:torus-zeta-n-m-j-n-asymptotics}) $C, C'>0$ depending only on $A_1, A_2, d$, and $\alpha$ such that
$
C n^{d-\alpha}\leq \zeta_n \leq C' n^{d-\alpha}
$
for all $n\geq 3$.
The next result is the analogue of Proposition~\ref{prop:HL-phase-transition}, and concerns the phase transition in this model.

\begin{prop}\label{prop:torus-phase-transition}
	Fix $0< \alpha <d$. 
	For $t>-1$, let
	$
	r_{t}^{\sss(n)}(x, y)
	:=
	(1+t) \cdot \zeta_n^{-1} \cdot J(x, y)
	$, 
	$x, y\in\bT_n^d$.
\begin{enumeratea}
		\item 
		Purely subcritical regime:
		There exists $C>0$ depending only on $d$ such that for any $\eps\in(0, 1)$,
		\begin{align*}
			\pr\big(\hskip0.2pt
			\big| \cC_1^{\sss(n)}\big( r^{\sss(n)}_{-\eps} \big)  \big|
			\leq
			C \eps^{-2}\log n
			\big)
			\to 1 \ \ \text{ as }\ \ n\to\infty.
		\end{align*}

		\item
		Purely supercritical regime:
		Fix $\eps>0$. 
		Then there exists $C>0$ depending only on $\eps,$ $\alpha,$ $A_1, A_2$, and $d$ such that
		\begin{align*}
			\pr\big(\hskip0.2pt
			\big| \cC_1^{\sss(n)}\big( r^{\sss(n)}_{\eps} \big)  \big|
			\geq
			C n^d
			\big)
			\to 1 \ \ \text{ as }\ \ n\to\infty.
		\end{align*}

\end{enumeratea}
\end{prop}

To better understand what happens near the critical point, fix $\lambda\in\bR$ and define, similar to \eqref{eqn:HL-def-kernel}, the kernel
\begin{align*}
\kappa_{\lambda}^{\sss (n)} (x, y)
:= 
\max\bigg\{
\frac{J(x, y)}{\zeta_n } +\frac{ \lambda }{ n^{4d/3} }
\ , \
0 \bigg\} 
\, , \  \  x \neq y\in\bT_n^d\, .
\end{align*}
(Similar to the comment made below \eqref{eqn:HL-def-kernel}, even when $\lambda<0$, we can ignore the maximum with zero for all large $n$.)
As was the case for the hierarchical lattice, a key step in our analysis of the critical window will be the study of the two-point function in the barely subcritical regime.
Let $\theta$ be as in \eqref{eqn:HL-theta-full-range}.
Define
$\rho^-=\rho^{\sss (n), -}:[0, \infty)\to[0, \infty)$ as
$
\rho^-(r) = \max\{\, \rho(r) - n^{-\theta} \, ,\, 0 \,\}
$.
Let
\begin{align}\label{eqn:torus-def-fp}
\kappa^{\sss (n), -}(x, y) 
:= \zeta_n^{-1}\cdot\rho^-( \|x-y\| ) \, ,
\ \ x, y\in\bT_n^d .
\end{align}
For $x, y\in\bT_n^d$, we let $\{x\connects y\}$ denote the event that there is a path in $\torus_n(\kappa^{\sss (n), -})$ connecting $x$ and $y$ (we always have the relation $x\connects x$).
Let
\begin{align}\label{eqn:95}
\pmtr := (A_1, A_2, d, \alpha, \theta)
\ \text{ and }\
\pmtr^{\ast} := (A_1, A_2, d, \alpha)\, .
\end{align}
In the rest of Section~\ref{sec:results-torus} and in Section~\ref{sec:torus-proofs}, we will adopt the convention about constants described right below \eqref{eqn:HL-def-parameters} with $\pmtr$ and $\pmtr^{\ast}$ as in \eqref{eqn:95}.
The next result concerns the two-point function in $\torus_n(\kappa^{\sss (n), -})$.

\begin{thm}\label{thm:torus-two-point}
The results in Theorem~\ref{thm:HL-two-point-}(a)--(d) continue to hold if $\Lambda_n$ is replaced by $\bT_n^d$ in their statements.
\end{thm}

Here, $|\Lambda_n|$ and $\diam(\Lambda_n)$ from Theorem~\ref{thm:HL-two-point-} are to be replaced by $|\bT_n^d| = n^d$ and $\diam(\bT_n^d) = d_n = \lfloor n/2 \rfloor$ respectively.
The following result concerns the scaling limit of the maximal components inside the critical window for this model.

\begin{thm}\label{thm:torus-scaling-limit}
Fix $0<\alpha<5d/6$ and $\lambda\in\bR$.
Then the conclusion of Theorem~\ref{thm:HL-scaling-limit} continues to hold if $\Lambda_n$ is replaced by $\bT_n^d$ in the statement.
\end{thm}

Let $\spls(\cdot)$, $\cL_{\circ}(\cdot)$, and $\cL^{\circ}(\cdot)$ be as introduced right before Theorem~\ref{thm:HL-surplus-convergence}.
The next result concerns the girth and the complexity of the critical maximal components.

\begin{thm}\label{thm:torus-surplus-convergence}
Fix $\alpha\in(0, 2d/3)$ and $\lambda\in\bR$.
Then the conclusion of Theorem~\ref{thm:HL-surplus-convergence} continues to hold if $\Lambda_n$ is replaced by $\bT_n^d$ in the statement.
\end{thm}

The analogues of open problems~\ref{conj:HL-5d/6<alpha<d} and \ref{conj:HL-short-cycles} remain open for the discrete torus as well.
A result analogous to the one claimed in open problem~\ref{conj:HL-summable-kernel} should also hold for effectively high dimensional long-range percolation on the discrete torus.
A natural guess is that the effective dimension in this context is high when either $d\geq 7$, or when $1\leq d\leq 6$ and $\alpha\in (d, 4d/3]$.

\begin{rem}
Fix $\lambda\in\bR$ and let $\mvxi(\lambda)$ be as in Theorem~\ref{thm:aldous-review} stated below.
In Theorem~\ref{thm:torus-scaling-limit}, considering the convergence of just the total measures of each maximal component yields
\begin{align}\label{eqn:909}
\big( n^{-2d/3}| \cC_i^{\sss (n)}(\kappa_{\lambda}^{\sss (n)}) | \, ;\ i\geq 1 \big)
\weakc
\mvxi(\lambda)
\ \ \text{ as }\ \ n\to\infty
\end{align}
whenever $d\geq 1$, $\alpha\in (0, 5d/6)$, and the kernel $J$ satisfies $J(x, y)\asymp \| x - y \|^{-\alpha}$ for $x\neq y\in\bT_n^d$.

The convergence in \eqref{eqn:909}, in the special case 
$d = 2$, $\alpha =1$, and $J(x, y) = C \| x - y \|^{-1}$, $x\neq y\in\bT_n^2$, 
was recently proved in \cite{turova-2-d-torus} using the exploration process technique mentioned in Section~\ref{sec:intro}.
In \cite{turova-2-d-torus}*{Section~1.3}, the authors write
``We can only conclude that treating the general model [instead of the special case $d = 2$, $\alpha =1$] at criticality remains a challenge," and predict that \eqref{eqn:909} should hold for any $d\geq 1$ and $\alpha\in (0, 2d/3)$.

Our Theorem~\ref{thm:torus-scaling-limit} shows that the convergence takes place in a more general sense--at the level of metric measure spaces--for any $d\geq 1$ and $\alpha\in (0, 5d/6)$.
Also, as mentioned above, one expects this convergence to hold for a larger set of values of $(d, \alpha)$.
\qed
\end{rem}

\noindent{\bf A comment on the presentation of proofs:}
The model of hierarchical percolation has radial symmetry, a property that percolation on the discrete torus lacks.
Most of the techniques employed in the proofs of our results on the hierarchical lattice that do not make use of this radial symmetry can be applied in a straightforward way while working with the discrete torus.
The proofs for both these models will be carried out in the following steps:
(i) First, we prove an upper bound on the two-point function in the barely subcritical regime. 
This upper bound is then used to bound the quantities $\Delta_n$ and $\tilde\Delta_n$ defined in \eqref{eqn:HL-def-Delta} below.
(ii) Use these bounds to then prove a lower bound on the two-point function.
(iii) Use the bounds from (i) together with the universality theorem stated in Theorem~\ref{thm:gen-2} below to obtain the critical scaling limit, and the girth and complexity results.
We will write down complete proofs of (i), (ii), and (iii) for the hierarchical lattice.
For the discrete torus, we will detail the argument for (i), and only outline the steps for (ii) and (iii) as they are similar to those for hierarchical percolation.

\section{Notation and limit objects}\label{sec:prelim}
Section~\ref{sec:gr-constr} lists some notation frequently used in the sequel.
Sections~\ref{sec:cont-limit-descp} and \ref{sec:smc-def} contain the necessary background on real trees and the multiplicative coalescent.
In Section~\ref{sec:1}, we describe the construction of the limiting object $\vCrit(\cdot)$ appearing in \eqref{eqn:lim-add-br-go}.
Section~\ref{sec:universality-res-recall} contains the statement of the universality theorem used in our proofs.

\subsection{Graph theoretic functionals and asymptotic notation}\label{sec:gr-constr}
We write $|A|$ or $\# A$ for the cardinality of a set $A$.
For $n\in\bZ_{>0}$, write $[n]$ for the set $\{1, 2, \ldots, n\}$.
For any rooted tree $\vt$, write $\height(\vt)$ for the height of $\vt$.
For a finite graph $G$, we write $V(G)$ and $E(G)$ for the set of vertices and the set of edges respectively.
We will write $|G|$ to mean $|V(G)|$. 
For any graph $G$ and $v\in V(G)$, we will denote the component of $G$ that contains $v$ by $\cC(v ; G)$.
As already mentioned right before the statement of Theorem~\ref{thm:HL-surplus-convergence}, for a connected graph $G$, we write $\spls(G) := |E(G)| - |G|+1$ for the number of surplus edges in $G$.
We view a connected component $\cC$ of $G$ as a metric space using the graph distance $d_{G}$.
When the graph $G$ is clear from the context, we suppress the subscript and simply write $d$. 
Often in our analysis, associated to a graph $G$, there will be a collection of positive vertex weights 
$\vw = \big\{w_v: v\in V(G)\big\}$. 
Two natural measures on the vertex set of $G$ are:
\begin{inparaenuma}
\item 
{\it Counting measure:} 
$\mu_{\text{ct}}(A): = |A|$, for $A \subseteq V(G)$.  
\item 
{\it Weighted measure:} $\mu_\vw(A) :=  \sum_{v \in A} w_v$, for $A \subseteq V(G)$. 
\end{inparaenuma}
For a finite connected graph $G$, we use $G$ for both the graph and the associated metric measure space.	
For two graphs $G_1$ and $G_2$, write $G_1\subseteq G_2$ to mean that $G_1$ is a subgraph of $G_2$.
For a metric space $(\fX, d)$, write $\diam(\fX):=\sup_{u, v\in\fX}d(u, v)$.
For a finite graph $G$, we let 
$\diam(G):=\max\big\{\diam(\cC)\, :\, \cC\text{ connected component of }G\big\}$.

We omit ceilings and floors when there is no confusion in doing so.
For two real sequences  $(a_n)$ and $(b_n)$, write $a_n\sim b_n$ to mean $a_n/b_n\to 1$ as $n\to\infty$.
For two sequences $(X_n)$ and $(Y_n)$ of random variables with $X_n, Y_n$  defined on the same probability space for each $n$,  write $X_n\sim Y_n$ to mean that $X_n/Y_n\weakc 1$ as $n\to\infty$.
For a probability distribution $\mu$,  write $X\sim\mu$ to mean that the random variable $X$ has law $\mu$. 
Write $X_n=O_P(Y_n)$ if the sequence $(Y_n^{-1}X_n)$ is tight,
and  $X_n=o_P(Y_n)$ if $Y_n^{-1}X_n \probc 0$ as $n\to\infty$.
We say that $X_n=\Theta_P(Y_n)$ if $X_n=O_P(Y_n)$ and $Y_n=O_P(X_n)$.
We say that a sequence $\cE_n,\, n\geq 1$, of events occurs with high probability (whp) if $\pr(\cE_n)\to 1$ as $n\to\infty$.
We will use the symbol ``$\disjt$" to denote disjoint occurrence of events; the definition and the relevant background on this can be found in \cite{grimmett-percolation-book}.

\subsection{Real trees and quotient spaces}
\label{sec:cont-limit-descp}
For $l>0$, write
$\cE_l$ for the space of continuous excursions on $[0,l]$. 
Fix $h,g \in \cE_l$, and a locally finite set $\cP \subseteq \bR_+\times \bR_+$. 
Define
\begin{equation*}
g \cap \cP := \big\{(x,y) \in \cP: 0 \leq x \leq l, \; 0 \leq y < g(x)  \big\}.
\end{equation*}
Let $\cT(h)$ be the real tree encoded by $h$ (see, e.g., \cite{evans-book,legall-book}) equipped with the push forward of the Lebesgue measure on $[0,l]$. 
For $(x,y) \in g \cap \cP$, let $r(x,y) := \inf\set{x': x' \geq x, \; g(x') \leq y}$.
Write $\cG(h,g,\cP)$ for the metric measure space obtained by identifying the pairs of points in $\cT(h)$ corresponding to the pairs of points $\{(x,r(x,y)) : (x,y) \in g \cap \cP\}$.
Thus $\cG(h,g,\cP)$ is obtained by adding a finite number of shortcuts to the real tree $\cT(h)$.

\vskip3pt

\noindent \textbf{Tilted Brownian excursions:} 
Let $(\ve_l(s), s \in [0,l])$ be a Brownian excursion of length $l$. 
Fix $\theta >0$ and let $\tilde \ve_l^\theta$ be an $\cE_l$-valued random variable such that for any bounded continuous function $f : \cE_l \to \bR$,
\begin{equation}
\label{eqn:tilt-exc-def}
\E[f(\tilde \ve_l^\theta)] = { \E\Big[f(\ve_l) \exp\Big(\theta\medint\int_0^l \ve_l(s)ds\Big) \Big] }\big/{\E\Big[\exp\Big(\theta\medint\int_0^l \ve_l(s)ds \Big)\Big]}.
\end{equation}
Write $\ve(\cdot) := \ve_1(\cdot)$, $\tilde \ve^\theta (\cdot) := \tilde \ve_1^{\theta} (\cdot)$, and $\tilde \ve_l (\cdot) := \tilde \ve_l^1 (\cdot)$.
In Section~\ref{sec:1}, we will use the above ingredients to construct the limiting random metric measure spaces in \eqref{eqn:lim-add-br-go}.

\subsection{Multiplicative coalescent and the random graph $\cG(\vx,q)$}
\label{sec:smc-def}
A multiplicative coalescent $\big(\vX(t);\, t\in \cA\big)$ is a Markov process with state space, 
\[
\ldown := \big\{
(x_1,x_2,\ldots): x_1\geq x_2\geq \cdots \geq 0, \sum_i x_i^2< \infty 
\big\}, 
\]
endowed with the metric inherited from $l^2$;
here, either $\cA=\bR$ or $\cA=[t_0, \infty)$ for some $t_0\in\bR$.
Its evolution described in words is as follows:
Fix $\vx\in \ldown$.
Then conditional on $\vX(t) = \vx$, each pair of clusters $i$ and $j$ merge at rate $x_i x_j$ to form a new cluster of size $x_i+x_j$. While this description makes sense for a finite collection of clusters (i.e., $x_i = 0$ for $i> K$ for some finite $K$), Aldous \cite{aldous1997brownian} showed that \chh{it also} makes sense for $\vx\in \ldown$.

An object closely related to the multiplicative coalescent is the random graph $\cG(\vx,q)$:  Fix the vertex set $[n]$, a collection of positive vertex weights $\vx=(x_i, i\in [n])$, and a parameter $q>0$. Construct the random graph $\cG(\vx, q)$ by placing an edge between $i\neq j\in [n]$ with probability
$1 - \exp(-q x_i x_j)$,
independently across pairs $\{i, j\}$ with $i\neq j$.
For a connected component $\cC$ of $\cG(\vx, q)$, define 
$\mass (\cC) = \sum_{i \in \cC} x_i={\mu_{\vx}(\cC)}$ 
{(recall the weighted measure $\mu_{\vx}$ from Section~\ref{sec:gr-constr})}. 
Rank the components in decreasing order of their masses, and let $\cC_i$ be the $i$-th maximal component. Consider a sequence of vertex weights $\vx = \vx^{\sss(n)}$ and a sequence $q = q^{\sss(n)}$.
Let 
\begin{align}\label{eqn:def-sigma_k-x_max}
\sigma_k : = \sum_{i \in [n]} x_i^k\ \ \text{ for } k \geq 1\, ,\ \ 
\text{ and }\ \ 
x_{\max} := \max_{i\in [n]} x_i\, .
\end{align}

\begin{ass}\label{ass:aldous-basic-assumption}
Assume that there exists a constant $\lambda \in \bR$ such that as $n \to \infty$,
\begin{equation*}
{\sigma_3}/{(\sigma_2)^3} \to 1, 
\ \ 
q - {(\sigma_2)^{-1}} \to \lambda, 
\ \ \text{ and }\ \ 
{x_{\max}}/{\sigma_2} \to 0.
\end{equation*}
\end{ass}
Note that these conditions imply that $\sigma_2 \to 0$. 
Fix $\lambda \in \bR$ and let $B$ denote standard Brownian motion.  
Define the processes $W_\lambda$ and $\tilde{W}_\lambda$ via
\begin{equation}
\label{eqn:parabolic-bm}
W_\lambda (t) := B(t) + \lambda t - t^2/2\, ,\ \ \text{ and }\ \  
\tilde W_\lambda (t) := W_\lambda (t) - \inf_{s \in [0,t]} W_\lambda (s)\, ,\ \ t\geq 0\, .
\end{equation}
An excursion of $\tilde W_\lambda$ is an interval $(l,r) \subset \bR^+$ such that $\tilde W_\lambda(l)=\tilde W_\lambda(r) = 0$ and $\tilde W_\lambda(t)>0$ for all $t \in (l,r)$. Write $r-l$ for the length of such an excursion. 
Aldous \cite{aldous1997brownian} showed that the lengths of the excursions of $\tilde W_\lambda$ can be arranged in decreasing order as
$\gamma_1(\lambda) > \gamma_2(\lambda) > \ldots > 0$.
Conditional on $\tilde W_\lambda$, let $\pois_i(\lambda)$, $i\geq 1$, be independent random variables where $\pois_i(\lambda)$ is Poisson distributed with mean equal to the area underneath the excursion corresponding to  $\gamma_i(\lambda)$. Recall the definition of surplus ($\spls$) of a connected graph from Section~\ref{sec:gr-constr}. 

\begin{thm}[\cite{aldous1997brownian}]\label{thm:aldous-review}
Under Condition \ref{ass:aldous-basic-assumption}, 	$\left(\mass( \cC_i); i \geq 1 \right) \weakc \mvxi(\lambda):= (\gamma_i(\lambda); i\geq 1)$,  as $n \to \infty$, with respect to the topology on $\ldown$.
Further, 
$\big(\big(\mass( \cC_i),\spls(\cC_i) \big);\, i \geq 1 \big) 
\weakc 
\mvXi(\lambda):=\big( \big(\gamma_i(\lambda), \pois_i(\lambda)\big);\, i\geq 1)$,  as $n \to \infty$, with respect to the product topology.
\end{thm}

\begin{rem}
In \cite{aldous1997brownian}*{Proposition 4}, Aldous only proves the first assertion in Theorem \ref{thm:aldous-review}, namely convergence of the rank-ordered masses. The second assertion in Theorem \ref{thm:aldous-review} is proved for the special case of the \erdos random graph \cite{aldous1997brownian}*{Corollary 2}.
However, the second convergence in Theorem \ref{thm:aldous-review} follows easily from \cite{aldous1997brownian}*{Proposition 10}.
\qed
\end{rem}

\subsection{Scaling limit of components in the critical \erdos random graph}\label{sec:1}
Now we can describe the scaling limit of the maximal components in
$\ERRG(n, 1/n+\lambda/n^{4/3})$ derived in \cite{addario2012continuum}. 
Recall the definitions of tilted Brownian excursions and the metric measure space $\cG(h,g,\cP)$ from Section \ref{sec:cont-limit-descp}. 
Let $\mvxi(\lambda) = (\gamma_i(\lambda); i\geq 1)$ be as in Theorem \ref{thm:aldous-review}. 
Conditional on $\mvxi(\lambda)$, let $\tilde \ve_{\gamma_i(\lambda)}$, $i\geq 1$, be independent tilted Brownian excursions with $\tilde \ve_{\gamma_i(\lambda)}$ having length $\gamma_i(\lambda)$.  
Let $\cP_i$, $i\geq 1$, be independent rate one Poisson processes on $\bR_+^2$ that are also independent of $\big(\tilde \ve_{\gamma_i(\lambda)};\,  i\geq 1\big)$. 
Define the sequence of random metric measure spaces $\vCrit(\lambda) := (\crit_i(\lambda)\, ;\, i\geq 1)$, where
\begin{equation}\label{eqn:limit-metric-def}
\crit_i(\lambda)
:= 
\cG\big(
2\tilde \ve_{\gamma_i(\lambda)}, \tilde \ve_{\gamma_i(\lambda)}, \cP_i
\big)\, , \ \ \ i\geq 1.
\end{equation}
We define
\begin{align}\label{eqn:spls-crit-defn}
\spls\big( \crit_i(\lambda) \big)
:=
\big| \tilde \ve_{\gamma_i(\lambda)} \cap \cP_i \big| \, , \ \ i\geq 1\, .
\end{align}

\subsection{The universality theorem}\label{sec:universality-res-recall}
 Recall the random graph $\cG(\vx,q)$ from Section~\ref{sec:smc-def}.
View $(\cC_i\, ;\, i \geq 1)$, ordered by mass,  as metric measure spaces equipped with the graph distance and the weighted measure $\mu_{\vx}$ (recall the weighted measure from Section~\ref{sec:gr-constr}). 
The first result describes the metric scaling limit of the maximal components.


\begin{ass}\label{ass:gen-1}
Assume that there exist $\eta_0 \in (0,\infty)$ and $r_0 \in (0,\infty)$ such that ${x_{\max}}/{\sigma_2^{3/2+\eta_0}} \to 0$ and  ${\sigma_2^{r_0}}/{x_{\min}} \to 0$
as $n \to \infty$, where $x_{\min}:=\min_{i\in [n]}x_i$.
\end{ass}

\begin{thm}[{\cite{sbssxw-graphon}*{Theorem 3.2}}]\label{thm:gen-1}
Under Conditions~\ref{ass:aldous-basic-assumption} and \ref{ass:gen-1}, 
the maximal components satisfy
$\big(\scl(\sigma_2, 1 ) \cC_i\, ;\,  i \geq 1 \big) \weakc \vCrit(\lambda)$ as $n\to\infty$. 
\end{thm}

To extend this result to one applicable to different models of random discrete structures, we replace each vertex $i \in [n]$ in the graph $\cG(\vx,q)$ with a metric measure space $(M_i, d_i, \mu_i)$; 
we refer to the spaces $(M_i, d_i, \mu_i)$, $i\in [n]$, as ``blobs.''
In Theorem~\ref{thm:gen-2}, we show that the metric measure space $\bar\cC_i$ now associated with $\cC_i$, under Conditions~\ref{ass:aldous-basic-assumption} and \ref{ass:gen-2}, converges to $\crit_i(\lambda)$ after proper rescaling, owing to macroscopic averaging of distances within blobs. 
To apply this result to a specific model of critical random discrete structure, we declare its components in the barely subcritical regime to be the blobs, and then couple the dynamics while going from the barely subcritical stage to the critical window with the multiplicative coalescent.
We now describe the ingredients required to incorporate the internal blob level geometry.

\vskip4pt

\noindent{\upshape (a)} {\bf Blob level superstructure: }  A simple finite graph $\cG$ with vertex set $[n]$ and vertex weight sequence $\vx:=(x_i, i \in [n])$. 

\vskip3pt

\noindent{\upshape (b)} {\bf Blobs:} A family of compact metric measure spaces 
$\vM := \{(M_i, d_i, \mu_i): i \in [n] \}$, 
where $\mu_i$ is a probability measure for each $i$.

\vskip3pt

\noindent{\upshape (c)} {\bf Blob to blob junction points:} A collection of points $\vX := \{X_{i,j}: i \in [n], j \in [n]\}$ such that $X_{i,j} \in M_i$ for all $i,j$.

\vskip3pt

Using these three ingredients,  we can define a metric measure space $\Gamma(\cG,\vx,\vM,\vX) = (\bar M, \bar d, \bar \mu)$ as follows: Let $\bar M := \bigsqcup_{i \in [n]} M_i$. Define the measure $\bar \mu$ as
\begin{equation}
\bar \mu( A ) = \sum_{i \in [n]} x_i \mu_i(A \cap M_i), \mbox{ for } A \subseteq \bar M.
\label{eqn:barmu-on-full-met}
\end{equation}
The metric $\bar d$ is obtained by using the intra-blob
distance functions $\big(d_i;\, i\in [n]\big)$ together with the graph distance on $\cG$ by putting an edge of length one between the pairs of vertices
$\{ \{X_{i,j}, X_{j,i} \}: \{i,j \} \mbox{ is an edge in } \cG \}$.
Thus, for $x, y \in \bar M$ with $x \in M_{j_1}$ and $y \in M_{j_2}$, 
\begin{equation}\label{eqn:44}
\bar d(x,y) = \inf\Big\{ k +  d_{j_1}(x, X_{j_1,i_1}) + \sum_{\ell=1}^{k-1} d_{i_\ell}(X_{i_\ell, i_{\ell-1}}, X_{i_\ell, i_{\ell+1}}) + d_{j_2}(X_{j_2, i_{k-1}}, y) \Big\},
\end{equation}
where the infimum is taken over $k\geq 1$ and all paths $(i_0, i_1,\ldots,i_{k-1}, i_k)$ in $\cG$ with $i_0 =j_1$ and $i_k = j_2$. 
Here, the infimum of an empty set is understood to be $+\infty$.
This corresponds to the case where $j_1$ and $j_2$ do not belong to the same component.

The above is a deterministic procedure for creating a new metric measure space. 
Now assume that we are provided with a parameter sequence $q^{\sss(n)}$, weight sequence $\vx^{\sss(n)} := (x_i^{\sss(n)} ; i \in [n])$, and  family of compact metric measure spaces $\vM^{\sss(n)} = \big((M_i^{\sss(n)}, d_i^{\sss(n)} , \mu_i^{\sss(n)}) ; i \in [n]\big)$.
As before, we will suppress the dependence on $n$. 
Let $\cG(\vx, q)$ be the random graph defined in Section~\ref{sec:smc-def} constructed using the weight sequence $\vx$ and the parameter $q$.   Let $(\cC_i ; i \geq 1)$ denote the  connected components of $\cG(\vx,q)$ ranked in decreasing order of their masses. 
Let $X_{i,j}$, $i,j \in [n]$,  be independent random variables (that are also independent of $\cG(\vx,q)$) such that for each fixed $i$, $X_{i,j}$, $j\in [n]$, are i.i.d. $M_i$-valued random variables that are $\mu_i$ distributed. 
Let $\vX = (X_{i,j},\ i,j \in [n])$.

We define $\bar \cG(\vx,q,\vM) := \Gamma\big( \cG(\vx, q), \vx, \vM, \vX \big)$.
Let $\bar \cC_i$ be the component in $\bar \cG(\vx, q, \vM)$ corresponding to the $i$-th maximal component $\cC_i$ in $\cG(\vx,q)$, viewed as a compact metric measure space as explained in \eqref{eqn:barmu-on-full-met} and \eqref{eqn:44}. 
Define the quantities
\begin{gather}
u_{i} := \int_{M_i \times M_i} d_i(x,y)  \mu_i (dx) \mu_i (dy), \qquad i\in [n],
\label{eqn:uik-def}\\
\mbox{$\tau := \sum_{i \in [n]} x_i^2 u_i$ and ${\diam_{\max}} := \max_{i \in [n]} \diam(M_i)$.}
\label{eqn:tau-dmax-def}
\end{gather}
\chh{Note that $u_i$ is the expectation of {the distance between} two blob-to-blob junction points in blob $i$ when these are sampled independently according to the measure $\mu_i$. Further,  $\tau/\sigma_2$ is the weighted average of these ``typical'' distances. }
\begin{ass}\label{ass:gen-2}
	There exist $\eta_0 \in (0,\infty)$ and $r_0 \in (0,\infty)$ such that
	{\upshape (a)} Condition \ref{ass:gen-1} holds, and
	{\upshape (b)}  as $n \to \infty$, $ (\diam_{\max}\cdot \sigma_2^{3/2-\eta_0})/{ (\tau + \sigma_2)} \to 0$ and  $ {\sigma_2 x_{\max} {\diam_{\max}}}/{\tau} \to 0$. 	
\end{ass}
\begin{thm}[{\cite{sbssxw-graphon}*{Theorem 3.4}}]\label{thm:gen-2}
Under Conditions~\ref{ass:aldous-basic-assumption} and \ref{ass:gen-2},
	\begin{equation*}
	\Big(\scl\Big(\frac{\sigma_2^2}{\sigma_2 + \tau}, 1 \Big) \bar \cC_i\, ;\, i \geq 1 \Big) \weakc \vCrit(\lambda)\, ,\ \mbox{ as } n \to \infty.
	\end{equation*}
\end{thm}

Further, we have the following result about the cycle structure in $\bar{\cC}_i$:

\begin{thm}\label{thm:cycle-structure}
Fix $i_*\geq 1$.
Assume that Conditions~\ref{ass:aldous-basic-assumption} and \ref{ass:gen-2} hold, and further,
the following holds whp:
$M_j$ is a real tree for each $j\in [n]$ for which $M_j\subseteq \bar{\cC}_{i_*}$.
Let $\cL_{\circ}( \bar{\cC}_{i_*} )$ (resp. $\cL^{\circ}( \bar{\cC}_{i_*} )$) denote the length of the shortest (resp. longest) cycle in $\bar{\cC}_{i_*}$.
Then as $n\to\infty$,
\[
\bigg( \frac{\sigma_2^2}{\sigma_2 + \tau} \bigg) \cdot \cL_{\circ}(\bar{\cC}_{i_*}) 
\weakc 
\cL_{\circ}\big(\crit_{i_*}(\lambda) \big) 
\ \ \text{ and }\ \
\bigg( \frac{\sigma_2^2}{\sigma_2 + \tau} \bigg) \cdot \cL^{\circ}(\bar{\cC}_{i_*}) 
\weakc 
\cL^{\circ}\big(\crit_{i_*}(\lambda) \big) \, .
\]
\end{thm}

The result in Theorem~\ref{thm:cycle-structure} is not explicitly stated in \cite{sbssxw-universal, sbssxw-graphon}.
However, this result follows directly from the proof of \cite{sbssxw-rank-1-scaling-limit}*{Theorem~7.3} and the proof of \cite{sbssxw-graphon}*{Display (64)}.

\section{Proofs: hierarchical lattice}\label{sec:HL-proofs}
We will prove
Theorems~\ref{thm:HL-two-point-}, \ref{thm:HL-scaling-limit}, \ref{thm:HL-surplus-convergence}, and Proposition~\ref{prop:HL-phase-transition} in Section~\ref{sec:HL-proofs}.
In Section~\ref{sec:HL-proof-branching-process}, we note certain properties of a collection of branching processes that will be useful in our proofs.
The proof of Theorem~\ref{thm:HL-two-point-} will be given in Section~\ref{sec:HL-proof-two-point}.
In Sections~\ref{sec:HL-proof-xmax-diamax}--\ref{sec:HL-proof-tau}, we gather all the ingredients necessary to prove Theorems~\ref{thm:HL-scaling-limit} and \ref{thm:HL-surplus-convergence}.
The proofs of Theorems~\ref{thm:HL-scaling-limit} and \ref{thm:HL-surplus-convergence} are then completed in Sections~\ref{sec:HL-proof-scaling-limit} and \ref{sec:HL-proof-surplus} respectively.
Finally, the proof of Proposition~\ref{prop:HL-phase-transition} is given in Section~\ref{sec:HL-proof-phase-transition}.
Let us recall here the convention about constants that was described right after \eqref{eqn:HL-def-parameters}.

We now introduce some notation that will be used throughout Section~\ref{sec:HL-proofs}.
For $A\subseteq\bH_L^d$, we will write $G_{\sss A}^-$ for $\HL_{\sss A} (\kappa^{\sss (n) , - }  )$ to simplify notation.
Suppose $G_{\sss \Lambda_n}^-$ has $m_n$ components.
We will write 
$\cC_i^-=\cC_i^{\sss (n), -}$ for the $i$-th largest component in $G_{\sss \Lambda_n}^-$, 
i.e., $\cC_i^- = \cC_i^{\sss (n), -} := \cC_i^{\sss (n)}(\kappa^{\sss (n), -})$, $1\leq i\leq m_n$.
In the proof of Theorem~\ref{thm:HL-scaling-limit}, we will apply Theorem~\ref{thm:gen-2} with
$\vx=(x_1, \ldots, x_{m_n}), q$, and $\vM=(M_1, \ldots, M_{m_n})$ given by
\begin{equation}\label{eqn:510}
	x_i = L^{-2nd/3}|\cC_i^-| \, ,\ \
	q = \lambda + \zeta_n^{-1} L^{n\big( \frac{4d}{3} - \theta \big)}  \, ,
	\ \text{ and }\
	M_i =  \scl\big(1, 1/|\cC_i^-| \big) \cC_i^-\, ,\ \ 1\leq i\leq m_n\, .
\end{equation}
(Thus, $x_1,\ldots, x_{m_n}$ are the component sizes in $G_{\sss \Lambda_n}^-$, suitably rescaled.
We will also use `$x$' to denote points in $\bH_L^d$.
There should not be any confusion in this regard as the meaning will be clear from the context.)
With this notation fixed for the rest of Section~\ref{sec:HL-proofs}, we let 
$\sigma_2, \sigma_3, x_{\max}, \diamax$, and $\tau$ be as in 
\eqref{eqn:def-sigma_k-x_max}, \eqref{eqn:tau-dmax-def}.

For $x\in\Lambda_n$, we will write $\cC^{ -}(x) = \cC^{\sss (n), -}(x)$ for the component of $x$ in 
$G_{\sss \Lambda_n}^-$, i.e., 
$\cC^-(x) = \cC\big(x; G_{\sss \Lambda_n}^- \big)$.
For $x, y\in\Lambda_n$, we will write $\{x\edge y\}$ for the event that there is an edge between $x$ and $y$ in $G_{\sss \Lambda_n}^-$.
As defined right below \eqref{eqn:HL-def-kappa-minus}, $\{x\connects y\}$ 
denotes the event that there is a path between $x$ and $y$ in $G_{\sss \Lambda_n}^-$, i.e., $y\in\cC^-(x)$.
Similarly, for $A\subseteq\Lambda_n$ and $x, y\in A$, we let 
$\{ x\connects y \text{ in } A\}$ be the event that there is a path in $G_{\sss A}^-$ connecting $x$ and $y$.
Thus, 
$\{ x\connects y \text{ in } \Lambda_n\}
=
\{ x\connects y\}$, 
and 
$\{z\in A \, :\, x\connects z \text{ in } A\}$
is the vertex set of the component of $x$ in $G_{\sss A}^-$, namely 
$\cC\big(x; G_{\sss A}^-\big)$.

Recall the function $\rho^-$ introduced around \eqref{eqn:HL-def-kappa-minus}.
Note that since $\alpha<\theta$, \eqref{eqn:1a} implies that there exists $n_0\geq 1$ depending only on $A_1, L, \alpha$, and $\theta$ such that for all $n\geq n_0$,
$\kappa^{\sss (n), -}(x, y) 
= \zeta_n^{-1}\big( \rho( \|x-y\| ) - L^{-n\theta} \big)
$
whenever $x-y\in\Lambda_n\setminus\{0\}$; 
i.e., we can ignore the maximum with $0$ in the definition of $\rho^-$ for our purposes.
We will often use the relation $\rho^-(r) = \rho(r) - L^{-n\theta}$, $L\leq r\leq L^n$, without explicitly mentioning that this is really true for $n\geq n_0$.
Define
\begin{align}\label{eqn:HL-def-fp}
	\fp_{r}^{\sss (n), -}
	:=
	1 - \exp\big( -\zeta_n^{-1}\cdot\rho^-(r) \big) 
	\, ,\ \ r>0\, .
\end{align}
Using \eqref{eqn:1a}, we see that
\begin{align}\label{eqn:HL-upper-bound-connection-prob}
\pr(x\edge y)
=
\fp_{\| x-y \|}^{\sss (n), -}
\leq
A_2\zeta_n^{-1} \| x-y \|^{-\alpha} , \ \ x\neq y\in\bH_L^d\, ,
\end{align}
a relation we will use frequently.

\subsection{A family of branching processes}\label{sec:HL-proof-branching-process}
For each $n\geq 1$, let $X_i^{\sss (n)}$, $1\leq i\leq n$, be independent random variables such that 
$X_i^{\sss (n)}
\sim
\text{Binomial}\big(L^{id} - L^{(i-1)d},\ \fp_{L^i}^{\sss (n), -}\big)$,
where $\fp_{\cdot}^{\sss (n), -}$ is as defined in \eqref{eqn:HL-def-fp}.
For any $n\geq 1$ and $1\leq j\leq n$, let $\cT_j^{\sss (n)}$ be a branching process tree where the number of children of each vertex has the same distribution as $\sum_{i=1}^j X_i^{\sss (n)}$.
We will write $\cW_k^{\sss (n)}$ for the number of vertices in the $k$-th generation of $\cT_n^{\sss (n)}$, $k\geq 0$.
Thus, $\cW_0^{\sss (n)} = 1$ and $\cW_1^{\sss (n)} \equald \sum_{i=1}^n X_i^{\sss (n)}$.
The trees $\cT_j^{\sss (n)}$, $1\leq j\leq n$, arise naturally while studying $G_{\sss \Lambda_n}^-$ as we explain next.

Fix $1\leq j\leq n$ and let $x\in\Lambda_j$.
Explore $\cC\big(x; G_{\sss \Lambda_j}^- \big)$ in a breadth-first manner to build its breadth-first spanning tree (BFS tree).
The exploration starts at $x$ which is designated as the root of the BFS tree.
In each step of this exploration, one looks for edges present in $G_{\sss \Lambda_j}^-$  between the current vertex being explored, say $y$, and all vertices of $\Lambda_j$ that have not been discovered up to that point.
The vertices found are added as children of $y$ in the BFS tree and they are ordered from left to right using some deterministic rule.
The vertices in each generation of the BFS tree are explored from left to right, and the exploration process moves to the next generation when all vertices in a generation have been explored.
The process terminates when all vertices in a generation have been explored, and there are no vertices in the next generation.

There is an obvious way of coupling this BFS tree with $\cT_j^{\sss (n)}$:
For each $y$ in the BFS tree and each $z$ in the BFS tree that was discovered before the exploration of the vertex $y$ (this includes the vertex $x$ if $y\neq x$), we give $y$ an extra child with probability 
$\fp^{\sss (n), -}_{\|y-z\|}$, and we do this independently for each such pair $(y, z)$ of vertices.
Then we identify each new vertex added to the BFS tree with the root of a copy of $\cT_j^{\sss (n)}$, where the copies of $\cT_j^{\sss (n)}$ are jointly independent.
(In other words, we run independent branching processes starting from each new vertex added to the BFS tree.)
The resulting tree will have the same root as the BFS tree and will contain the BFS tree as a subtree, and further, the resulting tree will have the same distribution as $\cT_j^{\sss (n)}$ (the latter claim uses the properties (i) and (ii) mentioned right below \eqref{eqn:HL-def-Lambda-n}).
In particular, 
\begin{align}\label{eqn:coupling}
	\diam\big( \cC^-(x) \big)
	\leq
	2\cdot\height(\cT_n^{\sss (n)})
	\ \ \text{ and }\ \
	|\cC(x; G_{\sss \Lambda_j}^-) | 
	\leq 
	|\cT_j^{\sss (n)}|
\end{align}
in this coupling.

The mean number of offsprings for each vertex in $\cT_j^{\sss (n)}$ is given by 
\begin{align}\label{eqn:HL-def-m-j-n}
	m^{\sss (n)}_j = (L^d - 1)\sum_{i=1}^j L^{(i-1)d}\fp_{L^i}^{\sss (n), -} .
\end{align}
The following asymptotics for $\zeta_n$ and $m_j^{\sss (n)}$, $1\leq j\leq n$, will be used repeatedly in our proof.

\begin{lem}\label{lem:HL-zeta-n-m-j-n-asymptotics}
Fix $\alpha\in (0, d)$ and $\theta\in (\alpha, \alpha + d/3)$.
Then the following hold:
\begin{enumeratea}
\item 
There exist constants $C, C'>0$ depending only on $\pmtr^{\ast}$ such that
$C L^{n(d-\alpha)} \leq \zeta_n \leq C' L^{n(d-\alpha)}$ for all $n$.
\item 
We have, 
$1-m_n^{\sss(n)} = \big(1 + O(L^{-nd})\big)\cdot \zeta_n^{-1} L^{n(d-\theta)}$ as $n\to\infty$.
In particular, $1-m_n^{\sss(n)} \asymp L^{-n(\theta - \alpha)}$.
\item 
We have, 
$m_j^{\sss(n)} \lesssim \zeta_n^{-1} L^{ j (d-\alpha)}$ for $1\leq j\leq n-1$.
\item
There exist $n_0\geq 2$ depending only on $\pmtr$ and $\delta\in(0, 1)$ depending only on $\pmtr^{\ast}$ such that
$m_j^{\sss (n)}\leq 1-\delta$
for all $n\geq n_0$ and $1\leq j\leq n-1$.
\end{enumeratea}
\end{lem}

\noindent{\bf Proof:}
It follows from \eqref{eqn:2} that $\zeta_n$ is the unique positive number that satisfies
\begin{align}\label{eqn:30}
	(L^d - 1)\cdot\sum_{i=1}^{n} 
	L^{(i-1)d} \big[ 1 - \exp\big( - \zeta_n^{-1}\rho(L^i) \big) \big]
	=
	1\, .
\end{align}
Using the relation \eqref{eqn:1a} and the fact that $1-e^{-u}\geq (1-e^{-1})u$ for $u\in (0, 1]$, we see that for any $x\geq A_2 L^{-\alpha}$,
\begin{align}\label{eqn:31}
	\sum_{i=1}^n L^{(i-1)d}\big[ 1 - \exp\big( -x^{-1}\rho(L^i) \big)\big]
	&
	\geq 
	x^{-1}(1-e^{-1})\sum_{i=1}^n L^{(i-1)d} A_1 L^{-i\alpha}
	\notag\\
	&
	\geq 
	C x^{-1} \sum_{i=1}^n L^{i(d-\alpha)}
	\geq 
	C' x^{-1} L^{n(d-\alpha)}\, .
\end{align}
A similar calculation using the relation \eqref{eqn:1a} and the fact that $1-e^{-u}\leq u$ for $u>0$ will show that for any $x>0$,
\begin{align}\label{eqn:32}
	\sum_{i=1}^n L^{(i-1)d}\big[ 1 - \exp\big( -x^{-1}\rho(L^i) \big)\big]
	\leq 
	C x^{-1} L^{n(d-\alpha)}\, .
\end{align}
Combining \eqref{eqn:30}, \eqref{eqn:31}, and \eqref{eqn:32} yields the claim in (a).

To prove (b), note that
\begin{align}\label{eqn:m_n-leq-1}
	1-m_n^{\sss (n)}
	&
	=
	(L^d - 1)\sum_{i=1}^{n} 
	L^{(i-1)d} \big[ 1 - \exp\big( - \zeta_n^{-1}\rho(L^i) \big) - \fp_{L^i}^{\sss (n), -}\big]
	\notag\\
	&\hskip20pt
	=
	(L^d - 1)\sum_{i=1}^{n} 
	L^{(i-1)d} \exp\big( - \zeta_n^{-1}\rho(L^i) \big)
	\cdot
	\big[ \exp\big( \zeta_n^{-1} L^{-n\theta} \big) - 1\big]
	\notag\\
	&\hskip40pt
	=
	(L^{nd} - 2)\big[ \exp\big( \zeta_n^{-1} L^{-n\theta} \big) - 1\big]
	\, ,
\end{align}
where the first step uses \eqref{eqn:HL-def-m-j-n} and \eqref{eqn:30}, 
the second step uses \eqref{eqn:HL-def-fp},
and the third step uses \eqref{eqn:2}.
Now, using the relation $\zeta_n\asymp L^{n(d-\alpha)}$, we get
\begin{align*}
	1-m_n^{\sss(n)}
	=
	(L^{nd} - 2)\big[  \zeta_n^{-1} L^{-n\theta}  + O\big( \zeta_n^{-2} L^{-2n\theta} \big)\big]
	=
	\zeta_n^{-1} L^{n(d - \theta)}\big( 1 + O(L^{-nd}) \big)
	\asymp
	L^{-n(\theta - \alpha)}\, .
\end{align*}

The claim in (c) follows from the observation that
\[
m_j^{\sss (n)}
\lesssim
\sum_{i=1}^j L^{id} \cdot \zeta_n^{-1} \rho^-(L^i)
\lesssim
\zeta_n^{-1} \sum_{i=1}^j L^{id} \cdot L^{-i\alpha}
\lesssim
\zeta_n^{-1} L^{j(d- \alpha)}
\]
for $1\leq j\leq n-1$, where the second step uses \eqref{eqn:1a}.

Turning finally to (d), we note that  
$m_j^{\sss (n)} \leq m_n^{\sss (n)} - a_n $
for any $1\leq j\leq n-1$, where
\[
a_n 
= 
(L^d - 1) L^{(n-1)d}\fp_{L^n}^{\sss (n), -}
\asymp
L^{nd} \cdot \zeta_n^{-1}\cdot L^{-n\alpha}
\, .
\]
The claim in (d) now follows from the relations $m_n^{\sss (n)} \leq 1$ (which is a consequence of the result in (b)) and $\zeta_n\asymp L^{n(d-\alpha)}$.
\qed

\subsection{Proof of Theorem~\ref{thm:HL-two-point-}}\label{sec:HL-proof-two-point}

We break up the proof into two parts: proving the upper bound and proving the lower bound.

\begin{prop}\label{prop:HL-two-point-bound}[Upper bound on the two-point function]
\begin{enumeratea}
\item
Let $\alpha\in (0, d)$ and $\alpha < \theta < \min\{ 2\alpha, \alpha + d/3 \}$.
Then for $x\in\Lambda_n\setminus\{0\}$,
\begin{align}\label{eqn:788}
	\pr(0 \connects x)
	\lesssim
	\begin{cases}
		\zeta_n^{-1} \| x\|^{-\alpha}\, ,\ \ \text{ if }\ \ L\leq \|x\| \leq L^{n a_0}\, ,
		\\[3pt]
		L^{-n( d + \alpha - \theta ) }\, , \text{ if }\ \ L^{n a_0}< \|x\| \leq L^n\, ,
	\end{cases}
\end{align}
where $a_0 = 2 - \theta\alpha^{-1}$.

\item 
Let $\alpha\in (0, d/3)$ and $2\alpha \leq \theta < \alpha + d/3 $.
Then
\begin{align}\label{eqn:789}
\pr(0 \connects x)
\lesssim
L^{-n( d + \alpha - \theta ) }
\ \ \text{ for }\ \ x\in\Lambda_n\setminus\{0\}.
\end{align}

\end{enumeratea}

\end{prop}

\begin{prop}\label{prop:HL-two-point-lower-bound}[Lower bound]
Suppose either $\alpha\in (0, 5d/6]$ and $\theta\in (\alpha, \alpha + d/3)$, or 
$\alpha\in (5d/6, d)$ and $\theta\in ( \alpha, 2d - \alpha )$.
Then
\begin{align}\label{eqn:790}
	\pr(0 \connects x)
	\gtrsim
	L^{-n( d + \alpha - \theta ) }
	\ \ \text{ for }\ \ x\in\Lambda_n\setminus\{0\}.
\end{align}
\end{prop}

Note that when $\alpha\in (0, d)$ and $\theta\in (\alpha, \alpha + d/3 )$,
$\pr(0 \connects x) 
\geq 
\pr(0 \edge x) 
\gtrsim 
\zeta_n^{-1} \| x \|^{-\alpha}$ 
for $x\in\Lambda_n\setminus\{0\}$,
which coupled with \eqref{eqn:790} yields
\begin{align}\label{eqn:791}
\pr(0 \connects x)
\gtrsim
\max\,\{\,
\zeta_n^{-1} \| x \|^{-\alpha} ,\,
L^{-n( d + \alpha - \theta ) }
\}
\ \ \text{ for }\ \ x\in\Lambda_n\setminus\{0\}
\end{align}
when $\alpha$ and $\theta$ are in the region specified in Proposition~\ref{prop:HL-two-point-lower-bound}.
It follows from the relation \eqref{eqn:3} that \eqref{eqn:791} gives the lower bound part of the claims made in Theorem~\ref{thm:HL-two-point-}(a)--(d).
This combined with the upper bounds provided by Proposition~\ref{prop:HL-two-point-bound} completes the proof of Theorem~\ref{thm:HL-two-point-}.
It remains to prove Propositions~\ref{prop:HL-two-point-bound} and \ref{prop:HL-two-point-lower-bound}.

\vskip5pt

\noindent{\bf Proof of Proposition~\ref{prop:HL-two-point-bound}:}
Let
\[
\annlarge_{\, j} := 
\Lambda_j\setminus\Lambda_{j-1}
\, ,\ \  j\geq 1. 
\]
(We use the symbol ``\,$\annlarge$\," as it resembles an annulus formed by two concentric Euclidean $\ell^{\infty}$-balls.)
For $z\in\annlarge_{\, j}$, $1\leq j\leq n$,
\begin{align}\label{eqn:89}
	&
	\pr( 0\connects z\text{ in }\Lambda_j )
	=
	\big( \#\,\annlarge_{\, j} \big)^{-1} \sum_{y\in\annsmall_{\, j}} \pr(0\connects y\text{ in }\Lambda_j)
	\leq
	\big( \#\,\annlarge_{\, j} \big)^{-1}  \sum_{y\in\Lambda_j\setminus\{0\}} \pr(0\connects y\text{ in }\Lambda_j)
	\notag\\
	&
	\hskip50pt
	\lesssim
	L^{-jd} \bE\big( \, \big| \cC\big(0; G_{\sss \Lambda_j}^- \big) \big|  - 1 \big)
	\leq
	L^{-jd} \bE\big( | \cT_j^{\sss(n)} |  - 1 \big)
	=
	L^{-jd} m_j^{\sss (n)} (1-m_j^{\sss(n)})^{-1} \, ,
\end{align}
where the first step uses symmetry, the third step uses the relation $\#\, \annlarge_{\, j} \asymp L^{jd}$ for $j\geq 1$, and the fourth step uses \eqref{eqn:coupling}.
Combining \eqref{eqn:89} with Lemma~\ref{lem:HL-zeta-n-m-j-n-asymptotics} yields
\begin{align}
	&\hskip30pt
	\pr(0\connects z)
	\lesssim
	L^{-n( d+\alpha-\theta )}  \ \text{ for } z\in\annlarge_{\, n}\, , \ \text{ and}
	\label{eqn:90}\\
	&
	\pr(0\connects z\text{ in }\Lambda_j)
	\lesssim
	\zeta_n^{-1} L^{-j\alpha} \ \ \text{ for } z\in\annlarge_{\, j}\, ,\  1\leq j\leq n-1\, .
	\label{eqn:90-a}
\end{align}

Now consider $x\in\annlarge_{\, i}$, where $1\leq i\leq n-1$.
Let 
$I_x := \min\{ j\in\{i, i+1, \ldots, n\} \, :\, 0\connects x\text{ in }\Lambda_j\}$ with the convention that $I_x = \infty$ if $x\notin\cC^-(0)$.
Then
\begin{align}\label{eqn:91}
\pr(0\connects x)
\leq
\pr(0\connects x\text{ in }\Lambda_i) 
+ 
\sum_{j=i + 1}^n\pr(I_x = j)
\lesssim
\zeta_n^{-1}L^{-i\alpha} 
+ 
\sum_{j=i + 1}^n\pr(I_x = j)
\end{align}
for $x\in\annlarge_{\, i}$, $1\leq i\leq n-1$, where the last step uses \eqref{eqn:90-a}.

Now, if $y\in\Lambda_{j-1}$ and $z\in\annlarge_{\, j}$, then $y-z\in\annlarge_{\, j}$.
Hence, for $x\in\annlarge_{\, i}$, $1\leq i\leq n-1$, and for $i+1\leq j\leq n$,
\begin{align}\label{eqn:920}
&
\pr( I_x = j )
\leq
\sum_{y\in\Lambda_{j-1}}\sum_{z\in\annsmall_{\, j} }
\pr\big(
\{0\connects y \text{ in } \Lambda_{j-1}\}
\disjt
\{y\edge z\}
\disjt
\{z\connects x \text{ in }\Lambda_j\}
\big)
\notag\\
&
\hskip19pt
\leq
\sum_{y\in\Lambda_{j-1}}\sum_{z\in\annsmall_{\, j} }
\pr\big( 0\connects y \text{ in } \Lambda_{j-1} \big)
\cdot
\pr( y\edge z )
\cdot
\pr\big( z\connects x \text{ in }\Lambda_j \big)
\notag\\
&
\hskip38pt
\lesssim
\zeta_n^{-1} L^{-j\alpha} 
\sum_{y\in\Lambda_{j-1}}
\pr\big( 0\connects y \text{ in } \Lambda_{j-1} \big)
\sum_{z\in\annsmall_{\, j} }
\pr\big( z\connects x \text{ in }\Lambda_j \big)
\notag\\
&
\hskip57pt
\leq
\zeta_n^{-1} L^{-j\alpha} 
\cdot
\bE| \cC(0; G_{\sss \Lambda_{j-1} }^-) |
\cdot
\bE| \cC(x; G_{\sss \Lambda_j}^-) | 
\lesssim
\zeta_n^{-1} L^{-j\alpha} 
\cdot
\bE|\cT_j^{\sss (n)}| \, ,
\end{align}
where the second step uses the BK inequality \cites{van-den-berg-kesten-bk-inequality, grimmett-percolation-book, reimer-bkr-inequality},
the third step uses \eqref{eqn:HL-upper-bound-connection-prob},
and the last step uses \eqref{eqn:coupling} and Lemma~\ref{lem:HL-zeta-n-m-j-n-asymptotics}(d).
Hence, for $x\in\annlarge_{\, i}$, $1\leq i\leq n-1$, 
\begin{align}\label{eqn:921}
\pr(0\connects x)
\lesssim
\zeta_n^{-1}L^{-i\alpha} 
+ 
\sum_{j=i+1}^{n-1} \zeta_n^{-1}L^{-j\alpha}
+
\zeta_n^{-1}L^{-n\alpha} \cdot L^{n(\theta - \alpha)}
\lesssim
\zeta_n^{-1}L^{-i\alpha} 
+ 
L^{-n (d+\alpha-\theta)}\, ,
\end{align}
where the first step uses \eqref{eqn:91}, \eqref{eqn:920}, Lemma~\ref{lem:HL-zeta-n-m-j-n-asymptotics}(b) and (d), and the second step follows from Lemma~\ref{lem:HL-zeta-n-m-j-n-asymptotics}(a).
The proof is completed upon noting that the bound in \eqref{eqn:921} is equivalent to \eqref{eqn:788} (resp. \eqref{eqn:789}) when $\alpha$ and $\theta$ are in the region specified in Proposition~\ref{prop:HL-two-point-bound}(a) (resp. Proposition~\ref{prop:HL-two-point-bound}(b)).
\qed

\vskip5pt

We will now prove the claimed lower bound in Proposition~\ref{prop:HL-two-point-lower-bound} with the help of the upper bound in Proposition~\ref{prop:HL-two-point-bound}.
Let
\begin{align}\label{eqn:HL-def-Delta}
	\Delta_n
	:=
	\hskip-1pt 
	\mathop{\sum\sum}\limits_{x\neq y\in\Lambda_n\setminus\{0\}} 
	\hskip-1pt 
	\pr( 0 \connects x ) \pr( 0 \connects y ) \pr( x\edge y )
	\ \ \text{ and }\ \ 
	\widetilde\Delta_n
	:=
	\hskip-1pt 
	\sum_{w\in\Lambda_n\setminus\{0\}} 
	\hskip-1pt 
	\pr( 0 \connects w ) \pr(0\edge w)\, .
\end{align}
The quantity $\Delta_n$ (resp. $\tilde\Delta_n$) is somewhat similar to the triangle diagram (resp. bubble diagram)
\cites{aizenman-newman-tree-graph-inequalities, heydenreich-hofstad-book-high-d-percolation, heydenreich-hofstad-sakai-mean-field}  arising in the study of high-dimensional percolation.
The difference is the presence of the factor $\pr(x\edge y)$ (resp. $\pr(0\edge w)$) instead of $\pr(x\connects y)$ (resp. $\pr(0\connects w)$).
We will now use Proposition~\ref{prop:HL-two-point-bound} to establish bounds on $\Delta_n$ and $\widetilde\Delta_n$.
These bounds will be needed in the later sections as well.

\begin{lem}\label{lem:2}
We have,
\begin{align*}
\Delta_n 
\lesssim
\begin{cases}
|\Lambda_n|^{-(d + 2\alpha - 2\theta)/d}\, , \ \ \text{ if } \ 
0<\alpha\leq 2d/3\ \text{ and }\ \alpha<\theta<\alpha + d/3\, , \\
|\Lambda_n|^{-3(d - \alpha)/d}\, , \ \ \text{ if }\ 
2d/3<\alpha<8d/9\ \text{ and }\ \alpha<\theta< 5\alpha/2 - d\, , \\
|\Lambda_n|^{-(d + 2\alpha - 2\theta)/d}\, , \ \ \text{ if } \ 
2d/3<\alpha< 8d/9\ \text{ and }\ 5\alpha/2 - d\leq\theta<\alpha + d/3\, , \\
|\Lambda_n|^{-3(d - \alpha)/d}\, , \ \ \text{ if }\ 
8d/9\leq \alpha< d\ \text{ and }\ \alpha<\theta<\alpha + d/3\, ,
\end{cases}
\end{align*}
and
\begin{align*}
\tilde\Delta_n 
\lesssim
\begin{cases}
|\Lambda_n|^{-(d + \alpha - \theta)/d}\, , \ \ \text{ if } \ 
0<\alpha\leq d/2\ \text{ and }\ \alpha<\theta<\alpha + d/3\, , \\
|\Lambda_n|^{-2(d - \alpha)/d}\, , \ \ \text{ if }\ 
d/2<\alpha<2d/3\ \text{ and }\ \alpha<\theta< 3\alpha - d\, , \\
|\Lambda_n|^{-(d + \alpha - \theta)/d}\, , \ \ \text{ if } \ 
d/2<\alpha< 2d/3\ \text{ and }\ 3\alpha - d\leq\theta<\alpha + d/3\, , \\
|\Lambda_n|^{-2(d - \alpha)/d}\, , \ \ \text{ if }\ 
2d/3\leq \alpha< d\ \text{ and }\ \alpha<\theta<\alpha + d/3\, .
\end{cases}
\end{align*}
\end{lem}

\noindent{\bf Proof:}
We focus on the bound on $\Delta_n$ when $0<\alpha\leq 2d/3$; 
the other cases can be handled similarly.
If $\alpha\in (0, d/3)$ and $2\alpha\leq\theta<\alpha + d/3$, then
\[
\Delta_n
\lesssim
L^{-2n( d + \alpha - \theta )}
\hskip-5pt
\sum_{ x\in\Lambda_n\setminus\{0\} } 
\sum_{ y\in\Lambda_n\setminus\{0, x\} }
\pr( x\edge y )
\leq
L^{-2n( d + \alpha - \theta )}
\hskip-5pt
\sum_{ x\in\Lambda_n\setminus\{0\} }  m_n^{\sss (n) }
\leq 
L^{-n ( d + 2\alpha - 2\theta )}\, ,
\]
as desired, where the first step uses \eqref{eqn:789}, and the last step uses the facts that $m_n^{\sss (n)}\leq 1$ (see \eqref{eqn:m_n-leq-1}) and $|\Lambda_n| = L^{nd}$.
Next, suppose $0<\alpha\leq 2d/3$ and $\alpha<\theta<\min\{ 2\alpha,\, \alpha + d/3\}$. 
Note that
\begin{gather}\label{eqn:567}
\zeta_n\Delta_n \asymp F_1 + F_2\, ,\ \text{ where}
\\
F_1
=
\sum_{0<\|y\|<\| x\|}\,  \pr( 0 \connects x ) \pr( 0 \connects y )  J(x, y)
\  \  \  \text{ and }\  \  \
F_2
=
\sum_{0<\|y\|=\| x\|}\,  \pr( 0 \connects x ) \pr( 0 \connects y )  J(x, y)\, . 
\notag
\end{gather}
Now,
\begin{align*}
	F_1
	\lesssim
	\sum_{i=1}^n \sum_{x\in\annsmall_{\, i}} \pr(0\connects x) \cdot L^{-i\alpha}
	\bigg( \sum_{j=1}^{i-1}  \sum_{y\in\annsmall_{\, j}} \pr(0\connects y) \bigg)\, ,
\end{align*}
where we have used \eqref{eqn:1a} to conclude that $J(x, y)\lesssim L^{-i\alpha}$ for $x\in\annlarge_{\, i}$ and $y\in\annlarge_{\, j}$, $j<i$.
Similarly,
\begin{align*}
	F_2
	&
	=
	\sum_{i=1}^n \sum_{x\in\annsmall_{\, i}} 
	\hskip-2pt
	\pr(0\connects x)^2
	\bigg[
	\sum_{j=1}^{i-1}  \sum_{y\in x+\annsmall_{\, j} } 
	\hskip-2pt
	J(x, y)
	+
	\hskip-2pt
	\sum_{y\in \annsmall_{\, i} \cap (x+\,\annsmall_{\, i} ) } 
	\hskip-2pt
	J(x, y)
	\bigg]
	\lesssim
	\sum_{i=1}^n \sum_{x\in\annsmall_{\, i}} 
	\hskip-2pt
	\pr(0\connects x)^2  L^{ i(d-\alpha) }
	\, .
\end{align*}
We can now use \eqref{eqn:788} to bound $F_1$ and $F_2$ via a direct calculation.
The claimed bound on $\Delta_n$ follows from \eqref{eqn:567} and an application of Lemma~\ref{lem:HL-zeta-n-m-j-n-asymptotics}(a).

While dealing with $\tilde\Delta_n$, we focus on proving the claimed bound when $0<\alpha\leq d/2$, as the other cases are similar.
When $\alpha\in (0, d/3)$ and $2\alpha\leq\theta<\alpha + d/3$,
\[
\tilde\Delta_n
\lesssim
L^{-n( d + \alpha - \theta )}\sum_{w\in\Lambda_n\setminus\{0\}} \pr(0\edge w)
=
L^{-n( d + \alpha - \theta )} m_n^{\sss (n)}
\leq
L^{-n( d + \alpha - \theta )} \, ,
\]
where the first step uses \eqref{eqn:789}, and the last step uses the relation $m_n^{\sss (n)}\leq 1$ from \eqref{eqn:m_n-leq-1}.
To get the desired bound when $0<\alpha\leq d/2$ and $\alpha<\theta<\min\{ 2\alpha,\, \alpha + d/3\}$, note that
\begin{align*}
\zeta_n \widetilde\Delta_n 
&
\asymp 
\sum_{w\in\Lambda_n\setminus\{0\}} \pr( 0\connects w )J(0, w)
\\
&
\lesssim
\sum_{1\leq i\leq a_0 n} L^{id} \cdot \zeta_n^{-1} L^{-i\alpha}\cdot L^{-i\alpha}
+
\sum_{a_0 n< i\leq n} L^{id} \cdot L^{-n( d + \alpha - \theta ) }\cdot L^{-i\alpha}\, ,
\end{align*}
where the last step uses \eqref{eqn:788} and \eqref{eqn:1a}.
The claimed bound on $\widetilde\Delta_n$ now follows from a direct calculation, and Lemma~\ref{lem:HL-zeta-n-m-j-n-asymptotics}(a).
\qed

\vskip5pt

Recall the coupling between the BFS tree of $\cC^-(x)$ and $\cT_n^{\sss (n)}$ explained right before \eqref{eqn:coupling}.
We now prove a result on how good the approximation of $\cC^-(0)$ via this branching process is.
This estimate will also be useful in later sections.

\begin{lem}\label{lem:3}
Consider $\alpha\in(0, d)$ and $\theta\in(\alpha, \alpha + d/3)$.
Then
\begin{align}\label{eqn:43}
0
\leq 
\bE\big( |\cT_n^{\sss(n)}| \big) - \bE\big( | \cC^-(0) | \big)
\lesssim
|\Lambda_n|^{2(\theta - \alpha)/d} \cdot \max\{ \Delta_n,\, \tilde\Delta_n \}\, .
\end{align}
Consequently,
\begin{align}\label{eqn:43-a}
0
\leq 
\bE\big( |\cT_n^{\sss(n)}| \big) - \bE\big( | \cC^-(0) | \big)
\lesssim
|\Lambda_n|^{-(d + 4\alpha - 4\theta)/d} 
\end{align}
if either 
(i) $\alpha\in (0, d/2]$ and $\theta\in(\alpha, \alpha + d/3)$, or
(ii) $ \alpha\in(d/2, 5d/6)$ and $2\alpha - d/2 < \theta <\alpha + d/3$.
On the other hand,
\begin{align}\label{eqn:43-b}
0
\leq 
\bE\big( |\cT_n^{\sss(n)}| \big) - \bE\big( | \cC^-(0) | \big)
\lesssim
|\Lambda_n|^{-2(d - \theta)/d} 
\end{align}
if either 
(i) $\alpha\in(d/2, 5d/6)$ and $\theta\in(\alpha, 2\alpha - d/2]$, or
(ii) $\alpha\in[5d/6, d)$ and $\theta\in(\alpha, 2d - \alpha)$.
In particular,
\begin{align}\label{eqn:42}
\bE|\cC^-(0)|
\gtrsim
|\Lambda_n|^{(\theta - \alpha)/d}
\end{align}
if either 
(i) $\alpha\in (0, 5d/6]$ and $\theta\in(\alpha, \alpha + d/3)$, or
(ii) $\alpha\in(5d/6, d)$ and $\theta\in(\alpha, 2d - \alpha)$.
\end{lem}

\noindent{\bf Proof:}
It is clear from \eqref{eqn:coupling} that 
$\bE\big( | \cC^-(0) | \big) \leq \bE\big( |\cT_n^{\sss(n)}| \big)$.
Now, in the coupling described above \eqref{eqn:coupling}, 
$ \big( 
|\cT_n^{\sss(n)}| - | \cC^-(0) | 
\big)$ 
is stochastically dominated by
\begin{equation}\label{eqn:88}
	\sum_{ y\in\cC^-(0)\setminus\{0\} } \
	\sum_{ z\in\cC^-(0)\setminus\{y\} }
	Z_{(y, z)} \cdot |\cT_{n; (y, z)}^{\sss(n)}|\, ,
\end{equation}
where conditional on $\cC^-(0)$, 
the random variables appearing in \eqref{eqn:88} are independent, 
$Z_{(y,z)}\sim\text{Bernoulli}(\fp^{\sss (n), -}_{\|y-z\|})$, and
$\cT_{n; (y, z)}^{\sss(n)} \equald \cT_n^{\sss (n)}$.
Thus, \eqref{eqn:88} together with an application of Lemma~\ref{lem:HL-zeta-n-m-j-n-asymptotics}(b) shows that
\begin{align}\label{eqn:87}
	&
	\bE\big( |\cT_n^{\sss(n)}| \big) - \bE\big( | \cC^-(0) | \big)
	\leq 
	\sum_{y\in\Lambda_n\setminus\{0\}}\,
	\sum_{z\in\Lambda_n\setminus\{y\}}
	\pr\big( 0\connects y,\, 0\connects z\big)
	\cdot \pr( y\edge z )
	\cdot (1 - m_n^{\sss(n)} )^{-1}
	\notag\\
	&
	\hskip40pt
	\lesssim
	L^{n(\theta - \alpha)}
	\sum_{y\in\Lambda_n\setminus\{0\}}\,
	\sum_{z\in\Lambda_n\setminus\{y\}}
	\pr\big( 0\connects y,\, 0\connects z\big)
	\cdot \pr( y\edge z )
	\leq
	L^{n(\theta - \alpha)} \big( F_1 + F_2 \big)\, ,
\end{align}
where
\begin{align*}
	&
	\hskip35pt
	F_1
	=
	\sum_{y\in\Lambda_n\setminus\{0\}}\,
	\sum_{z\in\Lambda_n\setminus\{y\}}
	\pr\big( \{ 0\connects z \}\disjt\{ z\connects y \} \big)
	\cdot \pr( y\edge z )
	\, ,\ \ \text{ and}
	\\
	&
	F_2
	=
	\sum_{y\in\Lambda_n\setminus\{0\}}\,
	\sum_{z\in\Lambda_n\setminus\{0, y\}}\,
	\sum_{x\in\Lambda_n\setminus\{y, z\}}
	\pr\big( \{ 0\connects x \}\disjt\{ x\connects y \}\disjt\{ x\connects z \} \big)
	\cdot \pr( y\edge z ) \, .
\end{align*}
Recall the quantities $\Delta_n$ and $\widetilde\Delta_n$ from \eqref{eqn:HL-def-Delta}.
Using the BK inequality, we see that
\begin{align}\label{eqn:61}
	F_1
	&
	\leq
	\sum_{z\in\Lambda_n} \pr(0\connects z) 
	\sum_{y\in\Lambda_n\setminus\{z\} } \pr(z\connects y)\cdot \pr( y\edge z )
	=
	\sum_{z\in\Lambda_n} \pr(0\connects z) \cdot \widetilde\Delta_n
	\notag\\
	&\hskip40pt
	=
	\bE\big( |\cC^-(0)| \big)\cdot \widetilde\Delta_n
	\leq
	\bE\big( |\cT_n^{\sss (n)}| \big)\cdot \widetilde\Delta_n
	=
	(1 - m_n^{\sss(n)})^{-1} \cdot \widetilde\Delta_n\, ,
\end{align}
and similarly,
\begin{align}\label{eqn:61-a}
	F_2
	&
	\leq
	\sum_{x\in\Lambda_n} \pr(0\connects x)
	\sum_{y\in\Lambda_n\setminus\{x\}} \,
	\sum_{z\in\Lambda_n\setminus\{x, y\}}
	\pr(x\connects y) \pr(x\connects z) \cdot \pr( y\edge z )
	\notag\\
	&
	\hskip30pt
	=
	\sum_{x\in\Lambda_n} \pr(0\connects x) \cdot\Delta_n
	\leq 
	(1 - m_n^{\sss(n)})^{-1} \cdot \Delta_n\, .
\end{align}
Combining \eqref{eqn:87}, \eqref{eqn:61}, and \eqref{eqn:61-a} together with an application of Lemma~\ref{lem:HL-zeta-n-m-j-n-asymptotics}(b) yields \eqref{eqn:43}.
Now \eqref{eqn:43-a} and \eqref{eqn:43-b} follow from\eqref{eqn:43} and Lemma~\ref{lem:2} via a direct calculation.
Finally, we note that 
$\bE|\cT_n^{\sss (n)}| 
=
( 1 - m_n^{\sss(n)} )^{-1}
\asymp 
L^{n(\theta - \alpha)}
$,
whereas the upper bounds on the right side of \eqref{eqn:43-a} and \eqref{eqn:43-b} are $o(L^{n(\theta - \alpha)})$ when $(\alpha, \theta)$ lies the regions specified there.
This implies \eqref{eqn:42}, and thus completes the proof.
\qed

\vskip5pt

\noindent{\bf Proof of Proposition~\ref{prop:HL-two-point-lower-bound}:}
Fix $(\alpha, \theta)$ such that either $\alpha\in (0, 5d/6]$ and $\theta\in (\alpha, \alpha + d/3)$, or 
$\alpha\in (5d/6, d)$ and $\theta\in ( \alpha, 2d - \alpha )$.
It follows from \eqref{eqn:1a}, \eqref{eqn:HL-def-kappa-minus}, Lemma~\ref{lem:HL-zeta-n-m-j-n-asymptotics}(a), and the relation $\alpha<\theta$ that there exists $C_{\ref{eqn:64}} = C_{\ref{eqn:64}} (\pmtr^{\ast}) > 0$ and 
$n_{\ref{eqn:64}} = n_{\ref{eqn:64}} (\pmtr)\geq 1$ such that
\begin{align}\label{eqn:64}
\min_{y\neq z\in\Lambda_n} \kappa^{\sss (n), -}(y, z)
\geq
C_{\ref{eqn:64}} L^{-nd}
\end{align}
for all $n\geq n_{\ref{eqn:64}}$.
Hence, for $n\geq n_{\ref{eqn:64}}$ and $x\in\Lambda_n\setminus\{0\}$,
\begin{align*}
&
\pr\big( 0\connects x\, \big|\, G_{\sss\Lambda_n\setminus\{x\} }^- \big)
=
1 - \prod_{ y\in\cC( 0 ;\, G_{\sss\Lambda_n\setminus\{x\} }^- ) }
\exp\big(
-\kappa^{\sss (n), -}(x, y)
\big)
\\
&
\hskip50pt
\geq
1 - \exp\big(  
- C_{\ref{eqn:64}} L^{-nd} \cdot |\cC( 0 ;\, G_{\sss\Lambda_n\setminus\{x\} }^- ) |
\big)
\geq
C_{\ref{eqn:64} } e^{-C_{\ref{eqn:64} } } 
L^{-nd} \cdot |\cC( 0 ;\, G_{\sss\Lambda_n\setminus\{x\} }^- ) |\, ,
\end{align*}
where the second step uses \eqref{eqn:64}, and the last step uses the observation that 
$ |\cC( 0 ;\, G_{\sss\Lambda_n\setminus\{x\} }^- ) | < |\Lambda_n| = L^{nd}$.
Taking expectation on both sides, we get
\begin{align}\label{eqn:64-a}
\pr( 0\connects x)
\gtrsim
L^{-nd} \cdot \bE  |\cC( 0 ;\, G_{\sss\Lambda_n\setminus\{x\} }^- ) | 
\ \ \text{ for }\ \ x\in\Lambda_n\setminus\{0\}
\, .
\end{align}

Next, for any $x\neq y\in\Lambda_n\setminus\{0\}$, 
\begin{align}\label{eqn:64-b}
\pr(0\connects y) 
&
\leq
\pr\big(0\connects y\text{ in }\Lambda_n\setminus\{x\} \big)
+
\pr\big( \{0\connects x\}\disjt\{ x\connects y \} \big)
\notag\\
&
\leq
\pr\big(0\connects y\text{ in }\Lambda_n\setminus\{x\} \big)
+
\pr(0\connects x) \cdot\pr( x\connects y )\, ,
\end{align}
where the second step uses the BK inequality.
Hence, for $x\in\Lambda_n\setminus\{0\}$,
\begin{align}\label{eqn:64-c}
&
\bE |\cC^-(0)|
=
1 + \pr(0\connects x) + \sum_{y\in\Lambda_n\setminus\{0, x\}} \pr( 0 \connects y)
\notag\\
&\hskip40pt
\leq
1 
+ 
\sum_{y\in\Lambda_n\setminus\{0, x\}} 
\pr\big( 0 \connects y\text{ in }\Lambda_n\setminus\{x\} \big)
+
\pr(0\connects x) \cdot\big[
1 + \sum_{y\in\Lambda_n\setminus\{x\}} \pr( x\connects y )
\big]
\notag\\
&\hskip80pt
=
\bE  |\cC( 0 ;\, G_{\sss\Lambda_n\setminus\{x\} }^- ) | 
+
\pr(0\connects x) \cdot\bE|\cC^-(0)| \, ,
\end{align}
where the second step uses \eqref{eqn:64-b}, 
and the third step uses the relation $\bE|\cC^-(x)| = \bE|\cC^-(0)|$. 
Now, \eqref{eqn:788} and \eqref{eqn:789} show that 
$\max_{x\in\Lambda_n\setminus\{0\}} \pr(0\connects x) = o(1)$, which combined with \eqref{eqn:64-c} yields
$\bE  |\cC( 0 ;\, G_{\sss\Lambda_n\setminus\{x\} }^- ) | 
\gtrsim
\bE|\cC^-(0)|
$
for $x\in\Lambda_n\setminus\{0\}$.
Plugging this into \eqref{eqn:64-a}, we see that for $x\in\Lambda_n\setminus\{0\}$,
\[
\pr( 0\connects x)
\gtrsim
L^{-nd}\cdot\bE |\cC^-(0)|
\gtrsim
L^{-n( d +\alpha - \theta )} \, ,
\]
where the last step follows from \eqref{eqn:42}.
This completes the proof.
\qed

\subsection{Tail bounds for $x_{\max}$ and $\diamax$}\label{sec:HL-proof-xmax-diamax}
We now start working toward the proofs of Theorems~\ref{thm:HL-scaling-limit} and \ref{thm:HL-surplus-convergence}.
So, we now fix a $\theta$ satisfying
\begin{align}\label{eqn:HL-def-theta}
	\alpha < \theta < 
	\begin{cases}
		4\alpha/3\, , \ \ \text{ if }\ 0<\alpha\leq d/2\, ,\\
		2\alpha - d/2\, , \ \ \text{ if }\ d/2<\alpha\leq 2d/3\, ,\\
		5\alpha/2 - d\, , \ \ \text{ if }\ 2d/3<\alpha< 5d/6\, .
	\end{cases}
\end{align}
Throughout Section~\ref{sec:HL-proof-xmax-diamax}, we work with $\alpha\in(0, 5d/6)$ and $\theta$ as in \eqref{eqn:HL-def-theta}.
Recall that $x_{\max}$ and $\diamax$ were introduced right below \eqref{eqn:510}.
The aim of this section is to prove the following tail bounds:

\begin{prop}\label{prop:HL-xmax-diamax-bound}
	For each $r>0$, there exists $C_{\ref{eqn:767}}(r)>0$ depending only on $r$ and $\pmtr^{\ast}$ such that
	\begin{gather}
		\pr\big(\diamax \geq C_{\ref{eqn:767}}(r) \cdot n L^{n(\theta - \alpha)}\big)
		\lesssim e^{-nr}
		\, , \ \ \text{ and}
		\label{eqn:767}\\
		\pr\big(x_{\max} \geq C_{\ref{eqn:767}}(r) \cdot n L^{2n(\theta - \alpha - d/3)}\big)
		\lesssim e^{-nr}\, .
		\label{eqn:768}
	\end{gather}
\end{prop}

\noindent{\bf Proof:}
It follows from Lemma~\ref{lem:HL-zeta-n-m-j-n-asymptotics}(b) that there exist 
$C_{\ref{eqn:898}}>0$ and $n_{\ref{eqn:898}} \geq 1$ depending only on $\pmtr^{\ast}$ and $\pmtr$ respectively such that 
\begin{align}\label{eqn:898}
m_n^{\sss(n)}
\leq
1-C_{\ref{eqn:898}} L^{-n(\theta - \alpha)}
\ \ \text{ for all } n\geq n_{\ref{eqn:898}} .
\end{align}
Let $\cW_k^{\sss (n)}$, $k\geq 0$, be as defined at the beginning of Section~\ref{sec:HL-proof-branching-process}.
Then
\begin{align}\label{eqn:9}
	\pr\big(\height(\cT_n^{\sss (n)}) \geq k\big)
	&
	=
	\pr(\cW_k^{\sss (n)} \geq 1)
	\leq
	\bE(\cW_k^{\sss (n)})
	=
	(m_n^{\sss (n)})^k
	\notag\\
	&
	\leq
	\big( 1- C_{\ref{eqn:898}} L^{-n(\theta - \alpha)}\big)^k
	\leq 
	\exp\big( - C_{\ref{eqn:898}} k L^{-n(\theta - \alpha)} \big)
\end{align}
for any $n\geq n_{\ref{eqn:898}}$ and $k\geq 1$, where the penultimate step uses \eqref{eqn:898}.
Consequently,
\begin{align*}
	\pr\big( \diamax\geq 2k \big)
	&
	\leq
	\sum_{x\in\Lambda_n}\pr\big(\diam(\cC^-(x)) \geq 2k\big)
	\\
	&
	\leq
	|\Lambda_n|\cdot \pr\big(\height(\cT_n^{\sss (n)}) \geq k\big)
	\leq 
	L^{nd}\cdot \exp\big( - C_{\ref{eqn:898}} k L^{-n(\theta - \alpha)} \big)\, ,
\end{align*}
where the second step uses \eqref{eqn:coupling}.
The claim in \eqref{eqn:767} now follows easily.

Turning to the proof of \eqref{eqn:768}, we define
\begin{align*}
	\varphi_n(t)
	:=
	\bE\big[ \exp({t\cW_1^{\sss (n)}}) \big]
	=
	\prod_{i=1}^{n} \big( 1 - \fp_{L^i}^{\sss (n), -} + e^t \fp_{L^i}^{\sss (n), -}  \big)^{L^{id} - L^{(i-1)d}}
	\, ,\ \ t\in\bR\, .
\end{align*}
Then for all $n$,
\begin{align}\label{eqn:10}
	\varphi_n(1)
	\leq
	\prod_{i=1}^{n} \big( 1 + 2 \fp_{L^i}^{\sss (n), -}  \big)^{L^{id} - L^{(i-1)d}}
	\leq
	\exp\bigg(\sum_{i=1}^{n} 2 \fp_{L^i}^{\sss (n), -} (L^{id} - L^{(i-1)d})\bigg)
	=
	\exp(2m_n^{\sss (n)})
	\leq 
	e^2\, ,
\end{align}
where the last step uses the relation $m_n^{\sss (n)} \leq 1$.
In particular, for any $n, k\geq 1$,
\begin{align}\label{eqn:10-a}
	\bE\big[ (\cW_1^{\sss (n)} )^k\big]
	\leq 
	k!\cdot\varphi_n(1)
	\leq
	k!\cdot e^2\, .
\end{align}
Now, for any $t\in [0, 1/2]$ and $n\geq n_{\ref{eqn:898}}$,
\begin{align}
	\varphi_n(t)
	&
	=
	1 + t\cdot \bE[\cW_1^{\sss (n)}] + \sum_{j\geq 2}  \frac{t^j}{j!} \cdot \bE\big[ (\cW_1^{\sss (n)} )^j \big] 
	\leq
	1 + t m_n^{\sss (n)} + t^2 \bE\big[ (\cW_1^{\sss (n)} )^2 \exp( t\cW_1^{\sss (n)}  ) \big]
	\notag\\
	&
	\leq
	1 + t m_n^{\sss (n)} + 
	t^2 \big( \bE\big[ (\cW_1^{\sss (n)} )^4 \big] \big)^{1/2} \varphi_n(2t)^{1/2}
	\leq
	1 + t \big( 1 - C_{\ref{eqn:898}} L^{-n(\theta - \alpha)} \big) + C_{\ref{eqn:12}} t^2\, ,
	\label{eqn:12}
\end{align}
where the last step uses \eqref{eqn:898}, \eqref{eqn:10}, \eqref{eqn:10-a}, and the relation $t\leq 1/2$.
Let 
\begin{align}\label{eqn:5}
	\eps_n
	:=
	\min\bigg\{
	\frac{C_{\ref{eqn:898}} }{2 C_{\ref{eqn:12}} L^{n(\theta - \alpha)}}
	\, , \,
	\frac{1}{2}
	\bigg\}\, .
\end{align}
Then \eqref{eqn:12} shows that for all $n\geq n_{\ref{eqn:898}}$,
\begin{align}\label{eqn:13}
	\varphi_n(\eps_n)
	\leq 
	1 + \eps_n \big( 1 - C_{\ref{eqn:13}} L^{-n(\theta - \alpha)} \big)
	\leq
	\exp\big( \eps_n \big( 1 - C_{\ref{eqn:13}} L^{-n(\theta - \alpha)} \big) \big)\, .
\end{align}
Let $V_i^{\sss (n)}$, $i\geq 1$, be an i.i.d. sequence with $V_1^{\sss (n)}\equald\cW_1^{\sss(n)}$.
Then for any $n\geq n_{\ref{eqn:898}}$ and $k\geq 1$,
\begin{align}\label{eqn:14}
	\pr(|\cT_n^{\sss (n)}| > k)
	&
	= 
	\pr\big(\min_{1\leq j\leq k} \sum_{i=1}^j (V_i^{\sss (n)} - 1) > -1 \big)
	\leq
	\pr\big( \sum_{i=1}^k V_i^{\sss (n)} \geq k \big)
	\notag\\
	&
	\leq
	\exp(-\eps_n k) \cdot\varphi_n(\eps_n)^k
	\leq
	\exp\big( - k\eps_n C_{\ref{eqn:13}} L^{-n(\theta - \alpha )}\big)\, ,
\end{align}
where the last step uses \eqref{eqn:13}.
Now the claimed bound in \eqref{eqn:768} follows upon noting that for any $a>0$,
\[
\pr(x_{\max} \geq a )
=
\pr\big( |\cC_1^-|\geq L^{2nd/3} a\big)
\leq
L^{nd} \pr\big( |\cC^-(0)|\geq L^{2nd/3} a\big)
\leq
L^{nd} \pr\big( |\cT_n^{\sss (n)}|\geq L^{2nd/3} a\big)
\, ,
\]
and using \eqref{eqn:14} and \eqref{eqn:5}.
\qed

\subsection{Asymptotics for $\sigma_2$}\label{sec:HL-proof-sigma-2}
Throughout Section~\ref{sec:HL-proof-sigma-2}, we work with $\alpha\in(0, 5d/6)$ and $\theta$ as in \eqref{eqn:HL-def-theta}.
Recall $q$ and $\sigma_2$ from \eqref{eqn:510} and from right below \eqref{eqn:510}.
The aim of this section is to establish the following convergence which is the most delicate one among all the convergences in Conditions~\ref{ass:aldous-basic-assumption}, \ref{ass:gen-1}, and \ref{ass:gen-2}.

\begin{prop}\label{prop:HL-sigma-2-asymptotics}
	We have, $q -( \sigma_2 )^{-1} \weakc \lambda $ as $n\to\infty$.
\end{prop}

It will be convenient to work with the susceptibility instead of $\sigma_2$.
For $k\geq 1$, we define the $k$-th susceptibility as
\begin{align}\label{eqn:HL-def-susceptibility}
	\bars_{k+1}^{\sss(n), -}
	:=
	L^{-nd}\sum_i |\cC_i^-|^{k+1}
	=
	L^{-nd}\sum_{x\in\Lambda_n} | \cC^-(x) |^k
	=
	\bE\big(\, |\cC^-(U)|^k \, \big|\, G_{\sss \Lambda_n}^- \big)\, ,
\end{align}
where $U$ is uniformly distributed over $\Lambda_n$ independent of $G_{\sss \Lambda_n}^-$.
Clearly,
\begin{align}\label{eqn:33}
	\sigma_k = L^{ -nd(2k-3)/3  } \cdot\bars_k^{\sss (n), -}\, ,\ \ k=2,3\, .
\end{align}
It follows from \eqref{eqn:HL-def-susceptibility} that
$
\bE \bars_2^{\sss(n), -} = \bE|\cC^-(U)| 
$,
and Lemma~\ref{lem:3} gives us control over $\bE|\cC^-(0)| = \bE|\cC^-(U)| $.
We now need a fluctuation bound on $\bars_2^{\sss(n), -} $.

\begin{lem}\label{lem:4}
	We have,
	$\var(\bars_k^{\sss (n), - }) 
	\lesssim
	L^{-nd + n (4k - 3) (\theta - \alpha) }$, 
	$k=2,3$.
\end{lem}

\noindent{\bf Proof:}
Let  $U$ be uniformly distributed over $\Lambda_n$ independent of $G_{\sss \Lambda_n}^-$.
Then for any $\ell\geq 1$,
\begin{align}\label{eqn:98}
	\bE\big( |\cC^-(0)|^\ell\cdot |\cC^-(U)|^\ell \big)
	&
	=
	\bE\big[\bE\big( |\cC^-(0)|^\ell\cdot |\cC^-(U)|^\ell\, \big|\, G_{\sss \Lambda_n}^- \big)\big]
	=
	\bE\big( |\cC^-(0)|^\ell\cdot\bars_{\ell+1}^{\sss (n), -} \big)
	\notag\\
	&
	=
	\bE\big( L^{-nd} \sum_{x\in\Lambda_n}|\cC^-(x)|^\ell \cdot \bars_{\ell+1}^{\sss (n), -} \big)
	=
	\bE\big[ \big(\bars_{\ell+1}^{\sss (n), -} \big)^2 \big]\, ,
\end{align}
where the second step uses \eqref{eqn:HL-def-susceptibility}, and the third step uses the fact that for each $x\in\Lambda_n$,
$\big( \cC^-(0) ,\, \bars_{\ell+1}^{\sss (n), -} \big)
\equald
\big( \cC^-(x) ,\, \bars_{\ell+1}^{\sss (n), -} \big)$.
For $0\in A\subseteq\Lambda_n$, let $\cE_A$ denote the event that $A$ is the vertex set of $\cC^-(0)$.
Then
\begin{align*}
	\bE\big( |\cC^-(U)|^{\ell}\, \big|\, \cE_A \big)
	=
	L^{-nd}\bigg[
	|A|^{\ell + 1} 
	+ 
	\bE\big( 
	\sum_{x\in\Lambda_n\setminus A} |\cC(x \, ;\, G_{\sss \Lambda_n\setminus A}^- )|^{\ell} 
	\big)
	\bigg]
	\leq 
	L^{-nd}|\cC^-(0)|^{\ell + 1} 
	+
	\bE\big( \bars_{\ell+1}^{\sss (n), -} \big)\, ,
\end{align*}
which in turn implies that
\begin{align}\label{eqn:45}
	\bE\big( |\cC^-(0)|^\ell\cdot |\cC^-(U)|^\ell \big)
	\leq
	L^{-nd}\cdot\bE\big( |\cC^-(0)|^{2\ell + 1} \big)
	+
	\bE\big( |\cC^-(0)|^{\ell} \big)\cdot\bE\big( \bars_{\ell+1}^{\sss (n), -} \big)\, .
\end{align}
By symmetry, 
$\bE\big( |\cC^-(0)|^{\ell} \big) = \bE\big( |\cC^-(x)|^{\ell} \big)$ for all $x\in\Lambda_n$, 
and consequently, 
\begin{align}\label{eqn:46} 
	\bE\big( |\cC^-(0)|^{\ell} \big)
	= 
	L^{-nd} \cdot \bE\big( \sum_{x\in\Lambda_n}  |\cC^-(x)|^{\ell} \big) 
	=
	\bE\big( \bars_{\ell+1}^{\sss (n), -} \big)\, .
\end{align}
We thus get from \eqref{eqn:98}, \eqref{eqn:45}, and \eqref{eqn:46} that
\begin{align}\label{eqn:99}
	\var\big(\bars_{\ell+1}^{\sss (n), -} \big)
	\leq 
	L^{-nd}\cdot\bE\big( |\cC^-(0)|^{2\ell + 1} \big)
	\leq
	L^{-nd}\cdot\bE\big( |\cT_n^{\sss (n)} |^{2\ell + 1} \big)\, ,
\end{align}
where the last step uses \eqref{eqn:coupling}.

We will now show that for all $k\geq 1$ there exists $C_{\ref{eqn:111}}(k)>0$ depending only on $k$ and $\pmtr^{\ast}$ such that 
\begin{align}\label{eqn:111}
	\bE(|\cT_n^{\sss (n)} |^{k})
	\leq
	C_{\ref{eqn:111}}(k) \cdot L^{n(\theta - \alpha) (2k -1)}
	\ \ \text{ for all } n\geq n_{\ref{eqn:898}}\, ,
\end{align}
which combined with \eqref{eqn:99} will complete the proof of Lemma~\ref{lem:4}.
We prove this claim via induction on $k$.
First note that \eqref{eqn:898} implies that \eqref{eqn:111} holds for $k=1$.
Let us assume that \eqref{eqn:111} holds for $k=1,\ldots,\ell -1$ for some $\ell\geq 2$.
Now,
\begin{align}\label{eqn:55}
	|\cT_n^{\sss (n)}| 
	\equald
	1 + \sum_{i=1}^{\cW_1^{\sss (n)}} |\cT_{n, i}^{\sss (n)}| \, ,
\end{align}
where $( \cT_{n, i}^{\sss (n)}\, ;\, i\geq 1 )$ is an i.i.d. sequence independent of $\cW_1^{\sss (n)}$ with 
$\cT_{n, i}^{\sss (n)} \equald \cT_n^{\sss (n)}$.
Hence,
\begin{align*}
	&\bE\big[ \big( |\cT_n^{\sss (n)}| - 1\big)^{\ell}\big]
	=
	\bE\big[ \big( \sum_{i=1}^{\cW_1^{\sss (n)}} |\cT_{n, i}^{\sss (n)}| \big)^{\ell}\big]
	\\
	&\hskip40pt
	\leq 
	\bE(\cW_1^{\sss (n)}) \cdot \bE( |\cT_n^{\sss (n)}|^{\ell} )
	+
	\sum_{j=2}^{\ell} \frac{1}{j!} \bE\big[ (\cW_1^{\sss (n)})^j \big]\cdot
	\sum\displaystyle_{\dagger} 
	\frac{\ell !}{k_1! \cdots k_j !} \bE\big[ |\cT_{n, 1}|^{k_1} \ldots |\cT_{n, j}|^{k_j} \big] \, ,
\end{align*}
where $\sum_{\dagger}$ denotes the sum over all $k_1, \ldots, k_j\geq 1$ satisfying $k_1 +\ldots +k_j=\ell$.
It thus follows from \eqref{eqn:10-a} and the induction hypothesis that for $n\geq n_{\ref{eqn:898}}$,
\begin{align}
	\bE\big[ \big( |\cT_n^{\sss (n)}| - 1\big)^{\ell}\big]
	&
	\leq 
	m_n^{\sss (n)}\cdot \bE( |\cT_n^{\sss (n)}|^{\ell} )
	+
	\sum_{j=2}^{\ell} e^2\cdot
	\sum\displaystyle_{\dagger} 
	\frac{\ell !}{k_1! \cdots k_j !} \big(L^{n(\theta-\alpha)}\big)^{2\ell - j} \prod_{i=1}^j C_{\ref{eqn:111}}(k_i) 
	\label{eqn:56}
	\\
	&
	\leq
	m_n^{\sss (n)}\cdot \bE( |\cT_n^{\sss (n)}|^{\ell} )
	+
	\big(L^{n(\theta-\alpha)}\big)^{2\ell - 2}
	\sum_{j=2}^{\ell} e^2\cdot
	\sum\displaystyle_{\dagger} 
	\frac{\ell !}{k_1! \cdots k_j !} \prod_{i=1}^j C_{\ref{eqn:111}}(k_i) \,. \notag
\end{align}
Expanding the term $\big( |\cT_n^{\sss (n)}| - 1\big)^{\ell}$ appearing on the left side of \eqref{eqn:56} and using \eqref{eqn:898} completes the proof of \eqref{eqn:111}.
\qed

\vskip6pt

\noindent{\bf Proof of Proposition~\ref{prop:HL-sigma-2-asymptotics}:}
Using \eqref{eqn:33}, we see that
\begin{align}\label{eqn:686}
	q-\frac{1}{\sigma_2}
	=
	\bigg( q-\frac{L^{nd/3}}{\bE |\cT_n^{\sss (n)}| } \bigg)
	+
	\bigg( \frac{L^{nd/3}}{\bE |\cT_n^{\sss (n)}| } - \frac{L^{nd/3}}{\bE\bars_2^{\sss (n), -} }\bigg)
	+
	\bigg( \frac{L^{nd/3}}{\bE\bars_2^{\sss (n), -} } - \frac{L^{nd/3}}{\bars_2^{\sss (n), -} }\bigg)
	=:
	F_1 + F_2 + F_3\, .
\end{align}
We bound the terms appearing on the right side of \eqref{eqn:686} one by one. 
First,
\begin{align}\label{eqn:686-a}
	F_1
	=
	q - L^{nd/3}\big(1 - m_n^{\sss (n)}\big)
	&
	=
	\lambda + \zeta_n^{-1} L^{n(\frac{4d}{3} - \theta)} 
	- 
	\big(1 + O(L^{-nd})\big)\zeta_n^{-1} L^{n(\frac{4d}{3} - \theta)} 
	\notag\\
	&
	=
	\lambda + O\bigg( \frac{L^{n(\frac{d}{3} - \theta)}}{L^{n(d-\alpha)}} \bigg)
	=
	\lambda + o(1)\, ,
\end{align}
where the second step uses Lemma~\ref{lem:HL-zeta-n-m-j-n-asymptotics}(b) and the definition of $q$ from \eqref{eqn:510}, and the third step uses Lemma~\ref{lem:HL-zeta-n-m-j-n-asymptotics}(a).

Now, as observed in \eqref{eqn:46}, $\bE( \bars_2^{\sss (n), -} ) = \bE( |\cC^-(0)| )$.
Further, 
$\bE|\cT_n^{\sss (n)}| = (1-m_n^{\sss (n)})^{-1}\asymp L^{n(\theta-\alpha)}$,
where the last step uses Lemma~\ref{lem:HL-zeta-n-m-j-n-asymptotics}(b).
Hence, we see from Lemma~\ref{lem:3} and from the choice of $\theta$ in \eqref{eqn:HL-def-theta} that
$0\leq \bE|\cT_n^{\sss (n)}| - \bE\bars_2^{\sss (n), -} = o(\bE|\cT_n^{\sss (n)}| )$.
In partcular, 
$\bE\bars_2^{\sss (n), -}\asymp L^{n(\theta-\alpha)}$,
and consequently
\begin{align}\label{eqn:686-b}
	|F_2|
	\lesssim
	L^{nd/3}\cdot\frac{\bE|\cT_n^{\sss (n)}| - \bE|\cC^-(0)| }{L^{2n(\theta - \alpha)}}
	=
	o(1)\, ,
\end{align}
where the last step uses Lemma~\ref{lem:3}, \eqref{eqn:HL-def-theta}, and the fact that $\alpha<5d/6$.

Finally, 
$| \bars_2^{\sss (n), -} - \bE\bars_2^{\sss (n), -} |
=
O_P\big( \big( \var(\bars_2^{\sss (n), -}) \big)^{1/2} \big)
=
o_P( \bE\bars_2^{\sss (n), -})$,
where the last step uses Lemma~\ref{lem:4}, \eqref{eqn:HL-def-theta}, and the fact that 
$\bE\bars_2^{\sss (n), -}\asymp L^{n(\theta-\alpha)}$.
In particular, $ (\bars_2^{\sss (n), -} )^{-1} = O_P(L^{-n(\theta-\alpha)})$.
Hence,
\begin{align}\label{eqn:686-c}
	|F_3|
	\lesssim
	L^{nd/3}\cdot O_P\bigg( \frac{ \big( \var(\bars_2^{\sss (n), -}) \big)^{1/2}  }{L^{2n(\theta - \alpha)}} \bigg)
	=
	o_P(1)\, ,
\end{align}
where the last step uses Lemma~\ref{lem:4} and \eqref{eqn:HL-def-theta}.

Combining \eqref{eqn:686}, \eqref{eqn:686-a}, \eqref{eqn:686-b}, and \eqref{eqn:686-c} completes the proof of Proposition~\ref{prop:HL-sigma-2-asymptotics}.
\qed

\subsection{Asymptotics for $\sigma_3$}\label{sec:HL-proof-sigma-3}
Throughout Section~\ref{sec:HL-proof-sigma-3}, we work with $\alpha\in(0, 5d/6)$ and $\theta$ as in \eqref{eqn:HL-def-theta}.
Recall $\sigma_3, \sigma_2$ from right below \eqref{eqn:510}.
The aim of this section is to establish the following convergence:

\begin{prop}\label{prop:HL-sigma-3-asymptotics}
	We have, $( \sigma_2 )^{-3}\sigma_3 \weakc 1 $ as $n\to\infty$.
\end{prop}

\noindent{\bf Proof:}
Recall the definition of $\cW_1^{\sss(n)}$ from Section~\ref{sec:HL-proof-branching-process}.
Then
\begin{align}\label{eqn:35}
	\var(\cW_1^{\sss(n)})
	&
	=
	\sum_{i=1}^n (L^{id} - L^{(i-1)d})\cdot \fp_{L^i}^{\sss (n), -} (1 - \fp_{L^i}^{\sss (n), -})
	\notag\\
	&
	=
	m_n^{\sss (n)} - \sum_{i=1}^n (L^{id} - L^{(i-1)d})\cdot ( \fp_{L^i}^{\sss (n), -} )^2
	=
	m_n^{\sss (n)} - O\big( m_n^{\sss (n)}\cdot \max_{1\leq i\leq n}  \fp_{L^i}^{\sss (n), -} \big) \, .
\end{align}
Since 
$\max\limits_{1\leq i\leq n}  \fp_{L^i}^{\sss (n), -} \lesssim\zeta_n^{-1}$, Lemma~\ref{lem:HL-zeta-n-m-j-n-asymptotics}(a) and (b) combined with \eqref{eqn:35} yield
\begin{align}\label{eqn:36}
	\var(\cW_1^{\sss(n)}) \to 1 \ \ \text{ as }\ \ n\to\infty\, .
\end{align}
Now, it follows from \eqref{eqn:55} that 
\begin{align*}
	\var(|\cT_n^{\sss (n)}|)
	&
	=
	\var\big( \cW_1^{\sss (n)} \cdot \bE( |\cT_n^{\sss (n)}| ) \big)
	+
	\bE\big( \cW_1^{\sss (n)} \cdot \var( |\cT_n^{\sss (n)}| ) \big)
	\\
	&
	=
	(1-m_n^{\sss (n)})^{-2} \var( \cW_1^{\sss (n)} )
	+
	m_n^{\sss (n)} \cdot \var( |\cT_n^{\sss (n)}| ) \, ,
\end{align*}
where the second step uses the relation 
$\bE |\cT_n^{\sss (n)}|  = (1-m_n^{\sss (n)})^{-1}$.
We thus get
$
\var(|\cT_n^{\sss (n)}|) 
=
\var( \cW_1^{\sss (n)} )\cdot (1-m_n^{\sss (n)})^{-3} 
$,
and consequently
\begin{align}\label{eqn:37}
	\bE( |\cT_n^{\sss (n)}|^2 )
	&
	=
	\var( \cW_1^{\sss (n)} )\cdot (1-m_n^{\sss (n)})^{-3} 
	+ 
	\big( \bE |\cT_n^{\sss (n)}| \big)^2
	\notag
	\\
	&
	=
	( 1 + o(1) )\cdot (1-m_n^{\sss (n)})^{-3} 
	+
	(1-m_n^{\sss (n)})^{-2} 
	=
	\frac{\zeta_n^3}{L^{3n (d-\theta)}} ( 1 + o(1) )
	\asymp
	L^{3n(\theta - \alpha)}\, ,
\end{align}
where the second step uses \eqref{eqn:36}, the third step uses Lemma~\ref{lem:HL-zeta-n-m-j-n-asymptotics}(b), and the final step uses Lemma~\ref{lem:HL-zeta-n-m-j-n-asymptotics}(a).

Next, assume that $\cT_n^{\sss (n)}$ and $\cC^-(0)$ are coupled as described around \eqref{eqn:coupling}.
Then for any sequence $(\omega_n\, ;\, n\geq 1)$ of positive numbers,
\begin{align}\label{eqn:38}
	&
	\bE( |\cT_n^{\sss (n)}|^2 ) - \bE\bars_3^{\sss (n), -}
	=
	\bE\big[ \big( |\cT_n^{\sss (n)}| + |\cC^-(0)| \big)\cdot \big( |\cT_n^{\sss (n)}| - |\cC^-(0)| \big) \big]
	\leq
	2\bE\big[ |\cT_n^{\sss (n)}| \cdot \big( |\cT_n^{\sss (n)}| - |\cC^-(0)| \big) \big]
	\notag
	\\
	&
	\hskip50pt
	\leq
	2\omega_n\cdot \bE\big( |\cT_n^{\sss (n)}| - |\cC^-(0)| \big) 
	+
	2\bE\big( |\cT_n^{\sss (n)}|^2 \cdot \ind_{ \{ |\cT_n^{\sss (n)}| \geq \omega_n \} } \big)
	\notag
	\\
	&
	\hskip50pt
	=
	2\omega_n\cdot \bE\big( |\cT_n^{\sss (n)}| - |\cC^-(0)| \big) 
	+
	\big[
	\omega_n^2 \pr( |\cT_n^{\sss (n)}| \geq \omega_n)
	+
	\medint\int_{\omega_n}^{\infty} 2t\pr( |\cT_n^{\sss (n)}| \geq t)\, dt
	\big]\, ,
\end{align}
where the first step uses \eqref{eqn:46} with $\ell=2$, and the second step uses the relation 
$ |\cC^-(0)| \leq |\cT_n^{\sss (n)}| $ from \eqref{eqn:coupling}.
Now, \eqref{eqn:14} and \eqref{eqn:5} show that there exists $C>0$ such that with the choice 
$\omega_n = C n L^{2n(\theta-\alpha)}$ 
in \eqref{eqn:38}, we get
\begin{align*}
	\bE( |\cT_n^{\sss (n)}|^2 ) - \bE\bars_3^{\sss (n), -}
	\leq
	2 C n L^{2n(\theta-\alpha)}\cdot \bE\big( |\cT_n^{\sss (n)}| - |\cC^-(0)| \big) + o(1)
	=
	o( L^{3n(\theta - \alpha)} )\, ,
\end{align*}
where the last step uses Lemma~\ref{lem:3} and \eqref{eqn:HL-def-theta}.
We thus get from \eqref{eqn:37} that
\begin{align}\label{eqn:39}
	\bE\bars_3^{\sss (n), -}
	=
	\frac{\zeta_n^3}{L^{3n (d-\theta)}} ( 1 + o(1) )
	\asymp
	L^{3n(\theta - \alpha)}\, .
\end{align}
Now, Lemma~\ref{lem:4} and \eqref{eqn:HL-def-theta} imply that
$
\var(\bars_3^{\sss (n),- })
\lesssim
L^{-nd + 9n(\theta - \alpha)} 
=
o(L^{6n(\theta - \alpha)})
$,
which combined with \eqref{eqn:39} yields
\begin{align}\label{eqn:40}
	\bars_3^{\sss (n), -}
	=
	\frac{\zeta_n^3}{L^{3n (d-\theta)}} ( 1 + o_P(1) )\, .
\end{align}

Finally, Proposition~\ref{prop:HL-sigma-2-asymptotics} and the definition $q$ in \eqref{eqn:510} show that
\begin{align}\label{eqn:41}
	(\sigma_2)^{-1} 
	=
	(1+o_P(1)) \cdot \zeta_n^{-1} L^{n (\frac{4d}{3} - \theta ) } \, .
\end{align}
Further, as noted in \eqref{eqn:33}, 
$\sigma_3 = L^{-nd} \bars_3^{\sss (n), -}$,
which combined with \eqref{eqn:40} and \eqref{eqn:41} completes the proof.
\qed

\subsection{Asymptotics for $\tau$}\label{sec:HL-proof-tau}
Throughout Section~\ref{sec:HL-proof-tau}, we work with $\alpha\in(0, 5d/6)$ and $\theta$ as in \eqref{eqn:HL-def-theta}.
Recall $\tau$ from right below \eqref{eqn:510}.
The aim of this section is to prove the following result:

\begin{prop}\label{prop:HL-tau-asymptotics}
	We have, $\tau = \zeta_n^2\cdot L^{-n(\frac{7d}{3} - 2\theta)} \cdot (1+o_P(1))$ as $n\to\infty$.
\end{prop}

We introduce some notation here that will be needed in the proof of the above proposition.
For any graph $\cH$, we will write $d_{\cH}(\cdot, \cdot)$ for the graph distance in $\cH$.
We will simply write $d^-$ for $d_{G_{\sss\Lambda_n}^-}$.
For $x\in A\subseteq\Lambda_n$, we let
\begin{gather*}
	Z_k^-(x ; A) := 
	\#\big\{
	y\in \cC(x ; G_A^-)\, :\, d_{G_A^-}(x, y) = k
	\big\} , \ 
	k\geq 0\, , \ \text{ and}
	\\[2pt]
	\cD_A^-(x)
	:=
	\sum_{y\in\cC(x ; G_A^-)} d_{G_A^-} (x, y)
	=
	\sum_{k\geq 1} k Z_k^-(x ; A) 
	\, .
\end{gather*}
We will simply write $\cD^-(x)$ for $\cD_{\sss\Lambda_n}^-(x)$.

For $x\in A\subseteq\Lambda_n$, let $\cT^{\sss (n)}(x; A)$ denote the following multitype branching process tree with type space $A$:
$\cT^{\sss (n)}(x; A)$ is rooted at a vertex of type $x$, and each vertex in $\cT^{\sss (n)}(x; A)$ of type $y\in A$ has one child of type $z$ with probability $\fp_{\| y-z \|}^{\sss (n), -}$ for each $z\in A\setminus\{y\}$, where $\fp_{\cdot}^{\sss (n), -}$ is as in \eqref{eqn:HL-def-fp}.
Let $\cW_k^{\sss (n)} (x ; A)$ denote the number of vertices in the $k$-th generation of $\cT^{\sss (n)}(x; A)$, $k\geq 0$.
(Thus, for any $1\leq j\leq n$ and $x\in\Lambda_j$, the tree obtained from $\cT^{\sss (n)}(x; \Lambda_j)$ by omitting the types of its vertices has the same distribution as $\cT_j^{\sss (n)}$.
Further,  
$\cW_k^{\sss (n)} (x ; \Lambda_n)
\equald
\cW_k^{\sss (n)}$, $k\geq 0$.)
By an argument similar to the one leading to \eqref{eqn:coupling}, the BFS tree of $\cC(x; G_A^-)$ (with the breadth-first exploration of $G_A^-$ starting at $x$) can be coupled with $\cT^{\sss (n)}(x; A)$ so that both trees have the same root and the former is a subtree of the latter.
Note that for each $y\in \cC(x ; G_A^-)$, the distance between $x$ and $y$ in the BFS tree is the same as $d_{G_A^-}(x, y)$.
Consequently,
$
Z_k^-(x ; A)
\leq
\cW_k^{\sss (n)} (x ; A)
\leq
\cW_k^{\sss (n)} (x ; \Lambda_n)
\equald
\cW_k^{\sss (n)} 
$, 
$k\geq 1$, 
and in particular,
\begin{align}\label{eqn:85}
	\cD_A^-(x) \leq \sum_{k\geq 1} k\cW_k^{\sss (n)}\, ,
\end{align}
where the inequalities hold in the sense of stochastic domination.
Note that if $x\in A\subsetneq\Lambda_n$, then
$\pr\big( \cD_A^-(x) < \cD^-(x) \big) >0$ and $\pr\big( \cD^-(x) < \cD_A^-(x) \big)>0$
in the obvious coupling between $G_A^-$ and $G_{\sss\Lambda_n}^-$.
In other words, we cannot directly compare $\cD_A^-(x)$ and $\cD^-(x)$.
Instead, the relation \eqref{eqn:85} will be useful in the proof as we will see.

\medskip

\noindent{\bf Proof of Proposition~\ref{prop:HL-tau-asymptotics}:}
It follows from \eqref{eqn:uik-def}, \eqref{eqn:tau-dmax-def}, and \eqref{eqn:510} that 
\begin{align}\label{eqn:86}
	\tau 
	&
	= 
	\sum_{i=1}^{m_n} x_i^2 u_i
	= 
	L^{-4nd/3} \sum_{i=1}^{m_n} |\cC_i^-|^2 \cdot u_i
	=
	L^{-4nd/3} \sum_{i=1}^{m_n}\sum_{x\in\cC_i^-}\sum_{y\in\cC_i^-} d^-(x, y)
	\notag\\
	&
	=
	L^{-nd/3}\sum_{x\in\Lambda_n} |\Lambda_n|^{-1} \sum_{y\in\cC^-(x)} d^-(x, y)
	=
	L^{-nd/3}\cdot\bE\big( \cD^-(U)\, \big|\, G_{\sss\Lambda_n}^- \big)\, ,
\end{align}
where $U$ is uniformly distributed over $\Lambda_n$ independent of $G_{\sss\Lambda_n}^-$.
Note that
\begin{align}\label{eqn:87-a}
	\bE\big[ \bE\big( \cD^-(U)\, \big|\, G_{\sss\Lambda_n}^- \big) \big]
	=
	\bE\big[ \cD^-(U) \big]
	=
	\bE\big[ \cD^-(0) \big]
	=
	\bE\big[ \sum_{k\geq 1} k Z_k^-(0 ; \Lambda_n) \big]\, ,
\end{align}
where the second step uses symmetry.
Now,
\begin{align}\label{eqn:92-a}
	\bE\big( \sum_{k\geq 1} k\cW_k^{\sss (n)} \big)
	=
	\sum_{k\geq 1} k (m_n^{\sss (n)})^k
	=
	m_n^{\sss (n)} (1 - m_n^{\sss (n)})^{-2}
	=
	\zeta_n^2 L^{-2n(d - \theta)} (1+o(1))
	\asymp
	L^{2n(\theta - \alpha)}
	\ , 
\end{align}
where the penultimate step uses Lemma~\ref{lem:HL-zeta-n-m-j-n-asymptotics}(b) and the last step uses Lemma~\ref{lem:HL-zeta-n-m-j-n-asymptotics}(a).
Assume that $\cT_n^{\sss (n)}$ and the BFS tree of $\cC^-(0)$ are coupled as described before \eqref{eqn:coupling}.
In this coupling, for any sequence $(\omega_n ;\, n\geq 1)$ of positive numbers,
\begin{align*}
	&
	\bE\big[ \sum_{k\geq 1} k\cW_k^{\sss (n)} - \sum_{k\geq 1} k Z_k^-(0 ; \Lambda_n)  \big]
	\\
	&\hskip15pt
	\leq
	\bE\big[ \height(\cT_n^{\sss (n)}) \cdot \sum_{k\geq 1} \big( \cW_k^{\sss (n)} - Z_k^-(0 ; \Lambda_n) \big) \big]
	=
	\bE\big[ \height(\cT_n^{\sss (n)}) \cdot\big( |\cT_n^{\sss (n)}| - |\cC^-(0)| \big) \big]
	\\
	&\hskip30pt
	\leq
	\omega_n\cdot\bE\big( |\cT_n^{\sss (n)}| - |\cC^-(0)| \big) 
	+
	\bE\big( 
	\height( \cT_n^{\sss (n)} ) 
	\cdot 
	\ind_{\{ \height(\cT_n^{\sss (n)}) \geq\omega_n \}}
	\cdot
	|\cT_n^{\sss (n)}|  
	\big)
	\\
	&\hskip45pt
	\leq
	\omega_n\cdot\bE\big( |\cT_n^{\sss (n)}| - |\cC^-(0)| \big) 
	+
	\big[ 
	\bE\big( \height(\cT_n^{\sss (n)})^3 \big)
	\cdot 
	\pr\big( \height(\cT_n^{\sss (n)}) \geq\omega_n \big)
	\cdot
	\bE\big( |\cT_n^{\sss (n)}|^3 \big)
	\big]^{1/3}  .
\end{align*}
Using \eqref{eqn:9} and \eqref{eqn:111}, we see that there exists $C>0$ such that with the choice 
$\omega_n = CnL^{n(\theta - \alpha)}$, we have the bound
\begin{align}\label{eqn:93-a}
	\bE\big[ \sum_{k\geq 1} k\cW_k^{\sss (n)} - \sum_{k\geq 1} k Z_k^{-}(0 ; \Lambda_n)  \big]
	\leq
	CnL^{n(\theta - \alpha)} \bE\big( |\cT_n^{\sss (n)}| - |\cC^-(0)| \big) + o(1)
	=
	o(L^{2n(\theta-\alpha)}) \, ,
\end{align}
where the last step follows from Lemma~\ref{lem:3} and \eqref{eqn:HL-def-theta}.
Combining \eqref{eqn:87-a}, \eqref{eqn:92-a}, and \eqref{eqn:93-a} yields
\begin{align}\label{eqn:900}
	\bE\big[ \bE\big( \cD^-(U)\, \big|\, G_{\sss\Lambda_n}^- \big) \big]
	=
	\bE\big[ \cD^-(0) \big]
	=
	\zeta_n^2 L^{-2n(d - \theta)} (1+o(1))
	\asymp
	L^{2n(\theta - \alpha)}\, .
\end{align}
We now need a fluctuation bound on $\bE( \cD^-(U)\, |\, G_{\sss \Lambda_n}^- )$.
To that end, note that
\begin{align}\label{eqn:94}
	\bE\big[ \big( \bE\big( \cD^-(U)\, \big|\, G_{\sss \Lambda_n}^- \big) \big)^2 \big]
	&
	=
	\sum_{x\in\Lambda_n} |\Lambda_n|^{-1}
	\bE\big[ \cD^-(x) \cdot \bE( \cD^-(U)\, |\, G_{\sss \Lambda_n}^- ) \big]
	\notag\\
	&=
	\bE\big[ \cD^-(0) \cdot \bE( \cD^-(U)\, |\, G_{\sss \Lambda_n}^- ) \big]
	=
	\bE\big[ \cD^-(0) \cdot \cD^-(U) \big]\, ,
\end{align}
where the second step follows from the observation that 
$
\big( \cD^-(x) \, , \, \bE( \cD^-(U)\, |\, G_{\sss \Lambda_n}^- ) \big)
$
has the same distribution as 
$
\big( \cD^-(0) \, , \, \bE( \cD^-(U)\, |\, G_{\sss \Lambda_n}^- ) \big)
$
for each $x\in\Lambda_n$.
Let $\cH$ be a connected graph with vertex set $V\subseteq\Lambda_n$, where $0\in V$. 
Then 
\begin{align}\label{eqn:100}
	&
	\bE\big( \cD^-(U) \, |\, \cC^-(0) = \cH \big)
	=
	\bE\big( \cD^-(U)\cdot\ind_{\{U\in V\}} \, |\, \cC^-(0) = \cH \big)
	+
	\bE\big( \cD^-(U)\cdot\ind_{\{U\notin V\}} \, |\, \cC^-(0) = \cH \big)
	\notag\\
	&\hskip60pt
	=
	|\Lambda_n|^{-1} \sum_{x\in V}\sum_{y\in V} d_{\cH}(x, y)
	+
	|\Lambda_n|^{-1}
	\hskip-4pt
	\sum_{x\in\Lambda_n\setminus V}\bE\big( \cD_{\sss \Lambda_n\setminus V}^- (x) \big)
	=: 
	F_1 + F_2\, .
\end{align}
On the event $\{\cC^-(0) = \cH\}$,
\begin{align}\label{eqn:101}
	F_1 
	&
	\leq 
	|\Lambda_n|^{-1} \sum_{x\in V}\sum_{y\in V} \big( d_{\cH} (0, x) + d_{\cH} (0, y) \big)
	\notag\\
	&
	=
	2 |\Lambda_n|^{-1} |V| \cdot \big( \sum_{x\in V}  d_{\cH} (0, x) \big)
	=
	2 L^{-nd} |\cC^-(0)|\cdot\cD^-(0) \, .
\end{align}
Further, \eqref{eqn:85} implies that
\begin{align}\label{eqn:102}
	F_2
	\leq
	|\Lambda_n|^{-1}\cdot  |\Lambda_n\setminus V|\cdot\bE\big( \sum_{k\geq 1} k\cW_k^{\sss (n)} \big)
	\leq
	\bE\big( \sum_{k\geq 1} k\cW_k^{\sss (n)} \big)\, .
\end{align}
Combining \eqref{eqn:94}, \eqref{eqn:100}, \eqref{eqn:101}, and \eqref{eqn:102} shows that
\begin{align*}
	\bE\big[ \big( \bE\big( \cD^-(U)\, \big|\, G_{\sss \Lambda_n}^- \big) \big)^2 \big]
	&
	=
	\bE\big[ \cD^-(0) \cdot \cD^-(U) \big]
	=
	\bE\big[ \cD^-(0) \cdot \bE\big( \cD^-(U) \, | \, \cC^-(0) \big) \big]
	\\
	&
	\leq
	2 L^{-nd} \bE\big( |\cC^-(0)| \cdot \big( \cD^-(0) \big)^2 \big)
	+
	\bE\big( \cD^-(0) \big) \cdot \bE\big( \sum_{k\geq 1} k\cW_k^{\sss (n)}  \big)\, ,
\end{align*}
and consequently,
\begin{align}\label{eqn:103}
	&
	\var\big( \bE\big( \cD^-(U)\, \big|\, G_{\sss \Lambda_n}^- \big) \big)
	\notag\\
	&
	\hskip30pt
	\leq 
	2 L^{-nd} \bE\big( |\cC^-(0)| \cdot \big( \cD^-(0) \big)^2 \big)
	+
	\bE\big( \cD^-(0) \big)
	\cdot
	\big[ 
	\bE\big( \sum_{k\geq 1} k\cW_k^{\sss (n)} \big)
	-
	\bE\big( \cD^-(0) \big)
	\big]
	\notag\\
	&
	\hskip60pt
	\leq
	2 L^{-nd} \bE\big( |\cT_n^{\sss (n)}| \cdot \big( \sum_{k\geq 1} k \cW_k^{\sss (n)} \big)^2 \big)
	+
	o(L^{4n(\theta - \alpha)})\, , 
\end{align}
where the last step follows from \eqref{eqn:900}, \eqref{eqn:93-a}, and the coupling described right before \eqref{eqn:coupling}.
Now, for any sequence $(\omega_n\, ;\, n\geq 1)$ of positive numbers,
\begin{align*}
	\bE\big( |\cT_n^{\sss (n)}| \cdot \big( \sum_{k\geq 1} k \cW_k^{\sss (n)} \big)^2 \big)
	&
	\leq
	\bE\big( \big( \height(\cT_n^{\sss (n)}) \big)^2\cdot  |\cT_n^{\sss (n)}|^3 \big)
	\\
	&
	\leq
	\omega_n^2\,\, \bE\big( |\cT_n^{\sss (n)}|^3 \big)
	+
	\bE\big( 
	\big( \height(\cT_n^{\sss (n)}) \big)^2 \cdot  
	|\cT_n^{\sss (n)}|^3 \cdot 
	\ind_{\{ \height(\cT_n^{\sss (n)}) \geq\omega_n  \}} 
	\big)\, .
\end{align*}
We can use \eqref{eqn:9} and \eqref{eqn:111} to select $C>0$ such that with the choice 
$\omega_n = CnL^{n(\theta - \alpha)}$, the bound in \eqref{eqn:103} yields
\begin{align}\label{eqn:103-a}
	\var\big( \bE\big( \cD^-(U)\, \big|\, G_{\sss \Lambda_n}^- \big) \big)
	&
	\lesssim
	L^{-nd}\cdot n^2 L^{2n( \theta - \alpha )} \bE\big( |\cT_n^{\sss (n)}|^3 \big)
	+
	o(L^{4n(\theta - \alpha)})
	\notag\\
	&
	\lesssim
	n^2 L^{-nd} L^{2n( \theta - \alpha )} L^{5n(\theta - \alpha)} + o(L^{4n(\theta - \alpha)})
	=
	o(L^{4n(\theta - \alpha)})\, ,
\end{align}
where the penultimate step uses \eqref{eqn:111}, and the final step uses \eqref{eqn:HL-def-theta}.
The proof of Proposition~\ref{prop:HL-tau-asymptotics} is completed upon combining \eqref{eqn:86}, \eqref{eqn:900}, and \eqref{eqn:103-a}.
\qed

\subsection{Proof of Theorem~\ref{thm:HL-scaling-limit}}\label{sec:HL-proof-scaling-limit}
Throughout Section~\ref{sec:HL-proof-scaling-limit}, we work with $\alpha\in(0, 5d/6)$ and $\theta$ as in \eqref{eqn:HL-def-theta}.
Let $\cG(\vx, q)$ and $\bar\cG(\vx, q, \vM)$ be as defined right before \eqref{eqn:def-sigma_k-x_max} and \eqref{eqn:uik-def} respectively corresponding to $\vx, q,$ and $\vM$ given by \eqref{eqn:510}.
Using the same notation as introduced around \eqref{eqn:def-sigma_k-x_max} and \eqref{eqn:uik-def}, we let $\cC_i$ denote the component of $\cG(\vx, q)$ with the $i$-th maximal mass where 
$\mass (\cC_i)=\sum_{j \in V(\cC_i)} x_j$, 
and we write $\bar\cC_i$ for the corresponding component in $\bar\cG(\vx, q, \vM)$.
It follows from Propositions~\ref{prop:HL-xmax-diamax-bound}, \ref{prop:HL-sigma-2-asymptotics}, \ref{prop:HL-sigma-3-asymptotics}, \ref{prop:HL-tau-asymptotics}, and the choice of $\theta$ in \eqref{eqn:HL-def-theta} that Conditions~\ref{ass:aldous-basic-assumption} and \ref{ass:gen-2} are satisfied by $\vx, q,$ and $\vM$.
Hence, Theorems~\ref{thm:aldous-review} and \ref{thm:gen-2} show that
\begin{gather}
\big(\big(\mass( \cC_i),\spls(\cC_i) \big);\, i \geq 1 \big) 
\weakc 
\mvXi(\lambda)	
\ \text{ as }\  n \to \infty,\ \text{ and}	
\label{eqn:352-a}\\
\big(
\scl\big( L^{-nd/3}, 1 \big) \bar \cC_i\, ;\, i \geq 1 
\big) \weakc \vCrit(\lambda)
\ \text{ as }\  n \to \infty,
\label{eqn:352}
\end{gather}
where we have used the fact that 
$
(\sigma_2 + \tau)^{-1}\sigma_2^2 
=
\tau^{-1} \sigma_2^2 \cdot (1 + o_P(1))
=
L^{-nd/3} (1 + o_P(1))
$.

Now, $\HL_{\sss \Lambda_n}(\kappa_{\lambda}^{\sss (n)})$ can be obtained from $G_{\sss\Lambda_n}^-$ by placing edges independently with probability 
\begin{align}\label{eqn:HL-def-ft-n}
\ft_n
:=
1 - \exp\big( -\big( \zeta_n^{-1} L^{-n\theta} + \lambda L^{-4nd/3} \big) \big)
\end{align}
between each $y\neq z\in\Lambda_n$ such that $\{ y, z \}\notin E(G_{\sss\Lambda_n}^-)$.
Thus, $\HL_{\sss \Lambda_n}(\kappa_{\lambda}^{\sss (n)})$ can be coupled with $G_{\sss\Lambda_n}^-$, $\cG(\vx, q)$, and $\bar\cG(\vx, q, \vM)$ via the following prescription:
\begin{inparaenumi}
\item
Construct $G_{\sss\Lambda_n}^-$, $\cG(\vx, q)$, and $\bar\cG(\vx, q, \vM)$.
\item
Conditional on step (i), for each $1\leq \ell\leq m_n$ and $y\neq z\in V(\cC_{\ell}^-)$ with $\{y, z\}\notin E(\cC_{\ell}^-)$, place an edge between $y$ and $z$ with probability $\vt_n$, and do this independently across different pairs $\{y, z\}$.
Write $Y_{\ell\ell}$ for the number of new edges thus created within $\cC_{\ell}^-$, $1\leq\ell\leq m_n$.
\item
Conditional on steps (i) and (ii), generate independent random variables $Y_{jk}$ for $1\leq j<k \leq m_n$ with $\{j, k\}\in E( \cG(\vx, q) )$, where $Y_{jk}$ is a 
$\text{Binomial}\big( |\cC_j^-| \cdot |\cC_k^-| \, ,\ \ft_n \big)$
random variable conditioned to be at least one.
\item
Conditional on steps (i), (ii), and (iii), for each $1\leq j<k \leq m_n$ with $\{j, k\}\in E( \cG(\vx, q) )$, let $\cA_{jk}$ denote a random sample of size $Y_{jk} - 1$ drawn without replacement from $(V(\cC_j^-)\times V(\cC_k^-))\setminus e_{jk}$, where $e_{jk}$ denotes the edge between $\cC_j^-$ and $\cC_k^-$ present in $\bar\cG(\vx, q, \vM)$;
assume that $\cA_{jk}$ are independent across $\{j, k\}\in E( \cG(\vx, q) )$.
Create an edge for each element in $\bigcup_{\star} \cA_{jk}$, where `$\star$' denotes union over $\{j, k\}\in E( \cG(\vx, q) )$.
\end{inparaenumi}
The resulting graph has the same distribution as $\HL_{\sss \Lambda_n}(\kappa_{\lambda}^{\sss (n)})$.

In the above coupling,  
$V(\bar\cC_i) = V\big( \cC_i^{\sss (n)} (\kappa_{\lambda}^{\sss (n)}) \big)$ 
and 
$E(\bar\cC_i) \subseteq E\big( \cC_i^{\sss (n)} (\kappa_{\lambda}^{\sss (n)}) \big)$
for each $i\geq 1$. 
In view of \eqref{eqn:352}, to complete the proof of Theorem~\ref{thm:HL-scaling-limit}, it is enough to show that for each $i\geq 1$,
\begin{gather}
\pr\big(
Y_{jk}\geq 2
\text{ for some }
1\leq j< k\leq m_n\text{ with }
\{j, k\}\in E(\cC_i)
\big)
\to 0 \ \text{ and}
\label{eqn:353}
\\
\pr\big(
\exists  \ell \in V(\cC_i)  \text{ such that }
Y_{\ell\ell}\geq 1
\big)
\to 0
\label{eqn:354}
\end{gather}
as $n\to\infty$.
For notational simplicity, we will prove \eqref{eqn:353} and \eqref{eqn:354} for $i=1$; the proof for the general case is identical.

To that end, note that 
\begin{align}\label{eqn:355}
\vt_n
\lesssim 
(1+|\lambda|) \cdot \zeta_n^{-1} L^{-n\theta}
\lesssim
(1+|\lambda|) \cdot L^{-n(d - \alpha + \theta)}\, ,
\end{align}
where we have used \eqref{eqn:HL-def-theta} and Lemma~\ref{lem:HL-zeta-n-m-j-n-asymptotics}(a).
Consequently, 
\begin{align}\label{eqn:356}
\ft_n\cdot \max_{1\leq \ell_1, \ell_2\leq m_n} |\cC_{\ell_1}^-|\cdot |\cC_{\ell_2}^-|
\leq
n L^{-n(d - \alpha + \theta)} (C_{\ref{eqn:767}}(1))^2 n^2 L^{4n(\theta - \alpha)}
\weakc
0\ \text{ as }\ n\to\infty\, ,
\end{align}
where the first step uses \eqref{eqn:355} and \eqref{eqn:768} and holds whp, and the second step uses \eqref{eqn:HL-def-theta}.
Let $\cF$ denote the $\sigma$-field generated by 
$G_{\sss\Lambda_n}^-$, $\cG(\vx, q)$, and $\bar\cG(\vx, q, \vM)$.
Then on the event 
$
\cE_n 
:=
\big\{
\ft_n\cdot \max\limits_{\sss 1\leq \ell_1,\, \ell_2\leq m_n} |\cC_{\ell_1}^-|\cdot |\cC_{\ell_2}^-|
\leq
n^4 L^{-n(d - 3\theta + 3\alpha)}
\big\}
$,
\begin{align}\label{eqn:357}
\pr\big( Y_{jk}\geq 2\, \big|\, \cF \big)
\leq
C_{\ref{eqn:357}} \cdot |\cC_j^-|\cdot |\cC_k^-|\cdot\ft_n
\end{align}
for any $1\leq j< k\leq m_n$ with $\{j, k\}\in E(\cC_1)$.
Let $\fT$ be any rooted spanning tree of $\cC_1$; let $\rho\in V(\cC_1)$ denote the root of $\fT$.
Write $\pa(j)$ for the parent of $j\in V(\cC_1)\setminus\{\rho\}$ in $\fT$.
Then, on the event $\cE_n$,
\begin{align}\label{eqn:358}
&
\sum_{\{j, k\}\in E(\cC_1)}	
\pr\big( Y_{jk}\geq 2\, \big|\, \cF \big)
\lesssim
\sum_{\{j, k\}\in E(\cC_1)}	 |\cC_j^-|\cdot |\cC_k^-|\cdot\ft_n
\notag\\
&\hskip35pt
\leq
\ft_n\cdot\big(
\sum_{\{\pa(j), j\}\in E(\fT)}	 |\cC_{\pa(j)}^-|\cdot |\cC_j^-|
+
\spls(\cC_1) \cdot  |\cC_1^-|^2
\big)
\notag\\[1pt]
&\hskip70pt
\leq
\ft_n |\cC_1^-|\cdot\big(
\sum_{\{\pa(j), j\}\in E(\fT)}	 |\cC_j^-|
+
\spls(\cC_1) \cdot  |\cC_1^-|
\big)
\notag\\[1pt]
&\hskip105pt
\leq
\ft_n |\cC_1^-|\cdot\big(
L^{2nd/3}\mass(\cC_1)
+
\spls(\cC_1) \cdot  |\cC_1^-|
\big)
\notag\\[1pt]
&\hskip140pt
=
\ft_n |\cC_1^-|\cdot\big(
O_P( L^{2nd/3} )
+
O_P(1) \cdot  |\cC_1^-|
\big)
=
o_P(1)\, ,
\end{align}
where the first step uses \eqref{eqn:357}, 
the penultimate step uses \eqref{eqn:352-a},
and the last step uses \eqref{eqn:355}, \eqref{eqn:768}, and \eqref{eqn:HL-def-theta}.
Since $\pr(\cE_n)\to 1$ by \eqref{eqn:356}, the convergence in \eqref{eqn:353} follows, for $i=1$, from \eqref{eqn:358}.
Turning to the proof of \eqref{eqn:354}, we first note that conditional on $\cF$, $Y_{\ell\ell}$ is stochastically dominated by a $\text{Binomial}( |\cC_{\ell}^-|^2,\, \ft_n )$ random variable, $1\leq \ell\leq m_n$.
Hence,
\begin{align*}
&
\sum_{\ell\in V(\cC_1)} \pr( Y_{\ell\ell} \geq 1 \, |\, \cF )
\leq
\sum_{\ell\in V(\cC_1)}  |\cC_{\ell}^-|^2\cdot \ft_n 
\leq
\ft_n \cdot |\cC_1^-|\sum_{\ell\in V(\cC_1)}  |\cC_{\ell}^-|
\\
&\hskip25pt
=
\ft_n \cdot |\cC_1^-| \cdot L^{2nd/3}\mass(\cC_1)
\lesssim
( 1 + |\lambda| ) \cdot L^{-n(d + \theta - \alpha)} \cdot n L^{2n (\theta - \alpha)} \cdot L^{2nd/3}\cdot O_P(1)
=
o_P(1)\, ,
\end{align*}
where the penultimate step uses \eqref{eqn:355}, \eqref{eqn:768}, \eqref{eqn:352-a}, and holds whp, and the last step uses \eqref{eqn:HL-def-theta}.
This yields \eqref{eqn:354} for $i=1$.
As mentioned before, the proofs of \eqref{eqn:353} and \eqref{eqn:354} for general $i$ are identical.
This completes the proof of Theorem~\ref{thm:HL-scaling-limit}.
\qed

\subsection{Proof of Theorem~\ref{thm:HL-surplus-convergence}}\label{sec:HL-proof-surplus}
Throughout Section~\ref{sec:HL-proof-surplus}, we work with $\alpha\in(0, 2d/3)$ and $\theta$ as in \eqref{eqn:HL-def-theta}.
We first collect two results that will be needed in the proof.
\begin{lem}\label{lem:34}
We have,
\[
\pr\big(\spls( \cC^-(0) ) \neq 0 \big)
\lesssim
\begin{cases}
L^{ -n( d + 2\alpha - 2\theta ) } ,\ \ \text{ if } 0<\alpha\leq d/2\, ,\\
L^{-n (2d - \alpha - \theta)} ,\ \ \text{ if } d/2<\alpha< 2d/3.
\end{cases}
\]
\end{lem}
\noindent{\bf Proof:}
Note that if $\spls( \cC^-(0) ) \neq 0 $, then there exist $x \neq y\in\Lambda_n$ with $y\neq 0$ such that the event $\{0\connects x\} \disjt \{x\connects y\} \disjt \{x\edge y\}$ occurs.
Hence,
\begin{align}\label{eqn:359}
&
\pr\big(\spls( \cC^-(0) ) \neq 0 \big)
\leq
\sum_{x\in\Lambda_n}\sum_{y\in\Lambda_n\setminus\{x\}}
\pr\big( \{0\connects x\} \disjt \{x\connects y\} \disjt \{x\edge y\} \big)
\notag
\\
&
\hskip30pt
\leq
\sum_{x\in\Lambda_n}\sum_{y\in\Lambda_n\setminus\{x\}}
\pr(0\connects x) \pr(x\connects y) \pr(x\edge y)
\notag
\\
&
\hskip60pt
\leq 
\sum_{x\in\Lambda_n} \pr(0\connects x) \cdot \tilde\Delta_n
=
\bE( |\cC^-(0)| )\cdot\tilde\Delta_n
\leq
(1-m_n^{\sss (n)})^{-1}\cdot \tilde\Delta_n \ ,
\end{align}
where the second step uses the BK inequality, the third step uses the definition of $\tilde\Delta_n$ from \eqref{eqn:HL-def-Delta}, and the last step uses \eqref{eqn:coupling}.
The claimed bound follows upon combining Lemma~\ref{lem:HL-zeta-n-m-j-n-asymptotics}(b) and Lemma~\ref{lem:2} with \eqref{eqn:359}.
\qed 

\vskip5pt

For $n\geq 2$, we say that $(v(i) \, ;\, i\in [n])$ is a size-biased permutation of $1, 2, \ldots, n$ using the weights $\mvy = (y_i > 0\, ;\, i\in [n])$ if
$\pr( v(1) = j ) = \big(\sum\limits_{\sss k\in [n]} y_k\big)^{-1} y_j$, $1\leq j\leq n$, 
and for $1\leq i\leq n -1$,
\[
\pr\big( v(i + 1) = j\, \big|\, v(1),\ldots, v(i) \big) 
= 
\bigg( \sum_{k\in [n]} y_k\ - \sum_{k=1}^i y_{v(k)} \bigg)^{-1} y_j 
\, , \ \ j\in [n]\setminus\{ v(1), \ldots, v(i) \}\, .
\]
The next lemma concerns a concentration result about the partial sums of a size-biased arrangement of a sequence.

\begin{lem}[\cite{sbssxw-rank-1-scaling-limit}*{Lemma 8.2}]\label{lem:size-biased-partial-sum}
Fix $n\geq 2$, and let $\mvy$ and $(v(i) \, ;\, i\in [n])$ be as above.
Let $\mva  := (a_i \geq 0;\ i \in [n])$ be such that 
$c_n:= \sum_{i \in [n]} y_i a_i/{\sum_{i\in [n]} y_i} > 0$.  
Let $y_{\max}:= \max_{i \in [n]} y_i$ and $a_{\max} := \max_{i \in [n]} a_i$. 
Let $\ell= \ell(n)  \in [n]$ be such that
\begin{equation}\label{eqn:ass-size-biased-partial-sum}
\big( \sum_{k\in [n]} y_k \big)^{-1}\ell y_{\max} \to 0\, ,
\ \text{ and }\ \   
a_{\max}/( \ell c_n) \to 0 \ \ \text{ as }\ \ n\to\infty.
\end{equation}
Then
\begin{align}\label{eqn:107-a}
\sup_{k \leq \ell }
\Big| \frac{\sum_{i = 1}^k a_{v(i)}}{\ell c_n} - \frac{k}{\ell} \Big| 
\weakc 0\ \ \text{ as }\ \ n\to\infty.
\end{align}
\end{lem}

\noindent{\bf Proof of Theorem~\ref{thm:HL-surplus-convergence}:}
For simplicity, we will only prove the claimed convergences for the maximal component $\cC_1^{\sss (n)}(\kappa_{\lambda}^{\sss(n)})$.
Recall the construction of $\HL_{\sss \Lambda_n}(\kappa_{\lambda}^{\sss (n)})$ described right below \eqref{eqn:HL-def-ft-n}.
As shown in \eqref{eqn:353} and \eqref{eqn:354}, 
$\cC_1^{\sss (n)}(\kappa_{\lambda}^{\sss(n)}) = \bar\cC_1$ in this construction whp.
Suppose we show that
\begin{align}\label{eqn:446}
\pr\big(
\spls(\cC_i^-) = 0\ \text{ for all }\ \cC_i^-\subseteq \bar\cC_1
\big)
\to 1
\ \ \text{ as }\ \ n\to\infty\, .
\end{align}
This would, in particular, imply that $\spls(\bar\cC_1) = \spls(\cC_1)$ whp.
Thus, the desired result will follow from \eqref{eqn:352-a} and Theorem~\ref{thm:cycle-structure}.
The rest of this section is devoted to the proof of \eqref{eqn:446}.

Let us briefly recall a construction used by Aldous in \cite{aldous1997brownian}*{Section~3.1}.
Recall $\vx$ and $q$ from \eqref{eqn:510}.
Let $\xi_{i,j}$, $i\neq j \in [m_n]$, be independent exponential random variables such that 
$\bE[\xi_{i,j}] = (q x_j)^{-1}$. 
Select a vertex $v(1)$ with $\pr(v(1) =k) = x_k/\sigma_1$ for $k\in [m_n]$.
We declare the children of $v(1)$ to be the vertices 
$\{i\in[m_n]\, :\, \xi_{v(1),i} \leq x_{v(1)} \}$, 
and place an edge between $v(1)$ and each of its children.
Write $c(1)$ for the number of children of $v(1)$, labeled as $v(2), v(3), \ldots, v(c(1)+1)$ in increasing order of the $\xi_{v(1), v(i)}$ values. 
Then we move to $v(2)$, and obtain the children of $v(2)$ as the vertices $i\notin\{ v(1),\ldots, v(c(1)+1)\}$ such that $\xi_{v(2), i}\leq x_{v(2)}$. 
Label them as $v\big( c(1)+2 \big), \ldots, v\big( c(1)+c(2)+1 \big)$ in increasing order of their $\xi_{v(2),i}$ values, and place an edge between $v(2)$ and each of its children.
Then we move to $v(3)$, and so on.
Proceed recursively until it is no longer possible to add new vertices to the tree being grown.
Then select a new vertex among the unexplored vertices with probability proportional to their $x$-values, and continue in a similar fashion until all $m_n$ vertices have been explored.
This produces a random forest $\vF$ with the following property:
one can construct a random graph with the same distribution as $\cG(\vx, q)$ by suitably placing edges within the components of $\vF$.
In particular, the random partition of $[m_n]$ into different components of $\vF$ has the same distribution as the partition of $[m_n]$ into different components of  $\cG(\vx,q)$.

The following two properties of this construction will be relevant to us:
\begin{inparaenuma}
\item
The sequence $(v(i) \, ;\, i\in [m_n])$ is a size-biased permutation of $1, 2, \ldots, m_n$ using the weights $\vx$.
This is a simple consequence of the properties of exponential random variables.
\item
Consider the component of $\vF$ with the maximal mass, and suppose the last vertex in this component to be explored is $v( N^{{\sss (n)}} )$.
(In other words, when the exploration of the maximal component of $\vF$ ends, a total of $N^{\sss(n)}$ vertices have been explored.)
Then it follows from the proof of \cite{aldous1997brownian}*{Proposition~4} that
$\sum_{i=1}^{N^{\sss (n)}} x_{v(i)}$ converges in distribution, as $n\to\infty$, to the right endpoint of the longest excursion of the process $\bar W_{\lambda}(\cdot)$ defined in \eqref{eqn:parabolic-bm}.
In particular, for any $\eps>0$, we can choose $K_{\eps}>0$ such that
\begin{align}\label{eqn:543}
\pr\big(
\sum_{i=1}^{N^{\sss (n)}} x_{v(i)} 
\geq
K_{\eps}
\big)
\leq
\eps \ \ \text{ for all }\ \ n.
\end{align}
\end{inparaenuma}

Let $R_n := \sigma_1/\sigma_2$.
Since $\sigma_1 = L^{-2nd/3}\sum_{i=1}^{m_n} |\cC_i^-| = L^{nd/3}$, it follows from Lemma~\ref{lem:HL-zeta-n-m-j-n-asymptotics}(a)  and \eqref{eqn:41} that
\begin{align}\label{eqn:544}
R_n = \Theta_P\big( L^{n( \frac{2d}{3} + \alpha - \theta ) }  \big) \, .
\end{align}
Now,
\begin{align}\label{eqn:545}
&
\pr\big( N^{\sss (n)} \geq 2 K_{\eps} R_n \big)
\leq
\pr\bigg( \sum_{i=1}^{N^{\sss (n)}} x_{v(i)} \geq K_{\eps} \bigg)
+
\pr\bigg( N^{\sss (n)} \geq 2 K_{\eps} R_n \, ,  \ \sum_{i=1}^{N^{\sss (n)}} x_{v(i)} \leq K_{\eps}  \bigg)
\notag\\
&\hskip50pt
\leq
\eps + \pr\bigg( \sum_{i=1}^{2 K_{\eps} R_n} x_{v(i)} \leq K_{\eps}  \bigg)
\leq
\eps 
+ 
\pr\bigg( 
\sup_{k\leq 2K_{\eps} R_n} \bigg| \sum_{i=1}^k x_{v(i)} - \frac{k}{R_n} \bigg| \geq K_{\eps} 
\bigg)\, ,
\end{align}
where the second step uses \eqref{eqn:543}.
To bound the second term in the last step of \eqref{eqn:545}, we apply Lemma~\ref{lem:size-biased-partial-sum} with $m_n$ replacing $n$, $y_i = a_i = x_i$, $i\in [m_n]$ (consequently, $c_n = 1/R_n$), and $\ell = R_n$.
It follows from \eqref{eqn:768}, \eqref{eqn:41}, and Lemma~\ref{lem:HL-zeta-n-m-j-n-asymptotics}(a) that \eqref{eqn:ass-size-biased-partial-sum} is satisfied.
Hence, \eqref{eqn:107-a} combined with \eqref{eqn:545} shows that
\begin{align}\label{eqn:546}
N^{\sss (n)} \leq n R_n\  \text{ whp. }
\end{align}
For $i\in [m_n]$, let $f_i := \ind_{ \{\spls (\cC_i^-) \neq 0 \} }$.
Then 
\begin{align*}
\pr\bigg( f_{v(i)} = 1 \,\big|\, G_{\sss\Lambda_n}^- , v(1), \ldots, v(i-1) \bigg)
= 
\frac{\sum_{j=1}^{m_n} x_j f_j - \sum_{j=1}^{i-1} x_{v(j)} f_{v(j)} }{ \sigma_1 - \sum_{j=1}^{i-1} x_{v(j)} }
\, ,\ \ 2\leq i\leq m_n\, .
\end{align*}
Hence, on the event 
\[
\cE_n
:=
\big\{
\sigma_2 > 2n x_{\max}\, , \ R_n\leq n L^{n( \frac{2d}{3} + \alpha - \theta ) } 
\big\}\, ,
\]
for any $1\leq i\leq nR_n$,
\begin{align}\label{eqn:547}
\pr\big( f_{v(i)} = 1 \,\big|\, G_{\sss\Lambda_n}^- \big)
\leq
\frac{\sum_{j=1}^{m_n} x_j f_j }{ \sigma_1 - x_{\max}\cdot nR_n }
\leq
\frac{2}{\sigma_1} \sum_{j=1}^{m_n} x_j f_j 
=
2\cdot \pr\big( \spls(\cC^- (U)) \neq 0 \, \big|\, G_{\sss\Lambda_n}^- \big)
\, ,
\end{align}
where the second step uses the definition of the event $\cE_n$, and $U$ is uniformly distributed over $\Lambda_n$ independent of $G_{\sss \Lambda_n}^-$.
Consequently,
\begin{align}\label{eqn:548}
&
\pr\big(
\cE_n\cap
\big\{ f_{v(i)} = 1 \text{ for some } 1\leq i\leq nR_n\}
\big)
\leq
\bE\bigg[
\ind_{\cE_n}\cdot 
\sum_{i=1}^{nR_n} \pr\big( f_{v(i)} = 1 \,\big|\, G_{\sss\Lambda_n}^- \big)
\bigg]
\notag\\
&\hskip60pt
\leq
\bE\bigg[   
2n^2 L^{n( \frac{2d}{3} + \alpha - \theta ) }
\cdot
 \pr\big( \spls(\cC^- (U)) \neq 0 \, \big|\, G_{\sss\Lambda_n}^- \big)
\bigg]
\notag\\
&\hskip120pt
=
2n^2 L^{n( \frac{2d}{3} + \alpha - \theta ) }
\cdot
\pr\big( \spls(\cC^- (0)) \neq 0\big)
=
o(1)\, ,
\end{align}
where the second step uses the definition of $\cE_n$ and \eqref{eqn:547}, and the last step uses \eqref{eqn:HL-def-theta} and Lemma~\ref{lem:34}.
It follows from \eqref{eqn:768}, \eqref{eqn:41}, Lemma~\ref{lem:HL-zeta-n-m-j-n-asymptotics}(a), and \eqref{eqn:544} that $\pr(\cE_n)\to 1$ as $n\to\infty$, which combined with \eqref{eqn:548} and \eqref{eqn:546} shows that
$
\pr\big(
f_{v(i)} = 1 \text{ for some } 1\leq i\leq N^{\sss (n)}
\big)
= o(1)
$.
This fact, in particular, implies \eqref{eqn:446}, thereby completing the proof.
\qed

\subsection{Proof of Proposition~\ref{prop:HL-phase-transition}}\label{sec:HL-proof-phase-transition}
Throughout Section~\ref{sec:HL-proof-phase-transition}, we work with $\alpha\in(0, d)$.
The parameter $\theta$ introduced in \eqref{eqn:HL-theta-full-range} plays no role in this section, so the threshold $n_0$ such that a relation ``$\lesssim$" holds for all $n\geq n_0$ will depend only on $\pmtr^{\ast}$ (instead of $\pmtr$) in this section.
Define
\begin{align}\label{eqn:676}
\tilde\zeta_n 
:= 
\sum_{x\in\Lambda_n} J(0, x)
=
(L^d - 1)\sum_{i=1}^{n} L^{(i-1)d} \rho(L^i)\, ,
\end{align}
where $J$ and $\rho$ are as introduced around \eqref{eqn:1a}.
For $\eps\in (0, 1)$, let $\tilde G_{\sss \Lambda_n}^{\sss (-\eps)}$ denote the random graph on vertex set $\Lambda_n$, where we place edges independently between each $x\neq y\in\Lambda_n$ with probability 
\[
\tilde r_{-\eps}^{\sss (n)}(x, y)
:=
(1-\eps) \cdot \tilde\zeta_n^{-1} J(x, y).
\]
Let $\tilde\cC_{i; n}^{\sss (-\eps)}$ denote the $i$-th largest component of $\tilde G_{\sss \Lambda_n}^{\sss (-\eps)}$, $i\geq 1$. We define
\begin{align}\label{eqn:677}
\tilde s_{2; n}^{\sss (-\eps)}
:=
L^{-nd}\sum_{i\geq 1} | \tilde\cC_{i; n}^{\sss (-\eps)} |^2\, .
\end{align}
 Let $\tilde\fT^{\sss (-\eps)}_{n}$ denote a branching process tree with offspring distribution equal to the law of $\sum_{i=1}^n \tilde X_{i; n}^{\sss (-\eps)}$, where 
$\tilde X_{i; n}^{\sss (-\eps)}$, $1\leq i\leq n$, are independent and 
\[
\tilde X_{i; n}^{\sss (-\eps)}
\sim
\text{Binomial}\big(
(L^d - 1)L^{(i-1)d} \, ,\, (1-\eps)\cdot \tilde\zeta_n^{-1} \rho(L^i)
\big).
\]
By an argument similar to the one leading to \eqref{eqn:coupling}, for any $x\in\Lambda_n$, there is a coupling in which 
\begin{align}\label{eqn:678}
|\cC( x; \tilde G_{\sss \Lambda_n}^{\sss (-\eps)} )|
\leq
|\tilde\fT^{\sss (-\eps)}_{n}|\, .
\end{align}


\medskip

\noindent{\bf Proof of Proposition~\ref{prop:HL-phase-transition}(a):}
Note that for any $\eps\in(0, 1)$ and $t\in [0, 1]$, 
\begin{align}\label{eqn:678-a}
&
\bE\big[ \exp\big( t \sum_{i=1}^n \tilde X_{i; n}^{\sss (-\eps)} \big) \big]
=
\prod_{i=1}^n \bigg( 1 + \tilde\zeta_n^{-1} (1-\eps) \cdot \rho(L^i) \cdot (e^t - 1) \bigg)^{ (L^d - 1) L^{(i-1)d } }
\notag\\
&\hskip40pt
\leq
\exp\bigg( \tilde\zeta_n^{-1} (1-\eps)\cdot (e^t - 1) \sum_{i=1}^n (L^d - 1) L^{ (i-1)d } \rho(L^i) \bigg)
\notag\\
&\hskip80pt
=
\exp\big(  (1-\eps)\cdot (e^t - 1) \big)
\leq
\exp\big( (1-\eps) t + C_{\ref{eqn:678-a}} t^2 \big)
\end{align}
for some universal constant $C_{\ref{eqn:678-a}} > 1/2$, where the third step uses \eqref{eqn:676}.
Hence, for any $\eps\in(0, 1)$ and $0\leq t\leq \eps/(2C_{\ref{eqn:678-a}})$,
$
\bE\big[ \exp\big( t \sum_{i=1}^n \tilde X_{i; n}^{\sss (-\eps)} \big) \big]
\leq
\exp\big( (1 - \eps/2 ) t\big)
$.
Now a chain of inequalities similar to the one in \eqref{eqn:14} will show that 
\begin{align}\label{eqn:678-b}
\pr\big( |\tilde\fT_n^{\sss (-\eps)}| >k \big) \leq \exp\big( - \eps^2 k /(4 C_{\ref{eqn:678-a}})\big)
\end{align}
for any $\eps\in(0, 1)$ and $k\geq 1$.
Consequently, for any $\eps\in(0, 1)$ and $k\geq 1$,
\[
\pr\big(  | \tilde\cC_{1; n}^{\sss (-\eps)} |  > k\big)
\leq
L^{nd} \pr\big(  \big| \cC\big( 0 ;  \tilde G_{\sss \Lambda_n}^{\sss (-\eps)} \big) \big|  > k\big)
\leq
L^{nd} \pr\big(   |\tilde\fT_n^{\sss (-\eps)}| >k \big)
\leq 
L^{nd} \exp\big( - \eps^2 k /(4 C_{\ref{eqn:678-a}})\big)
\, ,
\]
where the second step uses \eqref{eqn:678}.
Hence, for any $\eps>0$,
\begin{align}\label{eqn:678-c}
\pr\big(  | \tilde\cC_{1; n}^{\sss (-\eps)} |  \geq Cn\eps^{-2}\big) \to 0 \text{ as } n\to\infty
\end{align}
for a constant $C>0$ that depends only on $d$ and $L$.

Now, it follows from \eqref{eqn:30} and \eqref{eqn:676} that
\[
(L^d - 1)\cdot\sum_{i=1}^{n} 
L^{(i-1)d} \big[ 1 - \exp\big( - \zeta_n^{-1}\rho(L^i) \big) \big]
=
(L^d - 1)\cdot\sum_{i=1}^{n} 
L^{(i-1)d} \cdot \tilde\zeta_n^{-1}\rho(L^i) \, .
\]
Consequently,
\begin{align}\label{eqn:679}
0
\leq
\big( \zeta_n^{-1} - \tilde\zeta_n^{-1} \big) \cdot
\sum_{i=1}^{n} L^{(i-1)d} \rho(L^i)
\lesssim
\sum_{i=1}^{n} 
L^{(i-1)d} \cdot \big( \zeta_n^{-1} \rho(L^i) \big)^2\, .
\end{align}
Combining \eqref{eqn:1a} and Lemma~\ref{lem:HL-zeta-n-m-j-n-asymptotics}(a) with \eqref{eqn:679}, we see that
\begin{align}\label{eqn:680}
0
\leq
\zeta_n^{-1} - \tilde\zeta_n^{-1}
\lesssim
\begin{cases}
L^{-n(2d - \alpha)} \, , \ \text{ if } 0 < \alpha < d/2,\\
n L^{-3n(d - \alpha)} \, , \ \text{ if } \alpha = d/2,\\
L^{-3n(d - \alpha)} \, , \ \text{ if } d/2 < \alpha < d\, 
\end{cases}
\ \ \
=
o(\zeta_n^{-1})\, .
\end{align}
Consequently, there exists $n_0(\eps)\geq 1$ depending only on $\eps$ and $\pmtr^\ast$ such that for $n\geq n_0(\eps)$, we have
$
1 - \exp\big (- r_{-\eps}^{\sss (n)} (x, y) \big)
\leq
\tilde r_{-\eps/2}^{\sss (n)} (x, y) 
$
for all $x\neq y\in\bH_L^d$, where $r_{-\eps}^{\sss (n)}(\cdot\, , \cdot)$ is as in the statement of Proposition~\ref{prop:HL-phase-transition}.
Hence, for $n\geq n_0(\eps)$, there is a coupling between $\HL_{\sss\Lambda_n}(r_{-\eps}^{\sss (n)})$ and $\tilde G_{\sss \Lambda_n}^{\sss (-\eps/2)}$ in which the former is a subgraph of the latter.
This fact combined with \eqref{eqn:678-c} yields \eqref{eqn:4}.
\qed

\medskip

We now turn to the proof of Proposition~\ref{prop:HL-phase-transition}(b).
We will make use of the following result:

\begin{lem}[\cite{spencer2007birth}*{Theorem~3.1(2)}]\label{lem:spencer-wormald}
Suppose $c_0, c_1, c_2$, and $t$ are positive numbers with $t c_0>1$.
Let $H_n$, $n\geq 1$, be graphs (we view $H_n$ as the seed graph) with $|H_n| = n$ satisfying the following for all $n$:
\begin{enumeratei}
\item 
$n^{-1}\sum |\cC|^2 \geq c_0$, where the sum is over all connected components $\cC$ of $H_n$, and
\item 
$n^{-1}\#\big\{
v\in V(H_n)\, :\, | \cC( v; H_n ) | \geq s 
\big\}
\leq 
c_1 \exp(-c_2 s)$
for all $s\in\bZ_{>0}$.
\end{enumeratei}
Place edges with probability $t/n$ independently between all $u\neq v\in V(H_n)$ such that $\{u, v\}\notin E(H_n)$;
let $E_n^+$ denote the set of edges thus created.
Let $\cC_{\max}^+$ denote the largest component of the graph $(V(H_n), E(H_n)\cup E_n^+)$.
Then there exists $a_0>0$ depending only on $c_0, c_1, c_2$, and $t$ such that 
$ | \cC_{\max}^+ | \geq a_0 n $ whp.
\end{lem}

\noindent{\bf Proof of Proposition~\ref{prop:HL-phase-transition}(b)}
Fix $\delta\in (0, 1)$.
Let $U$ be uniformly distributed over $\Lambda_n$ independent of $\tilde G_{\sss \Lambda_n}^{\sss (-\delta)}$.
Consider the coupling between $\cC(U; \tilde G_{\sss \Lambda_n}^{\sss (-\delta)})$ and $\tilde\fT_n^{\sss (-\delta)}$ as mentioned around \eqref{eqn:678}.
From the construction of the BFS tree described right before \eqref{eqn:coupling}, it follows that in this coupling,
\begin{align}\label{eqn:681}
\pr\big(
|\cC(U; \tilde G_{\sss \Lambda_n}^{\sss (-\delta)})| \neq |\tilde\fT_n^{\sss (-\delta)}|
\big)
\leq
\bE\big( |\tilde\fT_n^{\sss (-\delta)}|^2 \big)
\max_{x, y\in\Lambda_n}\tilde r_{-\delta}^{\sss (n)} (x, y)
\leq
A_2\tilde\zeta_n^{-1} \bE\big( |\tilde\fT_n^{\sss (-\delta)}|^2 \big)\, ,
\end{align}
where the last step uses \eqref{eqn:1a}.
Similar to \eqref{eqn:HL-def-susceptibility}, we have 
$
\bE\big( \tilde s_{2; n}^{\sss (-\delta)} \big)
=
\bE\big( |\cC(U; \tilde G_{\sss \Lambda_n}^{\sss (-\delta)})| \big)
$.
Hence,
\begin{align}\label{eqn:682}
&
\bE\big( | \tilde \fT_n^{\sss (-\delta)} | \big) 
-
\bE\big( \tilde s_{2; n}^{\sss (-\delta)} \big)
=
 \bE\bigg( | \tilde \fT_n^{\sss (-\delta)} | 
 -
 |\cC(U; \tilde G_{\sss \Lambda_n}^{\sss (-\delta)})| 
 \bigg) 
 \leq
 \bE\bigg(
  | \tilde \fT_n^{\sss (-\delta)} | \cdot
  \ind_{\{  |\cC(U; \tilde G_{\sss \Lambda_n}^{\sss (-\delta)})| \neq |\tilde\fT_n^{\sss (-\delta)}| \}}
\bigg)
\notag\\
&
\hskip20pt
\leq
\bigg[
\bE\big( | \tilde \fT_n^{\sss (-\delta)} |^2 \big)
\cdot
\pr\big(  |\cC(U; \tilde G_{\sss \Lambda_n}^{\sss (-\delta)})| \neq |\tilde\fT_n^{\sss (-\delta)}| \big)
\bigg]^{1/2}
\lesssim
\tilde\zeta_n^{-1/2} \cdot \bE\big( | \tilde \fT_n^{\sss (-\delta)} |^2 \big)
\lesssim
\delta^{-3} \cdot \tilde\zeta_n^{-1/2} \, ,
\end{align}
where the first step holds in the coupling of \eqref{eqn:678}, the fourth step uses \eqref{eqn:681}, and the last step follows from arguments similar to the ones leading to \eqref{eqn:111}.
Now, from the definition of $\tilde\fT_n^{\sss (-\delta)}$, it follows that
$
\bE\big( |\tilde\fT_n^{\sss (-\delta)}| \big) = \delta^{-1}
$.
Since $\tilde\zeta_n\to\infty$ as $n\to\infty$, we see from \eqref{eqn:682} that
\begin{align}\label{eqn:683}
\bE\big( \tilde s_{2; n}^{\sss (-\delta)} \big)
=
\delta^{-1} + o(1) \ \text{ as } \ n\to\infty\, .
\end{align}
By arguments similar to the ones leading to \eqref{eqn:99} and \eqref{eqn:111}, we have
$
\var( \tilde s_{2;n}^{\sss (-\delta)} )
$
$
\leq 
L^{-nd} \bE\big( | \tilde\fT_n^{\sss (-\delta) } |^3 \big)
$
$
\lesssim
\delta^{-5} L^{-nd}
$, 
which combined with \eqref{eqn:683} yields
\begin{align}\label{eqn:684}
\delta\cdot \tilde s_{2; n}^{\sss (-\delta)} \weakc 1\ \ \text{ as }\ \ n\to\infty\, .
\end{align}
Using the inequalities $u\geq 1- e^{-u} \geq u(1-u)$, $u\geq 0$, we see that
\begin{gather*}
r_{-\delta}^{\sss (n)}(x, y)
\leq 
(1-\delta)\zeta_n^{-1} J(x, y)
=
\tilde r_{-\delta}^{\sss (n)} (x, y) \tilde\zeta_n\zeta_n^{-1}
\, , \ \ \text{ and }
\\[2pt]
r_{-\delta}^{\sss (n)}(x, y)
\geq 
(1-\delta)\zeta_n^{-1} J(x, y) \big( 1 - \zeta_n^{-1} J(x, y) \big)
\geq 
\tilde r_{-\delta}^{\sss (n)} (x, y) \tilde\zeta_n\zeta_n^{-1}\big( 1 - A_2\zeta_n^{-1}\big)\, ,
\end{gather*}
where $r_{-\delta}^{\sss (n)}(\cdot\, , \cdot)$ is as in the statement of Proposition~\ref{prop:HL-phase-transition}.
Now, $\tilde\zeta_n\zeta_n^{-1} \to 1$ by \eqref{eqn:680}, and $\zeta_n\to\infty$.
Hence, for any $0<\delta_1<\delta<\delta_2$, there is a natural coupling in which 
$
\tilde s_{2; n}^{\sss (-\delta_2)}
\leq
L^{-nd} \sum_i |\cC_i^{\sss (n)}( r_{-\delta}^{\sss (n)} ) |^2
\leq
\tilde s_{2; n}^{\sss (-\delta_1)}
$
for all $n$ greater than a threshold that depends only on $\delta, \delta_1, \delta_2$, and $\pmtr^\ast$.
Thus, \eqref{eqn:684} implies that
\begin{align}\label{eqn:685}
\delta \cdot L^{-nd} \sum_i |\cC_i^{\sss (n)}( r_{-\delta}^{\sss (n)} ) |^2
\weakc 
1 \ \ \text{ as }\ \ n\to\infty\, .
\end{align}

Next, for $y\in\Lambda_n$ and $s\geq 1$, define the event 
$
\cE_n(y; s)
:=
\big\{\hskip0.2pt
\big| \cC\big( y ; \HL_{\sss \Lambda_n}( r_{-\delta}^{\sss (n)} ) \big) \big| \geq s
\big\}
$.
Then for $y\in \Lambda_n\setminus\{0\}$,
\begin{align}\label{eqn:724}
&
\pr\big( \cE_n( 0 ; s ) \cap \cE_n( y ; s ) \big)
=
\pr\big( \cE_n( 0 ; s ) \cap \big\{ y\in \cC\big( 0 ; \HL_{\sss \Lambda_n}( r_{-\delta}^{\sss (n)} ) \big) \big\} \big)
\notag\\
&\hskip120pt
+
\pr\big( \cE_n( 0 ; s ) \cap \cE_n( y ; s ) 
\cap \big\{ y\notin \cC\big( 0 ; \HL_{\sss \Lambda_n}( r_{-\delta}^{\sss (n)} ) \big) \big\} \big)
\notag\\
&\hskip130pt
\leq
\pr\big( y\in \cC\big( 0 ; \HL_{\sss \Lambda_n}( r_{-\delta}^{\sss (n)} ) \big) \big)
+
\pr\big( \cE_n( 0 ; s ) \big) \cdot \pr\big( \cE_n( y ; s ) \big)\, .
\end{align}
Hence,
\begin{align}\label{eqn:725}
&
\var\big( \sum_{y\in\Lambda_n} \ind_{\cE_n(y ; s)} \big)
=
\sum_{y\in\Lambda_n} \var\big( \ind_{\cE_n(y ; s)} \big)
+
\sum_{y\neq z\in\Lambda_n} \cov\big( \ind_{\cE_n(y ; s)} , \ind_{\cE_n(z ; s)} \big)
\notag\\
&\hskip25pt
\leq
L^{nd} 
+ 
L^{nd} \sum_{y\in\Lambda_n\setminus\{0\}} \cov\big( \ind_{\cE_n(0 ; s)} , \ind_{\cE_n (y ; s)} \big)
\notag\\
&\hskip50pt
\leq
L^{nd}\big[
1 
+ 
\sum_{y\in\Lambda_n\setminus\{0\}} 
\pr\big( y\in \cC\big( 0 ; \HL_{\sss \Lambda_n}( r_{-\delta}^{\sss (n)} ) \big) \big)
\big]
=
L^{nd} \cdot \bE\big( |  \cC\big( 0 ; \HL_{\sss \Lambda_n}( r_{-\delta}^{\sss (n)} ) \big) | \big)\, ,
\end{align}
where the third step uses \eqref{eqn:724}.
Let $n_0(\cdot)$ be as introduced right below \eqref{eqn:680}.
Then for $n\geq n_0(\delta)$, 
\begin{align}\label{eqn:726}
\var\big( \sum_{y\in\Lambda_n} \ind_{\cE_n(y ; s)} \big)
\leq
L^{nd} \cdot \bE\big( | \tilde\fT_n^{\sss ( -\delta/2 ) } | \big)
=
2 \delta^{-1} L^{nd} 
\, ,
\end{align}
where the first step uses \eqref{eqn:725}, the definition of $n_0(\cdot)$, and \eqref{eqn:678}.
Further, for $n\geq n_0(\delta)$,
\begin{align}\label{eqn:129}
\bE\big( L^{-nd} \sum_{y\in\Lambda_n} \ind_{\cE_n(y ; s)} \big)
=
\pr\big( \cE_n(0 ; s) \big)
\leq
\pr\big( | \tilde\fT_n^{\sss ( -\delta/2 ) } | \geq s \big)
\leq
\exp\big( - C_{\ref{eqn:129}} \delta^2 s \big)\, ,
\end{align}
where the second step uses \eqref{eqn:678}, and the last step uses \eqref{eqn:678-b}.
Let
\begin{align}\label{eqn:0}
C_{ \ref{eqn:0} }
:=
\min\big\{ C_{ \ref{eqn:129} } \, ,\  (4 C_{ \ref{eqn:4} })^{-1} d\cdot\ln L  \big\} \, .
\end{align}
Then
\begin{align}\label{eqn:727}
&
\pr\big(
\exists s\geq 1 \text{ such that }
L^{-nd}\#\big\{
y\in\Lambda_n : \cE_n(y ; s) \text{ occurs}
\big\}
\geq 
2\exp\big( -C_{ \ref{eqn:0} }\delta^2 s \big)
\big)
\notag\\
&\hskip70pt
\leq
\sum_{s=1}^{C_{\ref{eqn:4}} n \delta^{-2}}
\pr\big(
L^{-nd} \sum_{y\in\Lambda_n} \ind_{ \cE_n(y ; s) } \geq 2\exp\big( -C_{ \ref{eqn:0} }\delta^2 s \big)
\big)
\notag\\
&\hskip140pt
+
\pr\big(
| \cC_1^{\sss (n)}(r_{-\delta}^{\sss (n)}) | > C_{\ref{eqn:4}} n \delta^{-2}
\big)
=:
F_1 + F_2\, .
\end{align}
For $n\geq n_0(\delta)$,
\begin{align}\label{eqn:728}
&
F_1
\leq
\sum_{s=1}^{C_{\ref{eqn:4}} n \delta^{-2}}
\pr\big(
L^{-nd} \sum_{y\in\Lambda_n} \big[ \ind_{ \cE_n(y ; s) } 
-
\pr( \cE_n(0 ; s) ) \big]
\geq 
\exp\big( -C_{ \ref{eqn:0} }\delta^2 s \big)
\big)
\notag
\\
&
\hskip50pt
\leq
\sum_{s=1}^{C_{\ref{eqn:4}} n \delta^{-2}}
\exp\big( 2C_{ \ref{eqn:0} }\delta^2 s \big)
\cdot
\frac{1}{L^{2nd}}
\cdot
\frac{2 L^{nd}}{\delta}
\leq
\frac{2 \exp\big( 2 C_{\ref{eqn:0}} C_{\ref{eqn:4}} n \big)}{\delta L^{nd} \cdot \big(\exp( 2 C_{\ref{eqn:0} } \delta^2)  - 1 \big) }
\notag
\\
&
\hskip100pt
\leq
2\big[ \delta L^{nd/2} \cdot \big(\exp( 2 C_{\ref{eqn:0} } \delta^2)  - 1 \big) \big]^{-1} \, ,
\end{align}
where the first step uses \eqref{eqn:129} and \eqref{eqn:0}, the second step uses \eqref{eqn:726}, and the last step uses \eqref{eqn:0}.
Combining \eqref{eqn:728} with \eqref{eqn:4} and \eqref{eqn:727} shows that the following events occur whp:
\begin{align}\label{eqn:729}
L^{-nd}\#\big\{
y\in\Lambda_n : 
\big| \cC\big( y ; \HL_{\sss \Lambda_n}( r_{-\delta}^{\sss (n)} ) \big) \big| \geq s
\big\}
\leq
2 \exp\big( -C_{\ref{eqn:0}} \delta^2 s\big) 
\ \text{ for all }s\geq 1\, .
\end{align}

Finally, to complete the proof of Proposition~\ref{prop:HL-phase-transition}(b), we apply Lemma~\ref{lem:spencer-wormald} with $\HL_{\sss \Lambda_n}( r_{-\delta}^{\sss (n)} )$ as the seed graph.
Note that for any $\eps>0$, we can construct $\HL_{\sss \Lambda_n}( r_{\eps}^{\sss (n)} )$ from $\HL_{\sss \Lambda_n}( r_{-\delta}^{\sss (n)} )$ by placing edges independently between each $x\neq y\in \Lambda_n$ that are not connected directly by an edge in $\HL_{\sss \Lambda_n}( r_{-\delta}^{\sss (n)} )$ with probability
\begin{align}\label{eqn:730}
&
1- \exp\big( -(\delta + \eps) \cdot \zeta_n^{-1} J(x, y) \big)
\geq 
(\delta + \eps) \cdot \zeta_n^{-1} J(x, y) \cdot\big( 1 - (1+\eps)\cdot\zeta_n^{-1} J(x, y) \big)
\notag
\\
&
\hskip20pt
\geq
C\cdot (\delta+\eps)\cdot L^{-n(d-\alpha)} A_1 L^{-n\alpha} 
\cdot\big( 1- C'\cdot(1+\eps)\cdot L^{-n(d-\alpha)} \big)
\geq
C_{\ref{eqn:730}}  (\delta+\eps)\cdot L^{-nd}  
\end{align}
for some $C_{\ref{eqn:730}} \in (0, 1)$, where we have used \eqref{eqn:1a} and Lemma~\ref{lem:HL-zeta-n-m-j-n-asymptotics}(a); 
this entire string of inequalities holds uniformly over $x\neq y \in\Lambda_n$ for $n$ greater than a threshold that depends only on $\eps$ and $\pmtr^\ast$.
Taking $0<\delta< \min\{1,\, C_{\ref{eqn:730}} \eps/2\}$ and using \eqref{eqn:685} and \eqref{eqn:729} in conjunction with Lemma~\ref{lem:spencer-wormald} yields the claim.
\qed

\section{Proofs: discrete torus}\label{sec:torus-proofs}
We will follow throughout Section~\ref{sec:torus-proofs} the convention about constants described right below \eqref{eqn:HL-def-parameters} with $\pmtr$ and $\pmtr^{\ast}$ as in \eqref{eqn:95}.
Further, the threshold $n_0 = n_0(\pmtr)$ implicit in a relation ``$\lesssim$" or ``$\asymp$" will be assumed to be at least $2^{10} + 1$.
As we had mentioned at the end of Section~\ref{sec:results-torus}, there are three main components in the proof.
The first one will be worked out in detail in Section~\ref{sec:torus-proof-two-point}.
The other two, being similar to the hierarchical model, will only be outlined.

\subsection{Proof of Theorem~\ref{thm:torus-two-point}}\label{sec:torus-proof-two-point}
For simplicity of notation, write $G_n^- = \torus_n(\kappa^{\sss (n), -})$.
For $A\subseteq\bT_n^d$, we will write $G_{\sss A}^-$ for the induced subgraph of $G_n^-$ on the vertex set $A$.
For $x\in\bT_n^d$, we let $\cC^-(x) = \cC(x ;\, G_n^-)$--the component of $x$ in $G_n^-$.
For $x, y\in\bT_n^d$, we will write $\{x\edge y\}$ for the event that there is an edge between $x$ and $y$ in $G_n^-$.
As defined right below \eqref{eqn:torus-def-fp}, $\{x\connects y\}$ denotes the event that there is a path between $x$ and $y$ in $G_n^-$, i.e., $y\in\cC^-(x)$.
For $A\subseteq\bT_n^d$ and $x, y\in A$, we let 
$\{ x\connects y \text{ in } A\}$ be the event that there is a path in $G_A^-$ connecting $x$ and $y$.

Similar to hierarchical percolation, the proof of Theorem~\ref{thm:torus-two-point} can be separated into the proof of an upper bound and the proof of a lower bound.
Recall $d_n$ from \eqref{eqn:torus-diameter}.

\begin{prop}\label{prop:torus-two-point-bound}[Upper bound]
	\begin{enumeratea}
		\item
		Let $\alpha\in (0, d)$ and $\alpha < \theta < \min\{ 2\alpha, \alpha + d/3 \}$.
		Then for $x\in\bT_n^d\setminus\{0\}$,
		\begin{align}\label{eqn:788-torus}
			\pr(0 \connects x)
			\lesssim
			\begin{cases}
				\zeta_n^{-1} \| x\|^{-\alpha}\, ,\ \ \text{ if }\ \ 1\leq \|x\| \leq n^{a_0}\, ,
				\\[3pt]
				n^{-( d + \alpha - \theta ) }\, , \text{ if }\ \ n^{a_0}< \|x\| \leq d_n\, ,
			\end{cases}
		\end{align}
		where $a_0 = 2 - \theta\alpha^{-1}$.

		\item 
		Let $\alpha\in (0, d/3)$ and $2\alpha \leq \theta < \alpha + d/3 $.
		Then
		\begin{align}\label{eqn:789-torus}
			\pr(0 \connects x)
			\lesssim
			n^{-( d + \alpha - \theta ) }
			\ \ \text{ for }\ \ x\in\bT_n^d\setminus\{0\}.
		\end{align}
		
	\end{enumeratea}

\end{prop}

\begin{prop}\label{prop:torus-two-point-lower-bound}[Lower bound]
	Suppose either $\alpha\in (0, 5d/6]$ and $\theta\in (\alpha, \alpha + d/3)$, or 
	$\alpha\in (5d/6, d)$ and $\theta\in ( \alpha, 2d - \alpha )$.
	Then
	\begin{align}\label{eqn:790-torus}
		\pr(0 \connects x)
		\gtrsim
		n^{-( d + \alpha - \theta ) }
		\ \ \text{ for }\ \ x\in\bT_n^d\setminus\{0\}.
	\end{align}
\end{prop}

Similar to the explanation given around \eqref{eqn:791}, combining \eqref{eqn:790-torus} with Proposition~\ref{prop:torus-two-point-bound} completes the proof of Theorem~\ref{thm:torus-two-point}.
Proposition~\ref{prop:torus-two-point-bound} is the analogue of Proposition~\ref{prop:HL-two-point-bound}.
However, since $\bT_n^d$ does not have rotational symmetry, the proof would have to be done a little differently.
We will now focus on proving Proposition~\ref{prop:torus-two-point-bound}.
For notational simplicity and to avoid having to deal with parity issues, we will assume that $n$ is an odd integer throughout the rest of Section~\ref{sec:torus-proof-two-point}.
Simple modifications of these arguments will yield the claims for even $n$.

We set up some notation in preparation for the proof of Proposition~\ref{prop:torus-two-point-bound}.
Recall the construction of $\bT_n^d$ described at the beginning of Section~\ref{sec:results-torus}, and for $1\leq i\leq d_n$ let $B_i$ denote the image of $[-i , i]^d$ in $\bT_n^d$ (thus, $B_{d_n} = \bT_n^d$).
For $n>2^{10}$, let 
$
k_n := 1 + \max\big\{ k\geq 1\, :\, 2^k \leq n/2^5 \big\} 
$.
Define $\Lambda_0 = \{ 0 \}$, $\Lambda_{k_n} = \bT_n^d$,
\begin{align}\label{eqn:656}
\Lambda_j = B_{2^j} \, , \ 1\leq j\leq k_n-1\, ,\ \text{ and }\ \
\annlarge_{\, j} = \Lambda_j\setminus\Lambda_{j-1} \, , \ 1\leq j\leq k_n\, .
\end{align}

Next, we introduce a family of branching processes.
Let
$
\fp_{r}^{\sss (n), -}
:=
1 - \exp\big( -\zeta_n^{-1}\cdot\rho^-(r) \big) 
$,
$r>0$, where $\rho^-$ was as introduced right above \eqref{eqn:torus-def-fp}.
For each $n> 2^{10}$, let $X_i^{\sss (n)}$, $1\leq i\leq d_n$, be independent random variables such that 
$X_i^{\sss (n)}
\sim
\text{Binomial}\big( (2i+1)^d - (2i-1)^d,\ \fp_i^{\sss (n), -}\big)$.
For any $n> 2^{10}$ and $1\leq j\leq d_n$, let $\cT_j^{\sss (n)}$ be a branching process tree where the number of children of each vertex has the same distribution as $\sum_{i=1}^j X_i^{\sss (n)}$.
An argument similar to the one leading to \eqref{eqn:coupling} will show that there exist couplings in which
\begin{gather}
| \cC(x;\, G_{\sss B_i}^-) | \leq |\cT_{2i}^{\sss (n)}|
\ \ \text{ for } x\in B_i\, ,\ 1\leq i\leq 2^{k_n - 1} ,\ \ \text{ and}
\label{eqn:coupling-torus-1}
\\
| \cC^-(x) | = | \cC(x;\, G_n^-) | \leq |\cT_{d_n}^{\sss (n)}|
\ \ \text{ for } x\in\bT_n^d\, ,
\label{eqn:coupling-torus-2}
\end{gather}
where \eqref{eqn:coupling-torus-1} uses the fact that $B_i \subseteq x + B_{2i}$ for any $x\in B_i$ and $1\leq i\leq 2^{k_n - 1} $.
For $1\leq j\leq d_n$, we let $m_j^{\sss (n)} := \sum_{i=1}^j \bE(  X_i^{\sss (n)} )$ be the mean number of offsprings per vertex in $\cT_j^{\sss (n)}$.
The following result is the analogue of Lemma~\ref{lem:HL-zeta-n-m-j-n-asymptotics} for the torus; we omit its proof.

\begin{lem}\label{lem:torus-zeta-n-m-j-n-asymptotics}
Fix $\alpha\in (0, d)$ and $\theta\in (\alpha, \alpha + d/3)$.
Then the following hold:
\begin{enumeratea}
\item 
There exist constants $C, C'>0$ depending only on $\pmtr^{\ast}$ such that
$C n^{(d-\alpha)} \leq \zeta_n \leq C' n^{(d-\alpha)}$ for all $n\geq 3$.
\vskip1pt
\item 
We have,
$1-m_{\sss d_n}^{\sss(n)} 
= 
\big(1 + O(n^{-d})\big)\cdot \zeta_n^{-1} n^{(d-\theta)}$ as $n\to\infty$.
In particular, $1-m_{\sss d_n}^{\sss(n)} \asymp n^{-(\theta - \alpha)}$.
\vskip1pt
\item 
For all $n>2^{10}$,
\begin{align}\label{eqn:135}
m_j^{\sss(n)} 
\leq 
C_{\ref{eqn:135}} \big( j/n \big)^{ (d-\alpha)}
\ \ \text{ for }\ \ 1\leq j\leq 2^{k_n}.
\end{align}
\item
There exist $n_0> 2^{10}$ depending only on $\pmtr$ and $\delta\in(0, 1)$ depending only on $\pmtr^{\ast}$ such that $m_j^{\sss (n)}\leq 1-\delta$ for all $n\geq n_0$ and $1\leq j\leq 2^{k_n }$.
\end{enumeratea}
\end{lem}

We will first complete the proof of Proposition~\ref{prop:torus-two-point-bound} assuming the next lemma, and then we will come back to its proof.

\begin{lem}\label{lem:torus-connect-inside}
Fix $\alpha\in (0, d)$ and $\theta\in (\alpha, \alpha + d/3)$.
Then 
\[
\pr\big( 0\connects z\ \text{ in }\ \Lambda_j \big)
\lesssim
\zeta_n^{-1} \| z \|^{-\alpha}
\ \text{ for }\ z\in\annlarge_{\, j} \, ,\ 1\leq j\leq k_n - 1\, .
\]
\end{lem}

\noindent{\bf Proof of Proposition~\ref{prop:torus-two-point-bound}:}
Consider $x\in\annlarge_{\, i}$, where $3\leq i\leq k_n$. 
If $0 \connects x$, let
$
I_x := \min\big\{ j\in\{ i, i+1, \ldots, k_n \} \, :\, 0\connects x\text{ in }\Lambda_j   \big\}
$, 
and let $I_x = \infty$ otherwise.
Then, for $i=j=k_n$ and for $3\leq i < j\leq k_n$, 
\begin{align}\label{eqn:11-a}
&
\pr( I_x = j )
\leq 
\sum_{y\in\annsmall_{\, j}} 
\sum_{z\in\Lambda_{j-1}} 
\pr\big(\{0\connects z\text{ in }\Lambda_{j-1}\} \disjt \{ z\edge y \} \disjt \{ x\connects y\text{ in }\Lambda_j \}
\big)
\notag\\
&
\hskip20pt
\leq
\sum_{y\in\annsmall_{\, j}} \,
\sum_{z\in\Lambda_{j-1}\setminus ( y + \Lambda_{ j - 3 } ) }  
\pr\big(0\connects z\text{ in }\Lambda_{j-1} \big) 
\pr( z\edge y )
\pr\big( x\connects y\text{ in }\Lambda_j \big)
\notag\\
&
\hskip40pt
+
\sum_{y\in\annsmall_{\, j}} \,
\sum_{z\in \Lambda_{j-1}\cap ( y + \Lambda_{ j - 3 } ) }  
\pr\big(0\connects z\text{ in }\Lambda_{j-1} \big) 
\pr( z\edge y )
\pr\big( x\connects y\text{ in }\Lambda_j \big)
=: 
F_1 + F_2\, ,
\end{align}
where the second step uses the BK inequality.
For $i=j=k_n$ and for $3\leq i < j\leq k_n$,
\begin{align}\label{eqn:11-b}
&
F_1
\lesssim
\zeta_n^{-1} 2^{-(j-3)\alpha}
\sum_{y\in\Lambda_j} 
\pr\big( x\connects y\text{ in }\Lambda_j \big)
\sum_{z\in\Lambda_{j-1} }  
\pr\big(0\connects z\text{ in }\Lambda_{j-1} \big) 
\notag\\
&
\hskip30pt
\leq
\zeta_n^{-1} 2^{-(j-3)\alpha}
\sum_{y\in\Lambda_j} 
\pr\big( x\connects y\text{ in }\Lambda_j \big)
\cdot
\bE |\cT_{2^j}^{\sss (n)}|
\notag\\
&
\hskip60pt
\leq
\zeta_n^{-1} 2^{-(j-3)\alpha}
\sum_{y\in\Lambda_j} 
\pr\big( x\connects y\text{ in }\Lambda_j \big)
\cdot
\big( 1 - m_{2^{k_n} }^{\sss (n)}  \big)^{-1}
\notag\\
&
\hskip90pt
\lesssim
\zeta_n^{-1} 2^{-j\alpha}
\sum_{y\in\Lambda_j} 
\pr\big( x\connects y\text{ in }\Lambda_j \big)
\, ,
\end{align}
where the first step uses \eqref{eqn:1a},
the second step uses \eqref{eqn:coupling-torus-1},
and the last step uses Lemma~\ref{lem:torus-zeta-n-m-j-n-asymptotics}(d).
Now, if $y\in\annlarge_{\, j}$, $4\leq j\leq k_n$, then 
$\Lambda_{j-1} \cap ( y + \Lambda_{j-3} ) \subseteq \annlarge_{\, j - 1}$.
Hence, for $i=j=k_n$ and for $3\leq i < j\leq k_n$,
\begin{align}\label{eqn:11-c}
F_2
&
\lesssim
\zeta_n^{-1} 2^{-j\alpha}
\sum_{y\in\Lambda_j} 
\pr\big( x\connects y\text{ in }\Lambda_j \big)
\sum_{z\in ( y + \Lambda_{ j - 3 } ) }  
\pr( z\edge y )
\notag\\
&
=
\zeta_n^{-1} 2^{-j\alpha}
\sum_{y\in\Lambda_j} 
\pr\big( x\connects y\text{ in }\Lambda_j \big)
\cdot
m_{2^{j-3}}^{\sss (n)}
\leq
\zeta_n^{-1} 2^{-j\alpha}
\sum_{y\in\Lambda_j} 
\pr\big( x\connects y\text{ in }\Lambda_j \big)
\, ,
\end{align}
where the first step uses Lemma~\ref{lem:torus-connect-inside}, and the last step uses Lemma~\ref{lem:torus-zeta-n-m-j-n-asymptotics}(d).

To get an upper bound on 
$
\sum_{y\in\Lambda_j} 
\pr\big( x\connects y\text{ in }\Lambda_j \big)
$, 
note that for $x\in\annlarge_{\, i}\,$, $3\leq i < j \leq k_n - 1$, 
$
\sum_{y\in\Lambda_j} 
\pr\big( x\connects y\text{ in }\Lambda_j \big)
=
\bE |\cC(x ;\, G_{\sss \Lambda_j}^- )|
\leq 
\bE |\cT_{2^{k_n}}^{\sss (n)}|
\lesssim
1
$,
where the second step uses \eqref{eqn:coupling-torus-1}, and the last step uses Lemma~\ref{lem:torus-zeta-n-m-j-n-asymptotics}(d).
Similarly, using \eqref{eqn:coupling-torus-2} and Lemma~\ref{lem:torus-zeta-n-m-j-n-asymptotics}(b) to handle the case $j = k_n$, we get
$
\sum_{y\in\bT_n^d} 
\pr\big( x\connects y\big)
\lesssim
n^{\theta - \alpha}
$
for $x\in\annlarge_{\, i}$, $3\leq i\leq k_n$.
Combining these facts with \eqref{eqn:11-a}, \eqref{eqn:11-b}, and \eqref{eqn:11-c} yields, for $x\in\annlarge_{\, i}$, $3\leq i\leq k_n - 1$,
\begin{align*}
\pr(0\connects x)
&
=
\pr(I_x = i) + \sum_{j=i+1}^{k_n-1} \pr(I_x = j) + \pr(I_x = k_n)
\\
&
\lesssim
\zeta_n^{-1} 2^{-i\alpha}
+ 
\sum_{j=i+1}^{k_n-1} \zeta_n^{-1} 2^{-j\alpha}
+
\zeta_n^{-1} 2^{-k_n\alpha}\cdot n^{\theta - \alpha}
\lesssim
\zeta_n^{-1} 2^{-i\alpha}
+ 
n^{-(d + \alpha - \theta)}\, ,
\end{align*}
where the second step uses Lemma~\ref{lem:torus-connect-inside} to bound $\pr(I_x = i)$, and the last step uses Lemma~\ref{lem:torus-zeta-n-m-j-n-asymptotics}(a) and the definition of $k_n$ given right above \eqref{eqn:656}.
It is easy to check that the resulting upper bound is equivalent to the bounds claimed in \eqref{eqn:788-torus} and \eqref{eqn:789-torus}.
Finally, for $x\in\annlarge_{\, k_n}$,
\begin{align*}
\pr(0\connects x)
=
\pr(I_x = k_n)
\lesssim
\zeta_n^{-1} 2^{-k_n\alpha} \cdot n^{\theta - \alpha}
\lesssim
n^{-(d + \alpha - \theta)}\, ,
\end{align*}
which completes the proof.
\qed

\vskip5pt

\noindent{\bf Proof of Lemma~\ref{lem:torus-connect-inside}:}
Let $j_0\in\bZ_{\geq 3}$ be such that
\begin{align}\label{eqn:12-a}
\frac{2^{\alpha} C_{\ref{eqn:135}}}{ \delta\cdot 2^{ j_0 (d - \alpha) } }
\leq
\frac{1}{2}\, ,
\end{align}
where $\delta$ is as in Lemma~\ref{lem:torus-zeta-n-m-j-n-asymptotics}(d).
It is easy to see that there exists $C_{\ref{eqn:12-b}}(j_0)>0$ depending only on $\pmtr^{\ast}$ such that for all $n > 2^{100 j_0}$, 
\[
\pr\big( 0\connects z\ \text{ in }\ \Lambda_j \big)
\leq
C_{\ref{eqn:12-b}}(j_0) \cdot \zeta_n^{-1} \| z \|^{-\alpha}
\ \text{ for }\ z\in\annlarge_{\, j}\, , \ 1\leq j\leq j_0\, .
\]
Let $n_0$ be as in Lemma~\ref{lem:torus-zeta-n-m-j-n-asymptotics}(d).
We will prove by induction on $j$ that for each $j\geq j_0 + 1$, there exists $C_{\ref{eqn:12-b}}(j)>0$ depending only on $\pmtr^{\ast}$ and $j$ such that 
(i) for each $n > \max\{ n_0,\, 2^{100 j_0} \}$ and each $j_0 + 1\leq j\leq k_n - 1$, 
\begin{align}\label{eqn:12-b}
\pr\big( 0\connects z\ \text{ in }\ \Lambda_j \big)
\leq
C_{\ref{eqn:12-b}}(j) \cdot \zeta_n^{-1} \| z \|^{-\alpha}
\ \text{ for all }\ z\in\annlarge_{\, j}\, ,
\end{align}
and further, 
(ii) $\sup_{j\geq j_0} C_{\ref{eqn:12-b}}(j) < \infty$.
This induction technique is similar to one used in the proof of \cite[Proposition~2.1]{hutchcroft-two-point-torus}.

Fix $\ell\geq j_0 + 1$, and assume that \eqref{eqn:12-b} holds for $j = \ell - 1$ 
for all $n > \max\{ n_0,\, 2^{100 j_0} \}$ such that $k_n \geq \ell + 1$.
Similar to \eqref{eqn:11-a}, we have, for $n > \max\{ n_0,\, 2^{100 j_0} \}$ and $z\in\annlarge_{\, \ell}$,
\begin{align}\label{eqn:12-c}
&
\pr( 0 \connects z\ \text{ in }\ \Lambda_{\ell} )
\leq
\sum_{y\in\annsmall_{\, \ell}} \,
\sum_{w\in\Lambda_{\ell-1}\setminus ( y + \Lambda_{ \ell - j_0 } ) }  
\pr\big(0\connects w\text{ in }\Lambda_{\ell - 1} \big) 
\pr( w\edge y )
\pr\big( z\connects y\text{ in }\Lambda_{\ell} \big)
\notag\\
&
\hskip25pt
+
\sum_{y\in\annsmall_{\, \ell}} \,
\sum_{w\in \Lambda_{\ell - 1}\cap ( y + \Lambda_{ \ell - j_0 } ) }  
\pr\big(0\connects w\text{ in }\Lambda_{\ell - 1} \big) 
\pr( w\edge y )
\pr\big( z\connects y\text{ in }\Lambda_{\ell} \big)
=: 
F_1 + F_2\, .
\end{align}
Now,
\begin{align}\label{eqn:12-d}
F_1
&
\leq
A_2 \zeta_n^{-1} 2^{-\alpha ( \ell - j_0 ) } 
\sum_{w\in\Lambda_{\ell-1}}  
\pr\big(0\connects w\text{ in }\Lambda_{\ell - 1} \big) 
\cdot
\sum_{y\in\Lambda_{\ell}} 
\pr\big( z\connects y\text{ in }\Lambda_{\ell} \big)
\notag\\
&
\leq
2^{j_0 \alpha} A_2 \zeta_n^{-1} \| z \|^{-\alpha}
\cdot\bE| \cC(0;\, G_{\sss\Lambda_{\ell - 1}}^- ) |
\cdot\bE| \cC(z;\, G_{\sss\Lambda_{\ell} }^-  ) |
\leq
2^{j_0 \alpha} A_2 \zeta_n^{-1} \| z \|^{-\alpha} \cdot\delta^{-2} ,
\end{align}
where the first step uses \eqref{eqn:1a}, the second step uses the relation $\| z\|\leq 2^{\ell}$ for $z\in\annlarge_{\, \ell}$, and the last step uses \eqref{eqn:coupling-torus-1} and Lemma~\ref{lem:torus-zeta-n-m-j-n-asymptotics}(d).
Next, 
since $\Lambda_{\ell - 1}\cap ( y + \Lambda_{ \ell - j_0 } ) \subseteq \annlarge_{\, \ell-1}$ for each $y\in\annlarge_{\, \ell}$,
\begin{align}\label{eqn:12-e}
F_2
&
\leq
C_{\ref{eqn:12-b}}(\ell - 1)\cdot \zeta_n^{-1} 2^{-\alpha(\ell - 1)}
\sum_{y\in\annsmall_{\, \ell}} 
\pr\big( z\connects y\text{ in }\Lambda_{\ell} \big)
\sum_{w\in ( y + \Lambda_{ \ell - j_0 } ) }  
\pr( w\edge y )
\notag\\
&
\leq
C_{\ref{eqn:12-b}}(\ell - 1)\cdot \zeta_n^{-1} \| z\|^{-\alpha} \cdot 2^{\alpha}
\cdot
\bE| \cC(z;\, G_{\sss\Lambda_{\ell} }^- ) |
\cdot
m_{2^{\ell - j_0}}^{\sss (n)}
\notag\\
&
\leq
C_{\ref{eqn:12-b}}(\ell - 1)\cdot \zeta_n^{-1} \| z\|^{-\alpha} \cdot 2^{\alpha}
\cdot
\delta^{-1} 
\cdot
C_{\ref{eqn:135}} \cdot 2^{-j_0 (d - \alpha)}
\leq
\frac{1}{2}\cdot C_{\ref{eqn:12-b}}(\ell - 1) \cdot \zeta_n^{-1} \| z\|^{-\alpha} \, ,
\end{align}
where the first step uses the induction hypothesis that \eqref{eqn:12-b} holds with $j = \ell -1$, the third step uses \eqref{eqn:coupling-torus-1}, Lemma~\ref{lem:torus-zeta-n-m-j-n-asymptotics}(d), 
and Lemma~\ref{lem:torus-zeta-n-m-j-n-asymptotics}(c) together with the relation $2^{\ell}\leq n$,
and the fourth step uses \eqref{eqn:12-a}.
Combining \eqref{eqn:12-c}, \eqref{eqn:12-d}, and \eqref{eqn:12-e}, we see that 
for all $n > \max\{ n_0,\, 2^{100 j_0} \}$ with $k_n \geq \ell + 1$,
\begin{gather}
\pr\big( 0\connects z\ \text{ in }\ \Lambda_{\ell} \big)
\leq
C_{\ref{eqn:12-b}}(\ell) \cdot \zeta_n^{-1} \| z \|^{-\alpha}
\ \text{ for each } z\in\annlarge_{\, \ell}\, ,\ \text{ where}
\notag\\
C_{\ref{eqn:12-b}}(\ell) 
:=
\frac{ 2^{j_0\alpha} A_2 }{ \delta^2 }
+
\frac{1}{2}C_{\ref{eqn:12-b}}(\ell-1) \, .
\label{eqn:recurrence}
\end{gather}
Further, the recurrence relation in \eqref{eqn:recurrence} implies that 
$\sup_{j\geq j_0} C_{\ref{eqn:12-b}}(j) < \infty$, 
which completes the proof.
\qed

\vskip5pt

Similar to \eqref{eqn:HL-def-Delta}, we define
\begin{align}\label{eqn:torus-def-Delta}
\Delta_n
:=
\hskip-1pt 
\mathop{\sum\sum}\limits_{x\neq y\in\bT_n^d\setminus\{0\} }  
\hskip-1pt 
\pr( 0 \connects x ) \pr( 0 \connects y ) \pr( x \edge y )
\ \text{ and }\ 
\widetilde\Delta_n
:=
\hskip-1pt 
\sum_{w\in\bT_n^d\setminus\{0\}} 
\hskip-1pt 
\pr( 0 \connects w ) \pr( 0 \edge w )\, .
\end{align}
The next result follows by arguing as in the proofs of Lemmas~\ref{lem:2} and \ref{lem:3}, and using the upper bounds in Proposition~\ref{prop:torus-two-point-bound}; we omit the details.

\begin{lem}\label{lem:2-torus}
The result in Lemma~\ref{lem:2} remains true if $\Lambda_n$ is replaced by $\bT_n^d$.
Further, the result in Lemma~\ref{lem:3} remains true if $\Lambda_n$ is replaced by $\bT_n^d$, and $\cT_n^{\sss(n)}$ is replaced by $\cT_{d_n}^{\sss(n)}$ (as defined above \eqref{eqn:coupling-torus-1}).
\end{lem}

Finally, the proof of Proposition~\ref{prop:torus-two-point-lower-bound} is similar to that of Proposition~\ref{prop:HL-two-point-lower-bound} and is omitted.

\subsection{Proof outline: Theorems~\ref{thm:torus-scaling-limit} and \ref{thm:torus-surplus-convergence}}\label{sec:torus-outline-proof}
First, let us note that the result in Proposition~\ref{prop:torus-phase-transition} follows along the same lines as the proof of Proposition~\ref{prop:HL-phase-transition}, so we omit its proof.

Turning to the proof of Theorem~\ref{thm:torus-scaling-limit}, we fix $\alpha\in (0, 5d/6)$ and $\theta$ as in \eqref{eqn:HL-def-theta}.
Let $G_n^-$ have $m_n$ components, and write $\cC_i^-=\cC_i^{\sss (n), -}$ for the $i$-th largest component in $G_n^-$, $1\leq i\leq m_n$.
We will apply Theorem~\ref{thm:gen-2} with
$\vx=(x_1, \ldots, x_{m_n}), q$, and $\vM=(M_1, \ldots, M_{m_n})$, where
\begin{equation}\label{eqn:510-torus}
x_i = n^{-2d/3}|\cC_i^-| \, ,\ \
q = \lambda + \zeta_n^{-1} n^{- \theta + 4d/3}  \, ,
	\ \text{ and }\
	M_i =  \scl\big(1, 1/|\cC_i^-| \big) \cC_i^-\, ,\ \ 1\leq i\leq m_n\, .
\end{equation}
With this notation, we let 
$\sigma_2, \sigma_3, x_{\max}, \diamax$, and $\tau$ be as in 
\eqref{eqn:def-sigma_k-x_max}, \eqref{eqn:tau-dmax-def}.
Then we claim that
\begin{gather}
q - (\sigma_2)^{-1} \weakc \lambda 
\ \text{ as } \ n\to\infty\, ,
\label{eqn:52-a}\\
\sigma_3\cdot (\sigma_2)^{-3} \weakc 1 \, ,\ \
\tau = \zeta_n^2\cdot n^{2\theta - 7d/3} \cdot (1+o_P(1))
\ \ \text{ as } \ n\to\infty\, ,
\label{eqn:52-b}
\end{gather}
and for each $r>0$, there exists $C_{\ref{eqn:52-c}}(r)$ depending only on $\pmtr^{\ast}$ and $r$ such that
\begin{gather}\label{eqn:52-c}
\pr\big(
\diamax \geq C_{\ref{eqn:52-c}}(r) \cdot n^{\theta - \alpha} \log n 
\big)
+
\pr\big(
x_{\max} \geq C_{\ref{eqn:52-c}}(r) \cdot n^{2(\theta - \alpha - d/3)} \log n
\big)
\lesssim n^{-r}\, .
\end{gather}
The proof of \eqref{eqn:52-a} uses Lemma~\ref{lem:2-torus} and follows steps similar to those in the proof of Proposition~\ref{prop:HL-sigma-2-asymptotics}.
The proofs of \eqref{eqn:52-b} and \eqref{eqn:52-c} are similar to those of Propositions~\ref{prop:HL-sigma-3-asymptotics}, \ref{prop:HL-tau-asymptotics}, and \ref{prop:HL-xmax-diamax-bound}.
Now the asymptotics in \eqref{eqn:52-a}, \eqref{eqn:52-b}, and \eqref{eqn:52-c} can be combined and we can argue in a manner similar to Section~\ref{sec:HL-proof-scaling-limit} to complete the proof of Theorem~\ref{thm:torus-scaling-limit}.

Finally, to prove Theorem~\ref{thm:torus-surplus-convergence}, fix $\alpha\in (0, 2d/3)$ and $\theta$ satisfying \eqref{eqn:HL-def-theta}.
Note that we have the following analogue of Lemma~\ref{lem:34}:
\[
\pr\big(\spls\big( \cC^-(0) \big) \neq 0 \big)
\lesssim
\begin{cases}
	n^{ -( d + 2\alpha - 2\theta ) } ,\ \ \text{ if } 0<\alpha\leq d/2\, ,\\
	n^{- (2d - \alpha - \theta)} ,\ \ \text{ if } d/2<\alpha< 2d/3.
\end{cases}
\]
We can use this estimate and imitate the argument presented in Section~\ref{sec:HL-proof-surplus} to complete the proof of Theorem~\ref{thm:torus-surplus-convergence}.

\section*{Acknowledgments}
We thank Shankar Bhamidi and Remco van der Hofstad for insightful discussions.
We also thank Gordon Slade for drawing our attention to the paper \cite{slade-plateau}.
SS was supported in part by the grant ANRF/ARGM/2025/001047/MTR from ANRF, and by the DST FIST program - 2021 [TPN--700661].


\bibliographystyle{imsart-number}
\bibliography{scaling}

\end{document}